\documentclass[12pt,twoside,final,psamsfonts]{amsart}
\usepackage[psamsfonts]{amssymb}
\usepackage{times,a4wide}
\usepackage{pictex}

\theoremstyle{plain}         
\newtheorem*{TheoremA}{Theorem A}
\newtheorem*{TheoremB}{Theorem B}
\newtheorem*{TheoremC}{Theorem C}
\newtheorem{theorem}{Theorem}
\newtheorem{proposition}[theorem]{Proposition}
\newtheorem{corollary}[theorem]{Corollary}
\newtheorem{lemma}[theorem]{Lemma}

\newtheorem*{claim*}{Claim}
\theoremstyle{remark}
\newtheorem{remark}{Remark}
\theoremstyle{definition}
\newtheorem{definition}{Definition}

\newcommand{\Int}{\operatorname{Int}}
\newcommand{\Samp}{\operatorname{Samp}}
\newcommand{\C}{{\mathbb C}}
\newcommand{\N}{{\mathbb N}}
\newcommand{\D}{{\mathbb D}}
\newcommand{\R}{{\mathbb R}}
\newcommand{\X}{{\mathcal X}}
\newcommand{\calf}{\mathcal{F}}

\newcommand{\db}{{\bar\partial}}
\newcommand{\iC}{{\int_{\C}}}

\newcommand{\Fpa}{{\mathcal F_{\phi,\omega}^p}}
\newcommand{\Fdosa}{{\mathcal F_{\phi,\omega}^2}}
\newcommand{\Finfa}{{\mathcal F_{\phi,\omega}^\infty}}
\newcommand{\ellpa}{{\ell_{\phi,\omega}^p}}

\newcommand{\tinu}{{\tilde \nu}}
\newcommand{\ellinfa}{{\ell_{\phi,\omega}^\infty}}

\renewcommand{\Re}{\mathop{Re}}

\title{Interpolating and sampling sequences for entire functions}
\author{Nicolas Marco}
\address{Dept. Matem\`atiques, Universitat Autonoma de Barcelona, 08193
Bellaterra (Barcelona), Spain}
\email{nmarco@mat.uab.es}
\thanks{The first author is supported by the Research Training Network
``Analysis and Operators'', with contract number HPRN-CT-2000-00116. 
The last two authors are supported by the DGICYT grant
 PB98-1242-C02-01 and by the CIRIT grant 1998SGR00052}
\author{Xavier Massaneda}
\address{Dept.\ Matem\`atica Aplicada i An\`alisi,
Universitat  de Barcelona, Gran Via 585, 08071 Bar\-ce\-lo\-na, Spain}
\email{xavier@mat.ub.es}
\author{Joaquim Ortega-Cerd\`a}
\address{Dept.\ Matem\`atica Aplicada i An\`alisi,
Universitat  de Barcelona, Gran Via 585, 08071 Bar\-ce\-lo\-na, Spain}
\email{quim@mat.ub.es}

\date{May 22, 2002}

\begin{document}

\begin{abstract}
We characterise interpolating and sampling sequences for the spaces of
entire functions $f$ such that $f e^{-\phi}\in L^p(\C)$, $p\geq 1$ (and some
related weighted classes), where $\phi$ is a subharmonic weight whose Laplacian
is a doubling measure. The results are expressed in terms of some
densities adapted to the metric induced by $\Delta\phi$. They generalise
previous results by Seip for the case $\phi(z)=|z|^2$, and by Berndtsson \&
Ortega-Cerd\`a and Ortega-Cerd\`a \& Seip for the case when $\Delta\phi$ is
bounded above and below.
\end{abstract}

\maketitle

\tableofcontents

\section{Introduction} 
In this paper we provide Beurling-type density conditions for  sampling and
interpolation in certain generalised Fock spaces. We consider a rather general
situation, with only mild regularity  conditions on the possible growth. Let
$\phi$ be a (nonharmonic) subharmonic function whose Laplacian $\Delta\phi$ is
a doubling measure (see definition and properties in
Section~\ref{doubling-measures}), and let $\omega$ denote a flat weight, that
is, a positive measurable function with slow growth (see details in
Section~\ref{flat-weights}). The spaces we deal  with are parametrised by an
index $p\in [1,\infty]$, as follows:
\begin{eqnarray*}
\Fpa&=&\bigl\{f\in H(\C) :
\|f\|_{\Fpa}^p=\int_{\C}|f|^p e^{-p\phi}\omega^p \rho^{-2}<\infty\bigr\}\qquad
1\leq p<\infty ,\\
\Finfa&=&\bigl\{f\in H(\C) : \|f\|_{\Finfa}=
\sup_{z\in\C}\omega(z)|f(z)|e^{-\phi(z)}<\infty\bigr\} .
\end{eqnarray*}
The function $\rho^{-2}$ is a regularised version of $\Delta\phi$, as described
in \cite{Ch}. More
precisely, if $\mu=\Delta\phi$ and  $z\in\C$, then $\rho_\phi(z)$ (or simply
$\rho(z)$ if no confusion can arise) denotes the positive radius such that 
$\mu(D(z,\rho(z))=1$. Such a radius  exists because doubling measures have no
mass on circles. 

Canonical examples of the weights considered are
$\phi(z)=|z|^\beta$, with  $\beta>0$, and $\omega=\rho^\alpha$, $\alpha\in\R$.

Two particular families of spaces seem of special interest. The first  one are
the usual weighted $L^p$-spaces of entire functions, obtained with
$\omega=\rho^{2/p}$. The second case arises when $\omega=1$; then the spaces
$\Fpa$ coincide with 
\[
\{f\in H(\C) : 
\int_{\C} |f|^p e^{-p\phi}\Delta\phi <\infty\} .
\]

Since functions $f$ in the spaces $\Fpa$ are determined by the growth of $|f|$,
their restriction to a sequence should be described as well in terms of growth.

Let $\Lambda\subset \mathbb C$ be a sequence and let
$v=\{v_\lambda\}_{\lambda\in\Lambda}$ be an associated sequence of values.

\begin{definition} A sequence $\Lambda$ is an \emph{interpolating
sequence for} $\Fpa$, $1\le p<\infty$ (denoted 
$\Lambda\in\Int \Fpa$),
if for every sequence of values 
$v$ such that
\begin{equation*} 
\|v\|_{\ellpa(\Lambda)}^p
=\sum_{\lambda\in\Lambda}   \omega^p(\lambda)  |v_\lambda|^p
e^{-p\phi(\lambda)}<\infty 
\end{equation*}
there exists $f\in \Fpa$ such that $f|\Lambda=v$.

Also, $\Lambda\in\Int \Finfa$ if for 
every sequence of values $v$ such that
\[ 
\|v\|_{\ellinfa(\Lambda)}=\sup_{\lambda\in\Lambda} 
\omega(\lambda) |v_\lambda| 
e^{-\phi(\lambda)}<\infty 
\]
there exists $f\in \Finfa$ such that $f|\Lambda=v$.
\end{definition}

An application of the open mapping theorem shows that when $\Lambda\in\Int\Fpa$
there is  $M>0$ such that for any $v \in \ellpa(\Lambda)$, there exists $f \in
\Fpa$ with $f|\Lambda =v$ and 
\begin{equation}\label{interpolation}
||f||_{\Fpa} \leq M||v||_{ \ellpa(\Lambda)}.
\end{equation}
The least possible $M$ in \eqref{interpolation} is called the 
\emph{interpolating constant} of $\Lambda$ 
and is denoted by $M_{\phi,\omega}^{p}(\Lambda)$, or $M(\Lambda)$ if no
confusion is possible.  

\begin{definition}
A sequence $\Lambda$ is a \emph{sampling
sequence for} $\Fpa$, $1\leq p< \infty$ (denoted 
$\Lambda\in\Samp \Fpa$), if there exists 
$C>0$ such that
for every $f\in \Fpa$ 
\begin{equation}\label{sampling1}
C^{-1}\|f|\Lambda\|_{\ellpa(\Lambda)} 
\leq \| f\|_{\Fpa}\leq 
C \|f|\Lambda\|_{\ellpa(\Lambda)} .
\end{equation}
Also, $\Lambda\in\Samp \Finfa$ if 
there exists $C>0$ such that
for every $f\in \Finfa$ 
\begin{equation}\label{sampling2}
\| f\|_{\Finfa}\leq  C \|f|\Lambda\|_{\ellinfa
(\Lambda)} .
\end{equation}
\end{definition}
The least constant $C$ verifying these inequalities is called the 
\emph{sampling constant} of $\Lambda$ and is denoted
$L_{\phi,\omega}^{p}(\Lambda)$, or simply $L(\Lambda)$.

The definitions of interpolating and sampling sequences in the spaces defined
by $L^\infty$ norms reflect the maximal growth for functions in the space, and 
are natural.  The definition for $p<\infty$ can be motivated in the following
way. Consider for instance the case $p=2$. The estimates of the normalised
Bergman kernel $k_{\phi,\omega}(\lambda,z)$ in $\Fdosa$  (see
Lemma~\ref{bergman-diagonal}) show that   $\langle
k_{\phi,\omega}(\lambda,\cdot),f\rangle \simeq f(\lambda)  \omega(\lambda)
e^{-\phi(\lambda)}$ for all $f\in\Fdosa$. Thus $\Lambda\in \Samp \Fdosa$ if and
only if  
\[
\| f\|_{\Fdosa}\simeq \sum_{\lambda \in\Lambda} |\langle
k_{\phi,\omega}(\lambda,\cdot),f \rangle|^2\quad \text{for all } f\in\Fdosa,
\] 
that is, if
and only if $\{k_{\phi,\omega}(\lambda,\cdot)\}_{\lambda\in \Lambda}$ is a
frame in  $\Fdosa$. Similarly, $\Lambda\in \Int  \Fdosa$ if and only if 
$\{k_{\phi,\omega}(\lambda,\cdot)\}_{\lambda\in \Lambda}$ is a Riesz  basis in
its closed linear span in $\Fdosa$. These are the standard problems of
interpolation and sampling in Hilbert spaces of functions with reproducing
kernels \cite{ShSh}. For $p\ne 2$ the previous definitions give the appropriate
notions of interpolation and sampling as well, in view of the pointwise growth
of functions in the spaces (see Lemma~\ref{estimacio-puntual} and
Remark~\ref{general-bergman}). 

Our description of interpolating and sampling sequences is expressed in terms
of certain Beurling-type densities adapted to the metric induced by
$\Delta\phi$, or more precisely, by its regularisation $\rho^{-2}(z) dz\otimes
d\bar z$. Before introducing the densities we need the notion of
$\rho$-separation.

\begin{definition} A sequence $\Lambda$ is 
$\rho$-\emph{separated} if there exists $\delta>0$ such that
\[
|\lambda-\lambda'|\geq\delta\max(\rho(\lambda),\rho(\lambda))
\qquad \lambda\neq \lambda'  ,
\]
\end{definition}
This is equivalent to saying that the points in $\Lambda$ are separated by a
fixed distance in the metric above (Lemma~\ref{distance}).

\begin{definition}  Assume that $\Lambda$ is a $\rho$-separated sequence and
denote $\mu=\Delta\phi$.

The \emph{upper uniform density of} $\Lambda$ with respect to
$\Delta\phi$ is
\[
\mathcal D^+_{\Delta\phi}(\Lambda)=\limsup_{r\to \infty}\sup_{z\in\C}
\frac{\# \bigl(\Lambda\cap \overline{D(z,r\rho(z))}\bigr)}
{\mu(D(z,r\rho(z)))}.
\]

The \emph{lower uniform density of} $\Lambda$ 
with respect to $\Delta\phi$ is
\[
\mathcal D^-_{\Delta\phi}(\Lambda)=\liminf_{r\to \infty}\inf_{z\in\C}
\frac{\# \bigl(\Lambda\cap \overline{D(z,r\rho(z))}\bigr)}
{\mu(D(z,r\rho(z)))}.
\] 
\end{definition}

The main theorems  are the following. Let $\Omega_\phi$ denote the class of
flat weights.
\begin{TheoremA} 
A sequence $\Lambda$ is sampling for $\Fpa$,  $p\in[1,\infty)$,
$\omega\in\Omega_\phi$, if and only if $\Lambda$ is a finite  union of
$\rho$-separated sequences containing a $\rho$-separated subsequence
$\Lambda'$  such that $\mathcal D^-_{\Delta\phi}(\Lambda')>1/2\pi$. A sequence
$\Lambda$ is sampling for $\Finfa$  if and only if $\Lambda$  contains a
$\rho$-separated subsequence $\Lambda'$  such that $\mathcal
D^-_{\Delta\phi}(\Lambda')>1/2\pi$.
\end{TheoremA}

\begin{TheoremB} A sequence $\Lambda$ is interpolating for $\Fpa$, 
$p\in[1,\infty]$,
$\omega\in\Omega_\phi$, if and only if $\Lambda$ is 
$\rho$-separated and
$\mathcal D^+_{\Delta\phi}(\Lambda)<1/2\pi$.
\end{TheoremB}

In particular, there are no sequences which are simultaneously sampling and
interpolating (it should be mentioned that this is not obtained as a corollary
of the theorems; it is actually an important ingredient of the proofs).

These results generalise previous work, beginning with the papers by Seip and
Seip-Wallst\'en \cite{Se92}, \cite{SeWa}. They described the interpolating and
sampling sequences for the classical Fock space in terms of the so-called
Nyquist densities. In the notation above this corresponds to 
$\phi(z)=|z|^2$ and $\omega\equiv 1$. This was extended in
\cite{LySe}, \cite{BeOr} and \cite{OrSe} to the case of entire functions $f$
such that $fe^{-\phi}\in L^p(\C)$, where $\phi$ is subharmonic  with bounded
Laplacian  $\varepsilon < \Delta\phi < M$. The description was given again in
terms of some Nyquist type densities. In these cases the function $\rho$ is
bounded above and below, hence the metric $\rho^{-2}(z)dz\otimes d\bar z$ is
equivalent to the Euclidean metric. In particular, $\rho(z)$ can be replaced
by the constant $1$ in the definition of the uniform densities. 

There are also some partial results in several complex variables. The classical
Fock space has been studied in \cite{MaPa00} and the weighted scenario in
\cite{Ln}. In this context there exist necessary or sufficient density
conditions, which do not completely characterise the sampling or interpolating
sequences. 

Interpolation problems for other spaces of functions related to these weights
have been considered by Squires and Berenstein and Li (see for instance
\cite{Sq83}, \cite{BeLi95} and the references therein).

The results mentioned above relied on the remarkable work by Beurling
\cite{Br} and on H\"ormander's  weighted $L^2$-estimates for the $\db$ equation
\cite{Ho}.  In our proofs we first extend Beurling's tools to the
context of certain spaces  which are non-invariant under translations. We need
as well  a H\"ormander type theorem giving precise estimates for the $\db$
equation in Banach norms other than $L^2$.

The plan of the paper is the following: In Section~\ref{doblants}  we study the
properties of doubling measures and introduce the flat weights. Recall that the
only assumption on our subharmonic weight $\phi$ is that the measure $\Delta
\phi$ is doubling. We will need a regularisation of $\phi$ and the construction
of a multiplier associated to $\phi$ (that is, an entire function $f$ such that
$|f|$ approximates $e^\phi$), very much in the spirit of  \cite{LyMa} and
\cite{Or}.  

In Section~\ref{basics} we state and prove some basic properties of functions
in $\Fpa$. The main result in this section is the following H\"ormander type
theorem.

\begin{TheoremC}\label{estimacions-dbar} 
Let $\phi$ be a subharmonic function
such that $\Delta\phi$ is a doubling measure. For any $\omega\in \Omega_\phi$, 
there is a solution $u$ to the equation $\db u = f$ such that $\| u
e^{-\phi}\omega\|_{L^p(\C)}  \lesssim \| f e^{-\phi}\omega\rho
\|_{L^p(\C)}$ for any $1\le p\le \infty$.  
\end{TheoremC}

We also include the estimates of the Bergman kernel that justify the notion of
interpolating and sampling sequences we have considered.  Finally , we  study
the invariance of our spaces under some appropriate scaled  translations. This
leads to the notion of weak limit and  the corresponding analysis analogous to
Beurling's. 

Section~\ref{basic-int-samp} is devoted to some preliminary (but important)
properties of interpolating and sampling sequences, including their behaviour
under weak limits. The main results in this section are some inclusion
relations between various spaces of interpolating and sampling sequences, and
the fact that there are no sequences which are simultaneously interpolating
and sampling for the same space of functions $\Fpa$.

In Section~\ref{sufisampling} we prove the sufficiency part of Theorem A.
We use again an approach similar to that of Beurling.

Section~\ref{necsampling} includes the proof of the necessity part of Theorem
A. For this we need once more Beurling's analysis, plus the non-existence of
sampling and interpolating sequences. We use some theorems that relate the
densities of sampling and interpolating sequences, following the ideas by
Ramanathan and Steger \cite{RaSt}.

Section~\ref{neceinterpolation} is devoted to the proof of the necessity part
of Theorem B. We use Ramanathan and Steger's theorem plus an original
argument that shows that the density inequality is actually strict. 

Finally, in Section~\ref{sufinterpolation} we deal with the sufficiency part of
Theorem B.  In the course of the proof, whose main tool is the multiplier, we
need to express the density in terms of rectangles instead of disks. The usual
argument of Landau~\cite{La} does not work, in view of the inhomogeneity of our
measures. Theorem~\ref{density-rectangles} takes care of this. 

A final word on notation: $C$ denotes a finite constant that may change in
value from one occurrence to the next.The expression  $f \lesssim g$ means that
there is a constant $C$ independent of the relevant variables such that $f\le C
g$, and $f\simeq g$ means that $f\lesssim g$ and $g\lesssim f$. 

\section{Subharmonic functions with  doubling Laplacian}\label{doblants}

In this chapter we recap some results on doubling measures and subharmonic
functions $\phi$ whose Laplacian $\Delta\phi$ is doubling. We start with 
regularity and integrability conditions on doubling measures. Next we show that
$\phi$ can be regularised, in the sense that there exists $\psi$ subharmonic
and regular for which the interpolation and sampling problems for $\Fpa$ and
$\mathcal F_{\psi,\omega}^p$ are equivalent. The final part is dedicated to the
construction of the multiplier. A useful application of this is the
existence of holomorphic ``peak functions'' with controlled growth.

\begin{definition}
A nonnegative Borel measure $\mu$ in $\C$ is called \emph{doubling} if  there
exists $C>0$ such that 
\[
\mu(D(z,2r))\leq C \mu(D(z,r))
\] 
for all $z\in\C$ and $r>0$. We denote by $C_\mu$ the minimum constant 
$C$ for which the inequality holds.
\end{definition}

Recall that when $\phi$ is subharmonic $\Delta\phi$ is a nonnegative Borel
measure, finite on compact sets.

For convenience we write $D^r(z)=D(z,r\rho(z))$ and
$D(z)=D^1(z)$. We will write $D_\phi^r(z)$ when we need to
stress that the radius depends on $\phi$.

Henceforth $dm$ denotes the Lebesgue measure in $\C$. We also use the
measure $d\sigma=dm/\rho^2$, which should be thought of as a 
doubling regularisation of
$\Delta\phi$ (see Theorem~\ref{regularitzacio}).

\subsection{Doubling measures}\label{doubling-measures}
Throughout this section we assume that $\mu$ is a positive doubling
measure non-identically zero.
We begin with a result of Christ \cite[Lemma 2.1]{Ch}.

\begin{lemma}\label{discs}
Let $\mu$ be a doubling measure in $\mathbb C$. There exists
$\gamma>0$ such that for any disks $D, D'$ of respective
radius $r(D)>r(D')$ with $D\cap D'\neq\emptyset$:
\[
\left(\frac{\mu(D)}{\mu(D')}\right)^\gamma\lesssim\frac{r(D)}{r(D')}
\lesssim \left(\frac{\mu(D)}{\mu(D')}\right)^{1/\gamma} .
\]
\end{lemma}
In particular, the support of $\mu$ has positive Hausdorff dimension.

\begin{remark}
This implies that 
there exist $k,\varepsilon>0$ such that
\begin{equation}\label{cota-inferior}
r^\varepsilon\lesssim \mu(D^r(z))\lesssim r^k \qquad z\in\C\ ,\ r>1  .
\end{equation}
Also, applying Lemma \ref{discs} and (\ref{cota-inferior}) 
to $D(0,|z|)$ and $D(z)$ we have, for $\rho(z)\leq |z|$
\[
\frac 1{|z|^{k/\gamma}}\lesssim
\bigl(\frac
1{\mu(D(0,|z|))} \bigr)^{1 / \gamma }\lesssim
\frac{\rho(z)}{|z|}\lesssim \bigl(\frac
1{\mu(D(0,|z|))} \bigr)^{ \gamma }\lesssim\frac 1{|z|^{\varepsilon\gamma}} .
\]
On the other hand, if $|z|<\rho(z)$, then $0\in D(z)$. 
Thus Lemma~\ref{discs} implies $\rho(z)\simeq \rho(0)$, hence $|z|<C$.
Therefore, there exist $\eta,C_0>0$ and $\beta\in (0,1)$ such that
\begin{equation}\label{cota-rho}
C_0^{-1} |z|^{-\eta}\le  \rho(z)\le C_0 |z|^{\beta} \qquad\qquad |z|>1 .
\end{equation}
\end{remark}

Let us study in more detail the relationship between $\rho(z)$ and 
$\rho(\zeta)$ for various $z,\zeta\in\C$.
A first observation is that $\rho(z)$ is a Lipschitz function. 
More precisely
\begin{equation}\label{lipschitz}
|\rho(z)-\rho(\zeta)|\leq |z-\zeta|\qquad z,\zeta\in\C .
\end{equation}
To see this there is no loss of generality in assuming that
$z,\zeta\in\R$, $\zeta<z$ and $\rho(\zeta)<\rho(z)$. Then 
$\zeta-\rho(\zeta)<z-\rho(z)$, since otherwise
$D(\zeta)\subset D(z)$, contradicting the fact that 
$\mu(D(z))=\mu(D(\zeta))=1$.

\begin{lemma}\label{christ} \cite[p.205]{Ch}. 
If $\zeta\notin D(z)$
then
\begin{equation*}
\frac{\rho(\zeta)}{\rho(z)}\lesssim
\left(\frac{|z-\zeta|}{\rho(\zeta)}\right)^{1-\delta}
\end{equation*}
for some $\delta\in (0,1)$ depending only on the doubling constant $C_\mu$.
\end{lemma}

As a consequence of Lemma~\ref{discs} and \eqref{cota-rho} 
we have
\begin{corollary}\label{rho}
For every $r>1$ there exists $\gamma>0$ such that if $\zeta\in
B(z,r)$ then
\begin{equation*}
\frac{1}{r^\gamma}\lesssim \frac{\rho(z)}{\rho(\zeta)} \lesssim r^\gamma.
\end{equation*}
\end{corollary}

It will be convenient to express some of the results in terms of the distance
$d_\phi$ induced by the metric $\rho^{-2}(z) dz\otimes d\bar z$.

\begin{lemma}\label{distance} There exists $\delta\in (0,1)$ such that 
for every $r>0$ 
there exists $C_r>0$ such that 
\begin{itemize}
\item[(a)] $\displaystyle C_r^{-1} \frac{|z-\zeta|}{\rho(z)}
\leq d_\phi(z,\zeta)\leq
C_r \frac{|z-\zeta|}{\rho(z)}\ $ if $\ |z-\zeta|\leq r\rho(z)$.
\item[(b)] $\displaystyle C_r^{-1} 
\left(\frac{|z-\zeta|}{\rho(z)}\right)^\delta
\leq d_\phi(z,\zeta)\leq
C_r \left(\frac{|z-\zeta|}{\rho(z)}\right)^{2-\delta}\ $ if 
$\ |z-\zeta|> r\rho(z)$.
\end{itemize}
\end{lemma}

This shows, in particular, that a sequence $\Lambda$ is $\rho$-separated if and
only if there exists $\delta>0$ such that  $\inf_{\lambda\neq\lambda'}
d_\phi(\lambda,\lambda')>\delta$.

\begin{proof}
By definition
\[
d_\phi(z,\zeta)=\inf \int_0^1 |\gamma'(t)| \rho^{-1}(\gamma(t)) dt ,
\]
where the infimum is taken over all piecewise $\mathcal C^1$ curves 
$\gamma:[0,1]\rightarrow\C$ with $\gamma(0)=z$ and $\gamma(1)=\zeta$.

The lower inequalities are contained in \cite[Lemma 3.1]{Ch} and its proof. 

The upper estimate in case (a) is immediate from Corollary~\ref{rho}. In case 
(b) take $\gamma(t)=z+t(\zeta-z)$ and use Lemma~\ref{christ};  then
\[
d_\phi(z,\zeta)\leq |\zeta-z|\int_0^1 \frac{dt}{\rho(\gamma(t))}
\lesssim\int_0^1\frac{ (t|\zeta-z|)^{1-\delta}}{(\rho(z))^{2-\delta}} dt
\lesssim\left(\frac{|z-\zeta|}{\rho(z)}\right)^{2-\delta} .
\]
\end{proof} 

From now on, given $z\in\C$ and $r>0$, we denote
\[
B(z,r)=\{\zeta\in\C : d_\phi(z,\zeta)<r\}.
\]

Doubling measures satisfy certain integrability conditions.

\begin{lemma}\label{integrals} Let $\mu$ be a doubling measure.
There exist $C>0$ and $m\in\N$ depending on $C_\mu$ such that
for any $r>0$
\begin{itemize}
\item[(a)] $\displaystyle \int_{D(z,r)}\log\bigl( \frac{2r}{|z-\zeta|}\bigr)
d\mu(\zeta)\leq C\; \mu(D(z,r))\quad\quad z\in\C $.
\item[(b)] $\displaystyle \sup\limits_{z\in\C}\int_{\C}
\frac{d\mu(\zeta)}{1+d_\phi^m(z,\zeta)} <\infty\ $.
\end{itemize}
\end{lemma}

\begin{proof} (a) is \cite[Lemma 2.3]{Ch}.

(b) According to Lemma~\ref{distance} it 
is enough to consider the integral on $|z-\zeta|\geq r \rho(z)$.
Applying Fubini's theorem we see that
\begin{multline*}
\int_{\zeta\notin 
D^r(z)}\bigl(\frac{\rho(z)}{|z-\zeta|}\bigr)^m d\mu(\zeta)=
\int_{\zeta\notin D^r(z)}
m\int_0^{\rho(z)/|z-\zeta|} 
t^{m-1}dt d\mu(\zeta) \\
= m\int_0^{1/r}
t^{m-1}\int\limits_{t<\rho(z)/|z-\zeta|<1/r}
d\mu(\zeta) dt
\leq m \int_0^{1/r} t^{m-1} \mu(D^{1/rt}(z) )\; dt .
\end{multline*}

Let $x_0=\log_2 C_\mu$, and for a given 
$t$ denote
$k(t)=\inf\{k\in\N : 1/rt\leq 2^k\}$. Then
\[
\mu(D^{1/rt}(z))\leq \mu(D^{2^{k(t)}}(z))\leq 2^{x_0 
k(t)}\leq
\bigl(\frac 2{rt}\bigr)^{x_0} ,
\]
hence the integral is bounded if $m>x_0$.

This and Lemma~\ref{distance}(b) show that the result holds for 
$m$ big enough.
\end{proof}

\begin{remark}\label{int-to-zero}
It is clear from the proof that 
\begin{itemize}
\item[(b')] $\displaystyle \lim\limits_{r\to\infty}\sup
\limits_{z\in\C}\int_{\zeta\notin B(z,r)}
\frac{d\mu(\zeta)}{d_\phi^m(z,\zeta)} =0$.
\end{itemize}
\end{remark}

There is a discrete version of the previous Lemma. 

\begin{lemma}\label{suma}
Let $\Lambda$ be a $\rho$-separated sequence. There exists
$m\in\N$ such that  
\[
\sup_{z\in\C}\sum\limits_{\lambda\in \Lambda}
\frac{1}{1+d_\phi^m(z,\lambda)}<\infty .
\]
\end{lemma}

\begin{proof}
By the separation and Lemma~\ref{distance}, it is enough to see that
for $m$ big enough
\[
\sup_{z\in\C}\sum\limits_{\lambda\notin B(z,r)}
\bigl(\frac{\rho(\lambda)}{|z-\lambda|}\bigr)^m<\infty .
\]
Take
$\delta>0$ such that the balls
$\{B(\lambda,\delta)\}_{\lambda\in\Lambda}$ are pairwise disjoint. By
 Corollary~\ref{rho}
\begin{eqnarray*}
\sum_{\lambda\notin B(z,r))}
\bigl(\frac{\rho(\lambda)}{|z-\lambda|}\bigr)^m \lesssim 
\sum_{\lambda\notin B(z,r)}
\int_{B(\lambda,\delta)} 
\bigl(\frac{\rho(\zeta)}{|z-\zeta|}\bigr)^m d\mu(\zeta) 
\lesssim 
\int\limits_{\lambda\notin B(z,r)}  \bigl(\frac{\rho(\zeta)}{|z-\zeta|}
\bigr)^m\; d\mu (\zeta) .
\end{eqnarray*} 
Lemma \ref{integrals}(b) implies that the integral is bounded.
\end{proof}

For later use, we state a refinement that follows similarly from Remark
\ref{int-to-zero}.

\begin{corollary}\label{suma-zero}
Let $\Lambda$ be a $\rho$-separated sequence. There exists
$m\in\N$ such that
\[
\lim_{r\to\infty}\sup_{z\in\C}\sum\limits_{\lambda\notin B(z,r)}
\frac{1}{d_\phi^m(z,\lambda)}=0 .
\]
\end{corollary}

We will need to partition the plane in rectangles of constant mass. We do that 
by adapting a general result of \cite{Yu} to our setting (see
also  \cite[Theorem 2.1]{Dr}).

\begin{theorem}\label{partition}
Let $\mu$ be a positive doubling measure non-identically zero.
There exists a ``partition'' of
$\C$ in  rectangles $R_k$ with sides parallel to the coordinate axis 
such that:
\begin{itemize}
\item[(a)] $\mu=\sum_k\mu_k$, where $\mu_k := \mu_{|R_k}$
satisfy $\mu_k(\C)=1$.
\item[(b)] $R_k$ are quasi-squares: there exists
$e>1$ depending only on $C_\mu$ such that
the ratio of sides of each $R_k$ lies in the interval $[1/e,e]$.
\item[(c)] There exists $C<0$ such that $C^{-1}\rho(a_k)\leq
diam(R_k) \leq C \rho(a_k)$, where $a_k$ denotes the centre of $R_k$.
\item[(d)] $\bigcup_k \overline{R_k}=\C$ and the interiors of $R_k$ are
 distinct.  
\end{itemize}
\end{theorem}
\begin{remark}
Dividing the original measure by $s \in \R^+$ we obtain a 
partition of $\C$ into quasi-squares of mass $s$. 
\end{remark} 
\begin{proof}

It is enough to partition the plane in quasi-squares of constant entire mass,
because by an stopping-time argument of \cite{Or} these can then be split
into quasi-squares of mass $1$. 

We construct our partition recursively.  
We start with a rectangle centred at $0$ of entire mass,
and with sidelengths $l\leq L$ so that $l\geq L/2$ and $ l^{1-\beta}\geq
12\sqrt 2 C_0$, where $\beta$ and $C_0$ are given in \eqref{cota-rho}  
(rectangle $ABCD$ in the picture).  Consider next a square
$Q_1$ centred at 0 of sidelength $3L$ ($A_1B_1C_1D_1$ in the  picture) and
define $R$ as the quasi-square with vertices $ABB'A'$, where $A'$ and $B'$ are
points on the same side of $Q_1$ taken so that $0 \notin R$. We want to make
$R$ a little bigger, to make sure that its mass is entire, and we want to do
that keeping control on the ratio of sides. Consider the rectangle  $AB\tilde
B\tilde A$, where $\tilde A, \tilde B$ are taken with  $|A \tilde A| = |B
\tilde B|=2|AA'|$. Denote by $R'$ the rectangle $A'B'\tilde B \tilde A$ added
to $R$. For $\lambda \in  R'$, 
\[  
\frac{\rho(\lambda)}{l} \leq \frac{6\sqrt 2 \rho(\lambda)}{|\lambda|} \leq 
\frac{6\sqrt{2}C_0}{|\lambda|^{1-\beta}} \leq  \frac{6\sqrt{2}C_0}
{l^{1-\beta}}
\leq \frac{1}{2} .
\]
Since the sides of $R'$ have length bigger or equal than $l$ we deduce that
$R'$ contains a disk of centre $\lambda$ and radius $\rho(\lambda)$, hence its
mass is at least $1$. This shows that there exists a rectangle $R_1$ 
($AA''B''B$ in the picture) of entire
mass between the original $R$ and the ``doubled'' $R'$.

We finish the first step of the process by constructing the analogous 
quasi-square $R_2$ of entire mass at the opposite side of $R$ ($CC''D''D$ in
the picture).

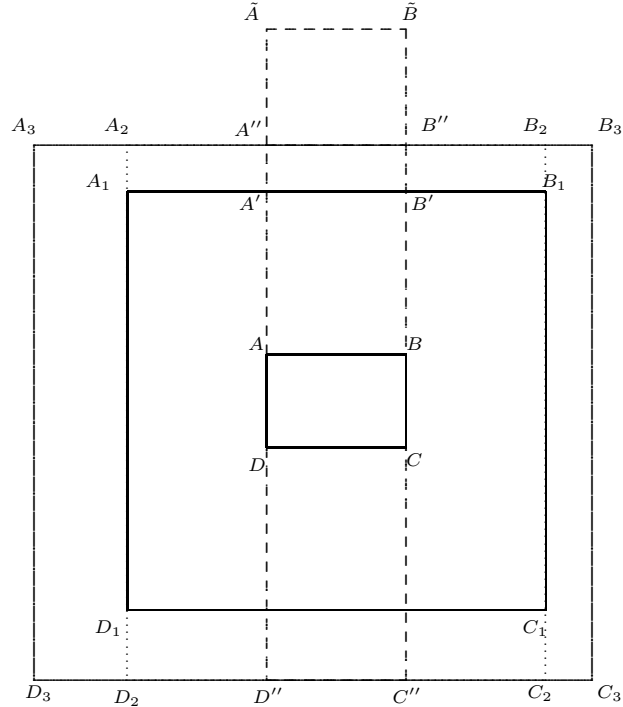
\begin{figure}[ht]\label{dibuix}
\begin{center}
\font\thinlinefont=cmr5
\begingroup\makeatletter\ifx\SetFigFont\undefined
\def\x#1#2#3#4#5#6#7\relax{\def\x{#1#2#3#4#5#6}}%
\expandafter\x\fmtname xxxxxx\relax \def\y{splain}%
\ifx\x\y   
\gdef\SetFigFont#1#2#3{%
  \ifnum #1<17\tiny\else \ifnum #1<20\small\else
  \ifnum #1<24\normalsize\else \ifnum #1<29\large\else
  \ifnum #1<34\Large\else \ifnum #1<41\LARGE\else
     \huge\fi\fi\fi\fi\fi\fi
  \csname #3\endcsname}%
\else
\gdef\SetFigFont#1#2#3{\begingroup
  \count@#1\relax \ifnum 25<\count@\count@25\fi
  \def\x{\endgroup\@setsize\SetFigFont{#2pt}}%
  \expandafter\x
    \csname \romannumeral\the\count@ pt\expandafter\endcsname
    \csname @\romannumeral\the\count@ pt\endcsname
  \csname #3\endcsname}%
\fi
\fi\endgroup
\mbox{\beginpicture
\setcoordinatesystem units <1.04987cm,1.04987cm>
\unitlength=1.04987cm
\linethickness=1pt
\setplotsymbol ({\makebox(0,0)[l]{\tencirc\symbol{'160}}})
\setshadesymbol ({\thinlinefont .})
\setlinear
%
%
\linethickness= 0.500pt
\setplotsymbol ({\thinlinefont .})
\setdashes < 0.1270cm>
{\plot  6.096 18.728  6.096 15.784 /
\plot  6.096 15.784  7.863 15.784 /
\plot  7.863 15.784  7.863 18.728 /
}%
%
%
\linethickness= 0.500pt
\setplotsymbol ({\thinlinefont .})
\setdots < 0.0953cm>
{\plot  4.331 15.784  9.631 15.784 /
\plot  9.631 15.784  9.631 22.557 /
\plot  9.631 22.557  4.331 22.557 /
\plot  4.331 22.557  4.331 15.784 /
}%
%
%
\linethickness= 0.500pt
\setplotsymbol ({\thinlinefont .})
\setdashes < 0.1270cm>
{\plot  6.096 19.907  6.096 22.557 /
\plot  6.096 22.557  7.863 22.557 /
\plot  7.863 22.557  7.863 19.907 /
}%
%
%
\linethickness= 0.500pt
\setplotsymbol ({\thinlinefont .})
\setsolid
{\putrule from  6.096 19.907 to  7.863 19.907
\putrule from  7.863 19.907 to  7.863 18.728
\putrule from  7.863 18.728 to  6.096 18.728
\putrule from  6.096 18.728 to  6.096 19.907
}%
%
%
\linethickness= 0.500pt
\setplotsymbol ({\thinlinefont .})
{\putrule from  9.631 21.969 to  4.331 21.969
\putrule from  4.331 21.969 to  4.331 16.667
\putrule from  4.331 16.667 to  9.631 16.667
\putrule from  9.631 16.667 to  9.631 21.969
}%
%
%
\linethickness= 0.500pt
\setplotsymbol ({\thinlinefont .})
{\plot  3.152 15.784 10.219 15.784 /
\plot 10.219 15.784 10.219 22.557 /
\plot 10.219 22.557  3.152 22.557 /
\plot  3.152 22.557  3.152 15.784 /
}%
%
%
\put{\SetFigFont{7}{8.4}{rm}{$A$}%
} [lB] at  5.861 19.967
%
%
\put{\SetFigFont{7}{8.4}{rm}{$B'$}%
} [lB] at  7.923 21.734
%
%
\put{\SetFigFont{7}{8.4}{rm}{$A''$}%
} [lB] at  5.685 22.676
%
%
\put{\SetFigFont{7}{8.4}{rm}{$B''$}%
} [lB] at  8.041 22.735
%
%
\put{\SetFigFont{7}{8.4}{rm}{$C_1$}%
} [lB] at  9.337 16.372
%
%
\put{\SetFigFont{7}{8.4}{rm}{$A_2$}%
} [lB] at  4.034 22.735
%
%
\put{\SetFigFont{7}{8.4}{rm}{$B_2$}%
} [lB] at  9.337 22.735
%
%
\put{\SetFigFont{7}{8.4}{rm}{$C_2$}%
} [lB] at  9.396 15.549
%
%
\put{\SetFigFont{7}{8.4}{rm}{$D_2$}%
} [lB] at  4.153 15.490
%
%
\put{\SetFigFont{7}{8.4}{rm}{$B_3$}%
} [lB] at 10.279 22.735
%
%
\put{\SetFigFont{7}{8.4}{rm}{$C_3$}%
} [lB] at 10.279 15.549
%
%
\put{\SetFigFont{7}{8.4}{rm}{$D_3$}%
} [lB] at  3.033 15.549
%
%
\put{\SetFigFont{7}{8.4}{rm}{$C$}%
} [lB] at  7.863 18.493
%
%
\put{\SetFigFont{7}{8.4}{rm}{$B$}%
} [lB] at  7.863 19.967
%
%
\put{\SetFigFont{7}{8.4}{rm}{$D$}%
} [lB] at  5.861 18.434
%
%
\put{\SetFigFont{7}{8.4}{rm}{$A_1$}%
} [lB] at  3.799 22.028
%
%
\put{\SetFigFont{7}{8.4}{rm}{$A'$}%
} [lB] at  5.743 21.734
%
%
\put{\SetFigFont{7}{8.4}{rm}{$B_1$}%
} [lB] at  9.572 22.028
%
%
\put{\SetFigFont{7}{8.4}{rm}{$D_1$}%
} [lB] at  3.918 16.372
%
%
\put{\SetFigFont{7}{8.4}{rm}{$A_3$}%
} [lB] at  2.857 22.735
%
%
\linethickness= 0.500pt
\setplotsymbol ({\thinlinefont .})
\setdashes < 0.1270cm>
{\plot  6.096 22.557  6.096 24.031 /
\plot  6.096 24.031  7.863 24.031 /
\plot  7.863 24.031  7.863 22.557 /
}%
%
%
\put{\SetFigFont{7}{8.4}{rm}{$\tilde A$}%
} [lB] at  5.802 24.149
%
%
\put{\SetFigFont{7}{8.4}{rm}{$\tilde B$}%
} [lB] at  7.804 24.149
%
%
\put{\SetFigFont{7}{8.4}{rm}{$D''$}%
} [lB] at  5.920 15.490
%
%
\put{\SetFigFont{7}{8.4}{rm}{$C''$}%
} [lB] at  7.688 15.490
\linethickness=0pt
\putrectangle corners at  2.857 24.301 and 10.279 15.433
\endpicture}

\end{center}
\end{figure}

Consider next the rectangle $Q_2$  limited by the segments $(A''B'')$,
$(C''D'')$, $(B_1C_1)$, $(D_1A_1)$  (the rectangle $A_2 B_2 C_2 D_2$ in 
the picture).
We iterate the  process above to each of the rectangles $B''B_2C_2C''$ and
$D''D_2A_2A''$, thus  obtaining two new quasi-squares $R_3=B''B_3C_3C''$ and 
$R_4=D''D_3A_3A''$ of entire mass. 

All in all, we obtain a new quasi-square $Q_3:=A_3B_3C_3D_3$ with ratio of
sides  lying in  $[1/2,2]$ which is a disjoint union of 5 quasi-squares of
entire mass. From here we repeat the process, taking $Q_3$ in place of the
original $R$, and continue recursively to obtain the ``partition'' of $\C$. By
construction we have (a), (b) and (d).

To prove (c) assume that $R$ is a quasi-square of mass 1, centre $a$ and
sidelengths $l$, $L$. Here $R\subset D(a,L\sqrt 2)$, hence  $\rho(a) \gtrsim L
\gtrsim diam(R)$. Also, $D(a,l) \subset R$ and $diam(R)\lesssim l \lesssim
\rho(a)$.  
\end{proof}

Lemmas \ref{discs} and \ref{christ} give  control on how big a disc
$D^r(\zeta)$ can be when $\zeta\in D^R(z)$. We will need another result along
the same lines.

Given a doubling measure $\mu$ and given 
$z\in\C$ and $0<r<R$, consider the associated regions
\begin{eqnarray*}
F_r(z,R)= \{\zeta : D^r(\zeta) \subset D^R(z)\}
\qquad\textrm{and}\qquad
G_r(z,R)=
\bigcup\limits_{\zeta \in D^R(z)} D^r(\zeta)\  .
\end{eqnarray*}
By definition $F_r(z,R)\subset D^R(z) \subset G_r(z,R)$.
Let $\gamma$ be the constant given by Lemma \ref{discs},  
and $\varepsilon,k$ the constants in \eqref{cota-inferior}.

\begin{lemma}\label{flors} Let $r>0$ be fixed. There exists $c>0$ such that 
if $\epsilon(R)=c(r^k/R^\varepsilon)^\gamma$, 
for all $z\in\C$ and $R>r$ we have
\begin{itemize}
\item[(a)] $G_r(z,R)\subset D^{R+\epsilon(R)}(z)$.
\item[(b)] $D^{R-\epsilon(R)}(z) \subset F_r(z,R)$.
\end{itemize}
\end{lemma}
  
\begin{proof}
Applying Lemma \ref{discs} to
$D^r(\zeta)$ and $D^R(z)$, and using \eqref{cota-inferior}, we have
\[
\Bigl(\frac{R^\varepsilon}{r^k}\Bigr)^\gamma\lesssim
\frac{R\rho(z)}{r\rho(\zeta)}\lesssim
\Bigl(\frac{R^k}{r^\varepsilon}\Bigr)^{1/\gamma} .
\]

(a) If $\zeta\in D^R(z)$ we have
$R\rho(z)+r\rho(\zeta)\leq R\rho(z)
(1+c(r^k/R^\varepsilon)^\gamma)$
for some $c>0$.

(b) $D^r(\zeta)\subset D^R(z)$ when $|\zeta-z|+r\rho(\zeta)\leq R\rho(z)$. For
$\zeta\in D^{R-\epsilon(R)}(z)$
\[
|\zeta-z|+r \rho(\zeta)\leq (R-\epsilon(R))
\rho(z)+c_1 R \rho(z)\Bigl(\frac{r^k}{R^\varepsilon}\Bigr)^\gamma  .
\]
Thus if
$(R-\epsilon(R))\rho(z)+c R \rho(z)( r^k/R^\varepsilon)^\gamma 
\leq R\rho(z)$ 
we have $D^{R-\epsilon(R)}(z)\subset F_r(z,R)$.  
\end{proof}

\begin{corollary}\label{flors-quadrats}
Let $\{R_k\}_k$ be a partition of $\C$, as in Theorem \ref{partition}.
Define
\begin{eqnarray*}
F(z,R)= \bigcup\limits_{k:R_k\subset D^R(z)} R_k    
\qquad\textrm{and}\qquad
G(z,R)=
\bigcup\limits_{k:R_k\cap D^R(z)\neq\emptyset} R_k \  .
\end{eqnarray*}
There exists a positive function $\epsilon(R)$ with 
$\lim_{R\to\infty} \epsilon(R)/R=\infty$ and such that
for all $z\in\C$ and $R>0$
\begin{itemize}
\item[(a)] $G(z,R)\subset D^{R+\epsilon(R)}(z)$.
\item[(b)] $D^{R-\epsilon(R)}(z) \subset F(z,R)$.
\end{itemize}
\end{corollary}

\begin{proof} As the previous Lemma, using Theorem \ref{partition}(c).
\end{proof}

We finish with a result showing that the measure of a disk cannot be
too concentrated near its border.

\begin{lemma}\label{corona-zero}
Let $\epsilon(r)$ be a positive function such 
that $\lim\limits_{r\to\infty} \epsilon(r)/r=0$. Then
\[
\lim_{r\to\infty}
\frac{\mu(D^{r+\epsilon(r)}(z))}{\mu(D^r(z))}=
\lim_{r\to\infty}
\frac{\mu(D^{r-\epsilon(r)}(z))}{\mu(D^r(z))}=1
\]
uniformly in $z\in\C$.
\end{lemma}

The proof is based in the following projection of the measure $\mu$.

\begin{lemma}\label{nu_z}
For every $z\in\mathbb C$ define the measure $\nu_z$ on $\mathbb R^+$ by
\[
\nu_z(A)=\mu(\{\zeta=z+re^{i\theta} : r\in A\})\qquad A\subset\mathbb R^+ .
\]
Then $\nu_z$ is doubling and there exists $K$ independent of $z$ such that
$C_{\nu_z}\leq K C_\mu$.
\end{lemma}

\begin{proof}
Given $x\in\mathbb R^+$ and $r>0$ let $I^r(x)=(x-r,x+r)\cap\R^+$. 
We want to see that
\[
\nu_z(I^{2r}(x))\leq K C_\mu\; \nu_z (I^r(x))
\]
for all $z\in\mathbb C$, $x\in\mathbb R^+$ and $r>0$.

Let $A_z^r(x)=\{\zeta=z+se^{i\theta} : s\geq 0\, , |s-x|<r\}$. By definition
$ \nu_z(I^{2r}(x))=\mu(A^{2r}_z(x))$.
Split $A^{2r}_z(x)$ into $k:=[\frac {2\pi}{4r}]$ sectors
\[
S_j=\bigl\{\zeta=z+se^{i\theta} : s\geq 0\, , |s-x|<2r\; ,\ 
(j-1)\frac{2\pi}k\leq \theta <j\frac{2\pi}k \bigr\}\qquad j=1,\dots,k .
\]
Being $\mu$ doubling  there exists $K>0$ such that
$\mu(S_j)\leq K C_\mu\; \mu(\tilde S_j)$, 
where $\tilde S_j$ is half the sector $S_j$, i.e.
\[
\tilde S_j=\bigl\{\zeta=z+se^{i\theta} : s\geq 0\, , |s-x|<r\; ,\ 
(j-1)\frac{2\pi}k+\frac{2\pi}{4k}<\theta<j\frac{2\pi}k-\frac{2\pi}{4k}
\bigr\} .
\]
Since the $\tilde S_j$'s are disjoint and
$\cup_j\tilde S_j\subset A_z^r(x)$, we get
\begin{eqnarray*}
\nu_z(I^{2r})(x)&=&\mu(A^{2r}_z(x))=\sum_{j=1}^k \mu (S_j)\leq 
K C_\mu \sum_{j=1}^k
\mu(\tilde S_j)
\leq  K C_\mu\; \mu(A_z^r(x))\\
&=& K C_\mu\; \nu_z(I^r(x)).
\end{eqnarray*}
\end{proof}

\begin{proof}[Proof of Lemma \ref{corona-zero}]
It is enough to see that
\[
\lim_{r\to\infty} 
\frac{\mu( D^{r+\epsilon(r)}(z)\setminus D^r(z))}{\mu (D^r(z))}=0
\]
uniformly in $z$.
By definition of $\nu_z$ we have
\begin{eqnarray*}
\frac{\mu( D^{r+\epsilon(r)}(z)\setminus D^r(z))}{\mu (D^r(z))}=
\frac{\nu_z \bigl((r\rho(z),(r+\epsilon(r))\rho(z)\bigr)}{\nu_z((0,r))} ,
\end{eqnarray*}
and by the corresponding version of Lemma \ref{discs} for doubling measures in
$\mathbb R^+$, and by Lemma \ref{nu_z}, there exists $K>0$
independent of $z$ such that
\[
\frac{\nu_z \bigl((r\rho(z),(r+\epsilon(r))\rho(z)\bigr)}{\nu_z((0,r))}
\leq K \left(\frac{\epsilon(r)}r\right)^\gamma .
\]
\end{proof}

\begin{remark}\label{corona-rectangles} 
An analogous result is true if in the definition of $\nu_z$ we
use, instead of a radial projection with respect to $z$, a projection
associated to quasi-squares of a fixed ratio $\alpha\in[e^{-1},e]$ ($e$ is the
constant of 
Theorem~\ref{partition}(b)).
Let $Q_\alpha^r(z)$ denote the rectangle with vertices
$z+r(1+i\alpha)$, $z+r(1-i\alpha)$, $z-r(1+i\alpha)$ and $z-r(1-i\alpha)$.
Given $z\in\mathbb C$ consider the measure $\nu_z$ in $\mathbb R$ such that 
\[
\nu_z(I^r(x))=\mu (Q_\alpha^{x+r}(z)\setminus Q_\alpha^{x-r}(z))
\]
on any interval $I^r(x)$.
As before, there exists $K>0$ independent of $z\in\mathbb C$ and 
$\alpha\in[e^{-1},e]$ such that $\nu_z$ is doubling with
$C_\nu\leq K C_\mu$. Therefore, if $R_\alpha^r(z):=Q_\alpha^{r\rho(z)}(z)$,
\[
\lim_{r\to\infty} 
\frac{\mu (R_\alpha^{r+\epsilon(r)}(z))}{\mu(R_\alpha^{r}(z))}=
\lim_{r\to\infty} 
\frac{\mu (R_\alpha^{r+\epsilon(r)}(z))}{\mu(R_\alpha^{r}(z))}=1
\]
uniformly in $z$.
\end{remark}

\subsection{Flat weights}\label{flat-weights}
In this section we describe the weights $\omega$ appearing in the spaces 
$\Fpa$.

\begin{definition}
A positive measurable 
function $\omega$ is called a \emph{flat weight for} $\phi$
if there exists  $C>0$ such that for all $z,\zeta\in\C$
\begin{equation}\label{fw1}
 |\log \omega(z)-\log \omega(\zeta)|\leq C
(1+\log^+ d_\phi(z,\zeta)) .
\end{equation}
\end{definition}

The class of flat weights associated to $\phi$ will be denoted by
$\Omega_\phi$.

Notice that the product $\omega\widehat{\omega}$ of two weights
$\omega,\widehat{\omega}\in\Omega_\phi$  belongs to $\Omega_\phi$ as well.
Also, if $\omega\in\Omega_\phi$ then $\omega^\alpha\in\Omega_\phi$, for all
$\alpha\in\R$.

Besides the obvious $\omega=1$, important examples of flat weights for $\phi$
are the functions $\omega=\rho^\alpha$, $\alpha\in\R$. This is seen applying
Lemma~\ref{christ} and Lemma~\ref{distance}.

Furthermore, the weights $\omega\in\Omega_\phi$ can be assumed to satisfy
\begin{equation}\label{fw2}
\bigl|1-\frac{\omega(z)}{\omega(\zeta)}\bigr|\leq
C d_\phi(z,\zeta)\qquad\textrm{if}\quad d_\phi(z,\zeta)\leq 1\ .
\end{equation}
If the original weight $\omega$ does not satisfy this condition, replace it by
the regularisation
\[
\tilde\omega(z)=\frac 1{\rho^2(z)}\int_{D(z)} \omega .
\]
It is clear, by \eqref{fw1}, that there exists $C>0$ such that
$C^{-1}\leq|\omega/\tilde\omega|\leq C$, hence the spaces of functions 
and sequences
associated to the weights $\omega$ and $\tilde\omega$ are the same.
On the other hand
\[
\bigl|\frac {\tilde\omega(\zeta)-\tilde\omega(z)}{\tilde\omega(\zeta)}\bigr|
\leq
\bigl|\frac 1{\rho^2(\zeta)}\bigl[ \int_{D(\zeta)}
\frac{\omega}{\tilde\omega(\zeta)}-
\int_{D(z)}\frac{\omega}{\tilde\omega(\zeta)}\bigr]\bigr|
+ \bigl|\frac 1{\rho^2(\zeta)}-\frac 1{\rho^2(z)}\bigr|
\int_{D(z)}\frac{\omega}{\tilde\omega(\zeta)}.
\]
Assuming that $d_\phi(z,\zeta)\leq 1$, from
\eqref{fw1}, \eqref{lipschitz} and Lemma~\ref{distance}(a) 
we deduce that
\begin{eqnarray*}
\bigl|\frac {\tilde\omega(\zeta)-\tilde\omega(z)}{\tilde\omega(\zeta)}\bigr|
&\lesssim & 
\frac{\sigma[(D(\zeta)\cup D(z))\setminus (D(\zeta)\cap D(z))]}{\rho^2(\zeta)}+
\frac{|\rho(z)-\rho(\zeta)||\rho(z)+\rho(\zeta)|}{\rho^2(\zeta)}\\
&\lesssim & \frac{\rho(z)|\zeta-z|}{\rho^2(\zeta)}+
\frac{|\zeta-z|}{\rho(\zeta)}\lesssim
d_\phi(\zeta,z).
\end{eqnarray*}

\subsection{Local behaviour and regularisation of $\phi$}
Let us start with a result comparing the values
of $\phi$ in a disk with the value on its centre.

\begin{lemma}\label{cota-disc}
For every $K>0$ there exists $A=A(K)>0$ such that for all $z\in\C$
\[
\sup_{w \in D^K(z)}|\phi(w)-\phi(z)-h_z(w)|\leq A ,
\]
where $h_z$ is a harmonic function in $D^K(z)$ with $h_z(z)=0$.
\end{lemma}

\begin{proof} The proof is as in \cite[Lemma 1]{OrSe}.
On each $D^K(z)$ decompose
\begin{equation}\label{green}
\phi(w)=\phi(z)+h_z(w)+\int_{D^K(z)} (G(w,\eta)-G(z,\eta))\; 
\Delta\phi(\eta)  ,
\end{equation}
where $G$ is the Green function of the disc $D^K(z)$ and $h_z$ is a
harmonic function in $D^K(z)$ such that $h_z(z)=0$. By 
Lemma~\ref{integrals}(a)
\[
\sup_{z\in\C} \int_{D^K(z)} \log\frac{K\rho(z)}{|z-\eta|}
\Delta\phi(\eta)<\infty
\]
and the result holds.
\end{proof}

We have seen in the previous section   that $\rho_\phi(z)$ is Lipschitz (see
(\ref{lipschitz})). Also, because of Lemma~\ref{discs}, $\phi$ is H\"older
continuous of some positive order on every bounded subset of $\C$ (see
\cite[Lemma 2.8]{Ch}). More regularity can be attained by taking a suitable
weight $\psi$ equivalent to $\phi$.

\begin{theorem}\label{regularitzacio}
Let $\phi$ be subharmonic with $\Delta\phi$ doubling. There exist 
$\psi\in \mathcal C^\infty (\C)$ subharmonic and $C>0$ such that 
$|\psi-\phi|\leq C$, $\Delta \psi$ is a doubling measure and
$\Delta\psi \simeq 1/\rho_\psi^2\simeq 1/\rho_\phi^2$. Moreover
$|\nabla (\Delta \psi) | \lesssim 1/\rho_\phi^3$.
\end{theorem}

Since the spaces of functions and sequences considered do not change if 
$\phi$ is replaced by $\psi$, from now on we will assume that $\phi\in\mathcal
C^\infty (\C)$,  $\Delta\phi\simeq 1/\rho^2$ and $|\nabla (\Delta \phi) |
\lesssim 1/\rho^3$.

In the proof of this result we will need to partition $\C$ and
discretize the measure.

\begin{lemma}\label{sigma} Let $\mu$ be a positive doubling measure
in $\C$.
Fix $m\in\N$.
There exist $k\in \N$ and $C>0$ such that for any partition $\{R_p\}_p$ 
as in Theorem~\ref{partition} with $\mu(R_p)=mk$ there are points
$\lambda_1^{(p)},\dots,\lambda_{mk}^{(p)}\in CR_p$ such that 
\begin{itemize}
\item[(a)] $\mu_p=\mu_{|R_p}$ and $\nu_p=\sum\limits_{j=1}^{mk} 
\delta_{\lambda_j^{(p)}}$
have the same first $m$ moments.
\item[(b)] $\Lambda=\{\lambda_j^{(p)}\}_{p,j}$ is a $\rho$-separated sequence.
\end{itemize}
\end{lemma}

\begin{proof} 
By Lemma 5 of \cite{Or}, there exists $k \in \N$ such that for
all measure   $\mu_p$ supported in a rectangle $R_p$ with total mass $mk$,
there are points $\sigma_1^{(p)},\ldots \sigma_k^{(p)}\in R_p$ such that
$\mu_p$ and  $m\sum_{j=1}^k \delta_{\sigma_j^{(p)}}$ have the same first $m$
moments.

In order to get a separated sequence replace each $\sigma_j^k{(p)}$ by $m$  
points $\gamma_{j,l}^{(p)}=\sigma_{j}^{(p)}+\tau_{j}^{(p)}  e^{i 2\pi l/m}$, 
$l=0,\dots, m-1$, lying on a circle around  $\sigma_{j}^{(p)}$. Since for all
polynomials $p$ of degree less than $m-1$  

\begin{equation*}
m\; p(\sigma_{j}^{(p)}) = \sum_{l=0}^{m-1}  p( \gamma_{j,l}^{(p)})  ,
\end{equation*}
the measures $\mu_p$ and $\sum_{j,l} \delta_{\gamma_{j,l}^{(p)}}$ have still  
the same first $m$ moments.  We will be done as soon as we see that the
$\tau_{j}^{(p)}$ can be  chosen uniformly bounded and so that 
$\Lambda=\{\gamma_{j,l}^{(p)}\}$ is $\rho$-separated. For this we use a
Besicovitch's lemma: the family $\{R_p\}_p$ can be split in $q$ families
$\{R_{p}^{1}\}_{p \in I_1},\dots \{R_{p}^{q}\}_{p \in I_q}$ such that two
rectangles of the same family are far apart, in the sense that $MR_p^l \cap
MR_{p'}^l =\emptyset$, $p\neq p'$, for some large constant $M$. For the first
family $\{R_{p}^{1}\}_{p \in I_1}$, it is easy to choose $\tau_j^{(p)}$ 
such that the resulting sequence $\Gamma_1=\{ \gamma_{j,l}^{(p)} : p \in
I_{1}; j=1,\dots,k ; l=0,\dots m-1 \} $  is $\rho$-separated. Next we choose
$\tau_j^{(p)}$,  $p\in I_2$, so that $\Gamma_2 \cap \Gamma_1$ is
$\rho$-separated, where  $\Gamma_2=\{ \gamma_{j,l}^{(p)} : p \in I_{2};
j=1,\dots,k ; l=0,\dots, m-1 \}  $. Choosing $\tau_j^{(p)}$ recursively in this
way we obtain $\Lambda=\Gamma_1 \cup \ldots \cup \Gamma_q$  $\rho$-separated.  
\end{proof}

\begin{proof}[Proof of Theorem~\ref{regularitzacio}]
For any  $M$ (to be chosen later) consider $k\in\N$ as in Lemma~\ref{sigma} 
and a partition $\{R_p\}_p$ as in Theorem~\ref{partition}.
Take then the sequence $\Lambda=\{\lambda_j^{(p)} \}_{j,p}$ given by
Lemma~\ref{sigma}. Recall that $\lambda_j^{(p)}\in CR_p$,
$\mu(R_p)=Mk$ and that the measures   $\mu_p$
and  $\nu_p=M\sum\limits_{j=1}^k \delta_{\lambda_j^{(p)}}$ have the same first
$M$ moments. 

By Theorem \ref{partition}(c) there exists $r>0$ such that  $C R_p\subset
D^r(\lambda_j^{(p)})$ for any $p\in\N$ and $i\leq k$.  Furthermore, 
by construction of $\{R_p\}_p$ there exists
$q\in\N$ such that any $z\in\C$ lies in at most $q$ disks 
$D^r(\lambda_j^{(p)})$.

We now regularise $\nu_p$ by setting
\[
\tilde\nu_p=\sum_{j=1}^{Mk}
\frac{\X\bigl(\frac{|z-\lambda_j^{(p)}|}{r\rho(\lambda_j^{(p)})}\bigr)}
{\int \X\bigl(\frac{|z-\lambda_j^{(p)}|}{r\rho
(\lambda_j^{(p)})}\bigr)} ,
\]
where $\X$ is a smooth non-negative cut-off function of one real variable such
that $\X(t)=1$ if $|t|<1$, $\X(t)=0$ if $|t|>2$ and $|\X'|$ is bounded.

Notice that $\tilde\nu_p$ and $\mu_p$ have the first $M$ moments. Indeed, by
the mean value property
\[
\int_{\C} z^l d\tilde\nu_p=\sum_{i=1}^{Mk} (\lambda_j^{(p)})^l\qquad l=0,
\dots,M-1 .
\]
Define $\tilde\nu=\sum\limits_{p=1}^\infty \tilde\nu_p$ and
\[
\psi(z)=\phi(z)+\frac 1{2\pi}\int_{\C}
\log |z-\zeta| (\tilde\nu-\Delta\phi)(\zeta) .
\]
We claim that $\tilde\nu$ is a doubling measure. 
The proof of this fact is a bit technical and will be deferred to the end.

By definition $\Delta\psi=\tilde\nu$. Also, $\tilde\nu(z)$ is a sum of at most
$q$ terms of order $ 1/\rho^2(\lambda_j^{(p)})$, with $z\in
D^r(\lambda_j^{(p)})$. Therefore $\Delta\psi\simeq 1/\rho_\phi^2$ and
 $|\nabla (\Delta\psi) | \lesssim
1/\rho_\phi^3$. 
In particular
\[
\int_{D_\phi(z)}\Delta\psi(\zeta)\simeq \int_{D_\phi(z)}
\frac {dm(\zeta)}{\rho_\phi^2(\zeta)} 
\simeq 1 ,
\]
hence $\rho_\phi\simeq\rho_\psi$.

Let us  show next that $|\psi-\phi|\leq C$ for some $C>0$.

Let $a_p$ denote the centre of $R_p$.
Assume $z\in R_{p_0}$ and let $I_{p_0} = \{ p\in \N : 
d_\phi(a_p,a_{p_0}) \leq 10r \}$. 
Remark that for $p\notin I_{p_0}$, 
$\zeta \in supp(\tinu_p) $ and $z\in supp(\tinu_{p_0})$ 
we have $ d_\phi(z,\zeta) \simeq d_\phi(a_p,a_{p_0})$. 
Indeed, this follows from 
\begin{equation*}
|\zeta-a_p|\leq 3 r\rho(a_p) \leq \frac{3 }{10} |a_p-a_{p_0}|, 
\end{equation*}
the analogous estimate for $|z-a_{p_0}|$ and Lemma~\ref{distance}. This yields 
\begin{equation} \label{majoration}
\int_\C d_\phi^{-M}(z,\zeta)  \tinu_p(\zeta)  
\lesssim \int_\C d_\phi^{-M}(z,\zeta) \mu_p(\zeta)
\qquad z\in I_{p_0}\; ,\ p\notin I_{p_0} .
\end{equation}
We split
\[
2\pi(\psi(z)-\phi(z))=\sum_{p\in I_{p_0}}\int_{\C}\log|z-\zeta|
(\tilde\nu_p-\mu_p)(\zeta)
+\sum_{p\notin I_{p_0}}\int_{\C}\log|z-\zeta|(\tilde\nu_p-\mu_p)(\zeta)
\]
and estimate each sum separately.

Let $p_M$ denote the $M$-th Taylor polynomial of
$\log(|z-\zeta|/\rho(z))$. Since $\tilde\nu_p-\mu_p$ have vanishing
moments of order less or equal to $M$, we can estimate
\begin{eqnarray*}
I_1:&=&\bigl|\sum_{p\notin I_{p_0}}
\int_{\C}\log|z-\zeta| (\tilde\nu_p-\mu_p)(\zeta)\bigr|\\
&=&\bigl|\sum_{p\notin I_{p_0}}
\int_{\C}(\log\frac{|z-\zeta|}{\rho(z)}-p_M(\zeta)) 
(\tilde\nu_p-\mu_p)(\zeta)\bigr|\leq
\sum_{p\notin I_{p_0}} \int_{\C}\bigl(\frac{\rho(z)}{|z-\zeta|}\bigr)^M
(\tilde\nu_p + \mu_p)(\zeta) .
\end{eqnarray*}
Taking $M$ big enough and using (\ref{majoration}) and
Lemmas~\ref{distance}(a) and \ref{integrals}(b), 
\begin{equation*} 
I_1\lesssim  \int_{\C\setminus D^r(z) }
\bigl(\frac{\rho(z)}{|z-\zeta|}\bigr)^M \mu(\zeta) \lesssim
\int_{\C\setminus B(z,Cr) }\frac{d\mu(\zeta)}{d_\phi^{\delta M}(z,\zeta)}
\leq C.
\end{equation*}

For the remaining term we use again the moment condition 
together with the fact that for $p \in  I_{p_0}$ there exists
$\gamma$ such that $\cup\{supp(\tilde\nu_p ),\ p\in I_{p_0}\} \subset
D_\phi^\gamma(z)$. Thus
\begin{eqnarray*}
I_2:&=&\bigl|\sum_{p\in I_{p_0}}
\int_{\C}\log|z-\zeta|  (\tilde\nu_p-\mu_p)(\zeta)\bigr|=
\bigl|\sum_{p\in I_{p_0}}
\int_{\C}\log\bigl(\frac{2\gamma\rho(z)}{|z-\zeta| }\bigr)
(\tilde\nu_p-\mu_p)(\zeta)\bigr|\\
&\lesssim& \int_{D^{\gamma}(z)} \log\bigl(\frac{2\gamma\rho(z)}{|z-\zeta|
}\bigr)(\tilde\nu+\mu)(\zeta) . \qquad\qquad\ 
\end{eqnarray*}

By Lemma~\ref{integrals}(a) this is finite.

We prove now that $\tinu$ is doubling.  We first show that it is doubling for
big balls, i.e. there exist $R_0>0$ and a constant $C$ depending only on the
doubling constant $C_{\Delta\phi}$ of $\Delta\phi$ such that for all $R>R_0$ we
have $\tinu(D^R(a)) \leq C \tinu(D^{R/2}(a))$.

As in Corollary \ref{flors-quadrats}, define   
\begin{eqnarray*}
F(a,R)  =  \bigcup_{p: R_p\subset D^R(a)} R_p 
\qquad\textrm{and}\qquad
G(a,R) = \bigcup_{p: R_p\cap D^R(a)  \neq \emptyset } R_p.
\end{eqnarray*}
Since $\tinu(R_p) \simeq \int_{R_p} d\sigma/\rho^2
\simeq \mu(R_p)$,
we see that
$\tinu(F(a,R)) \simeq \mu(F(a,R))$ and $\tinu(G(a,R)) 
  \simeq \mu(G(a,R))$.
By Corollary \ref{flors-quadrats}, also $
 D^{R-\epsilon( R)}  (a)  \subset  F(a,R)$ and
$D^{R+\epsilon( R)}(a)  \supset G(a,R)$.
This and the fact that $\mu$ is doubling yield   
\begin{equation*}
\tinu(D^R(a)) \leq \tinu(G(a,R)) \simeq  \mu(G(a,R)) \leq 
\mu(D^{R+\epsilon( R)}(a)) \leq C_{\Delta\phi}\; 
\mu(D^{1/2(R+\epsilon( R))}(a)),
\end{equation*}
and
\begin{equation*}
\tinu(D^{R/2}(a)) \geq   \tinu(F(a,R/2)) \simeq \mu(F(a,R/2))\geq 
\mu(D^{R/2-\epsilon( R/2)}(a)).
\end{equation*}
Therefore 
\begin{equation*}
\tinu(D^R(a))\leq C_{\Delta\phi}\; \tinu(D^{R/2}(a))
\frac{\mu(D^{1/2(R+\epsilon( R))}(a))}
{\mu(D^{R/2-\epsilon( R/2)}(a))}.
\end{equation*}
Lemma \ref{corona-zero} shows that the quotient converges to $1$ as $R$ goes
to infinity uniformly in $a$, so there exists $R_0$ such that
$\tinu(D^R(a))\leq 2C_{\Delta\phi}\; \tinu(D^{R/2}(a))$
for all $R\geq R_0$.

Corollary \ref{rho} implies that  $\tinu \simeq
1/\rho^2(a)$ on  $D^R(a)$ when $R\leq R_0$, so we deduce that 
$\tinu(D^R(a))\lesssim  \tinu(D^{R/2}(a))$.
\end{proof}

\subsection{The multiplier}
A basic tool in our approach is the use of the so-called multiplier:
an entire function $g$ such that $ |g|\simeq e^\phi$ outside a neighbourhood
of the zeros of $g$.

\begin{theorem}\label{multiplier}
Let $\phi$ be a subharmonic function such that $\Delta\phi$ is 
a doubling measure. 
There  exists an entire function $g$ such that
\begin{itemize}
\item[(a)] The zero-sequence $\mathcal Z(g)$ of $g$
is $\rho_\phi$-separated and 
$\displaystyle \sup\limits_{z\in\C}d_\phi(z, \mathcal Z(g))<\infty$.
\item[(b)] $|g(z)|\simeq e^{\phi(z)}
d_\phi(z, \mathcal Z(g))\ $ for all $z\in\C$.
\end{itemize} 
The function $g$ can be chosen so that, moreover, 
it vanishes on a prescribed 
$z_0\in\C$. 
We say that $g$ is a multiplier associated to $\phi$.
\end{theorem}
\begin{proof}
Take a partition $\{R_p\}$ of $\C$ with $\mu(R_p)=2\pi mN$ and consider the
sequence  $\Lambda$ given by Lemma~\ref{sigma}. For the sake of clarity we
write $R_p$ instead of $C R_p$ ($C$ is the constant of Lemma~\ref{sigma}). Note
that now $\{R_p\}_p$ is not a partition,  although there exists a uniform
constant $q$ such that all points of $\C$ lie in at most $q$ quasi-squares
$R_p$. As in Lemma~\ref{sigma},  denote $\mu_p=(1/2\pi) \mu_{|R_p}$ and let
$\nu_p$ be the sum of the $\lambda\in\Lambda$ associated to $R_p$. Recall that
$\mu_p$ and $\nu_p$ have the same first $m$ moments.

Let $g$ be a holomorphic function satisfying 
\begin{equation*}
\log|g| = \phi- \frac 1{2\pi}\int_\C\log|z-\zeta| 
(\Delta\phi-2\pi \sum_{\lambda \in \Lambda} \delta_\lambda),
\end{equation*}
which exists because the Laplacian of the term at the right hand side is a sum 
of Dirac masses. By definition $\mathcal Z(g)=\Lambda$, and the previous
construction ensures that (a) holds. 

Let us prove (b). Assume that $z \in R_{p_0}$ and let $I_{p_0}$ denote the set
of indices $p$ such that $R_p\cap  R_{p_0}\neq \emptyset $. As in the previous
proof, split
\begin{eqnarray*}
\log|g(z)|-\phi(z) = -\int_\C\log|z-\zeta| 
(\frac{\Delta\phi}{2\pi}-\sum_{\lambda \in 
\Lambda} \delta_\lambda)
= S_1(z)+S_2(z) ,
\end{eqnarray*}
where 
\begin{equation*}
S_1(z):= \sum_{p\notin I_{p_0}}   \int_\C\log|z-\zeta| (\nu_p-\mu_p)
\end{equation*}
and 
\begin{equation*}
S_2(z):= \sum_{p\in I_{p_0}}   \int_\C\log|z-\zeta| (\nu_p-\mu_p).
\end{equation*}
Again as in the proof of Theorem~\ref{regularitzacio}, using the Taylor
expansion of $\log|z-\zeta|$ together with  the moment condition one sees that 
$|S_1(z)|$ is bounded. 

For the second sum notice that there exits $\gamma>0$ 
such that $\cup_{p\in I_{p_0}} R_p \subset D^\gamma(z)$. Hence, denoting
$|z-\Lambda|=\inf\limits_{\lambda\in\Lambda} |z-\lambda|$, we get
\begin{eqnarray*}
S_2(z) &=& \sum_{p\in I_{p_0}}   \int_\C \log\frac{2\gamma\rho(z)}
{|z-\zeta|}(\mu_p-\nu_p)
\leq \int_{ D^\gamma(z)}  \log \frac{2\gamma\rho(z)}{|z-\zeta|}d\mu(\zeta) - 
\log\frac{2\gamma\rho(z)}{|z-\Lambda|} \\
&\leq& C_2 + \log\frac{|z-\Lambda|}{\rho(z)}.
\end{eqnarray*}
On the other hand, using the $\rho$-separation of $\Lambda$  
\begin{eqnarray*}
-S_2 (z)&\leq &  \sum_{p \in I_{p_0}} \sum_{\lambda \in R_p}  \log 
\frac{2\gamma\rho(z)}{|z-\lambda|} \leq \log\frac{\rho(z)}{|z-\Lambda|} 
+ C(\delta)\cdot \#\bigl(\Lambda\cap \cup_{p\in I_{p_0}} R_p\bigr).
\end{eqnarray*}
Since $\# I_{p_0}$ is uniformly 
bounded, this and the estimate of $S_1$ give:
\begin{equation*}
 \log\frac{|z-\Lambda|}{\rho(z)} - C \leq \log|g(z)|-\phi(z) 
 \leq \log\frac{|z-\Lambda|}{\rho(z)} + C'.
\end{equation*}
The result is then immediate from Lemma~\ref{distance}(a).
\end{proof}

Next we state a useful application of the multiplier. It is a result about
peak functions. These functions attain value $1$ at a given point and decay very
fast away from the point. They are very useful in the estimates of the Bergman
kernel and in the construction of solutions to the $\db$ equation.  Another
proof of  the following Lemma, using estimates for the $\db$-equation, can be
found in an Appendix. This second proof  is along the lines of 
\cite[Theorem~2.1]{FoSi89}, where a related result is proved.

\begin{theorem}\label{weight}
Take $\varepsilon>0$, $\omega\in\Omega_\phi$ and $m\in\N$. There exists $C>0$
such that for all $\eta\in\C$ there is an entire function $P_\eta$
with $P_\eta(\eta)=1$ and
\[
|P_\eta(z)|\leq C
e^{\varepsilon(\phi(z)-\phi(\eta))}
\frac{\omega(\eta)}{\omega(z)} 
\frac{1}{1+d_\phi^{m}(z,\eta)} .
\]
\end{theorem}

\begin{proof}
Let $h$ be a multiplier for $\varepsilon\phi$ (constructed as in 
Theorem~\ref{multiplier}) with zero sequence $\Sigma=\{\sigma_k\}_k$ and such
that $\{\eta\}\cup\Sigma$ is  $\rho$-separated. In particular  $|h(z)|\simeq
e^{\varepsilon\phi(z)} d_\phi(z,\Sigma) $. It follows from the construction of
the multiplier that for each $M\in\N$ there exists  $r>0$ such that $\#
(\Sigma\cap  B(\lambda,r))\gtrsim M$ for all $\eta\in\C$. Given 
$\sigma_1,\dots,\sigma_M\in \Sigma\cap  B(\lambda,r)$ define
\[
P_\eta(z)=c_\eta\frac{h(z)}{(z-\sigma_1)\cdots 
(z-\sigma_M)}\frac{\rho^M(\eta)}{e^{\varepsilon\phi(\eta)}} ,
\]
where $c_\eta$ is chosen so that $P_\eta(\eta)=1$.

Let us see first that there exists $c>0$ independent of $\eta$
with $c^{-1}\leq c_\eta\leq c$. Since 
$|\eta-\sigma_i|\simeq \rho(\eta)$, then 
\[
\frac 1{c_\eta}=\frac{h(\eta)}{(\eta-\sigma_1)\cdots 
(\eta-\sigma_M)}\frac{\rho^M(\eta)}{e^{\varepsilon\phi(\eta)}}
\simeq 
\frac{e^{\varepsilon\phi(\eta)}
d_\phi(\eta,\Sigma) }{\rho^M(\eta)} 
\frac{\rho^M(\eta)}{e^{\varepsilon\phi(\eta)}}=
d_\phi(\eta,\Sigma) \simeq 1 .
\]
We split the estimate of $|P_\eta(z)|$ into several regions. Let 
$\varepsilon>0$ be such that 
that the balls $B(\sigma_i,\varepsilon)$ and $B(\eta,\varepsilon)$
are pairwise disjoint. Consider $K>0$ with $\cup_{i=1}^M
B(\sigma_i,\varepsilon)\subset B(\eta,K)$.

i) $z\in \bigcup_{i=1}^M 
B(\sigma_i,\varepsilon)$. For $z\in 
B(\sigma_i,\varepsilon)$ we have $\rho(z)\simeq\rho(\eta)\simeq
\rho(\sigma_i)$,  $d_\phi(z,\Sigma)\simeq|z-\sigma_{i}|/\rho(\sigma_i)$ and
$d_\phi(z,\sigma_j)\gtrsim 1$, $j\neq i$. Thus 
\[
|P_\eta(z)|\lesssim
\left|\frac{h(z)}{z-\sigma_{i}} \right|  \rho(\eta) 
e^{-\varepsilon\phi(\eta)} \simeq e^{\varepsilon(\phi(z)-\phi(\eta))} .
\]

ii) $z\in B(\eta,K) \setminus\bigcup_{i=1}^M 
B(\sigma_i,\varepsilon)$. 
Here $\rho(z)\simeq\rho(\eta)$ and
$|z-\sigma_i|\gtrsim \rho(\eta)$, so
\[
|P_\eta (z)|\lesssim 
\frac{e^{\varepsilon\phi(z)} d_\phi(z,\Sigma) }{\rho^M(\eta)}
\frac{\rho^M(\eta)}{ e^{\varepsilon\phi(\eta)}}
\lesssim e^{\varepsilon(\phi(z)-\phi(\eta))} .
\]

iii) $z\notin B(\eta,K)$. Here $d_\phi(z,\sigma_i)\simeq 
d_\phi(z,\eta)$, so
\begin{equation*}
|P_\eta(z)|\lesssim \frac{e^{\varepsilon\phi(z)}
d_\phi(z,\Sigma)}{|z-\eta|^M}
\frac{\rho^M(\eta)}{e^{\varepsilon\phi(\eta)}}
\lesssim e^{\varepsilon(\phi(z)-\phi(\eta))}
\bigl(\frac{\rho(\eta)}{|z-\eta|}\bigr)^M  .
\end{equation*}
This and Lemma~\ref{distance}(b) solve the case $\omega=1$. 

For arbitrary $\omega\in\Omega_\phi$ there exists
$\gamma>1$  such that if $d_\phi(z,\eta)\geq 1$ then
\[ 
d_\phi^{-\gamma}(z,\eta)\lesssim \frac{\omega(\eta)}{\omega(z)}\lesssim
d_\phi^{\gamma}(z,\eta) .
\]
Thus the result follows from the previous construction 
taking $M$ big enough and using again
Lemma~\ref{distance}(b).
\end{proof}

\section{Basic properties of functions in $\Fpa$}\label{basics}

Here we study the behaviour of functions in $\Fpa$ and related topics. We
prove the estimates with norms $\|\cdot\|_{\Fpa}$ on the solutions to the
$\db$ equation (Theorem C) and provide estimates of the Bergman Kernel of
$\Fdosa$ on the diagonal. We also introduce a scaled translation  in the plane
that gives rise to a translated weight and to an isometry between the spaces
of functions for the original and the translated weight. This will be used
when studying the properties of weak limits (Section~\ref{weak-limits}).

\subsection{Pointwise estimates}
Let us first see what is the natural growth 
of functions in $\Fpa$. Recall that $d\sigma=dm/\rho^2$.

\begin{lemma}\label{estimacio-puntual} Let $1\leq p<\infty$ and
$\omega\in\Omega_\phi$. 
For any $r>0$ there exists $C>0$ 
such that for any $f\in H(\C)$ and $z\in\C$:
\begin{itemize}
\item[(a)] $|f(z)|^p e^{-p\phi(z)}\leq \displaystyle\frac
{C}{\omega^{p}(z)}\int_{D^r(z)}|f|^p e^{-p\phi} \omega^{p} d\sigma$.
\item[(b)] $|\nabla(|f|e^{-\phi})(z)|\leq
\displaystyle\frac {C}{\omega(z)\rho(z)}\bigl(\int_{D^r(z)}|f|^p
e^{-p\phi} \omega^p d\sigma \bigr)^{1/p}$.
\item[(c)] If $R>r$ then  $\ |f(z)|^p e^{-p\phi(z)}\leq
\displaystyle\frac {C_{R}}{\omega^{p}(z)}
\int_{D^R(z)\setminus D^r(z)}|f|^p e^{-p\phi} \omega^{p}d\sigma$.
\end{itemize}
\end{lemma}

\begin{proof} 
Let $H_z$ be a holomorphic function with $\Re H_z=h_z$, where $h_z$ is the
harmonic function in $D^r(z)$ given in Lemma \ref{cota-disc}.\\
(a) is proved as in \cite[Lemma 1]{OrSe}:
\begin{eqnarray*}
|f(z)|^p e^{-p\phi(z)}&=&|f(z)e^{-H_z(z)}|^p e^{-p\phi(z)}\\
&\leq& \frac{C}{\rho^2(z)} \int_{D^r(z)}|f(\zeta)|^pe^{-p(h_z(\zeta)
+\phi(z))}
\simeq  \frac{1}{\omega^{p}(z)} \int_{D^r(z)}|f|^p e^{-p\phi} 
\omega^{p} d\sigma.
\end{eqnarray*} 
(b) First let us see that  
$|\partial\phi/\partial \zeta- \partial h_z/\partial \zeta| 
\lesssim 1$ on $D^r(z)$. 
By (\ref{green}), if $\zeta\in D^r(z)$
\[
|\frac{\partial\phi}{\partial \zeta}(\zeta)-
\frac{\partial h_z}{\partial \zeta}(\zeta)|=
\bigl| \frac{\partial}{\partial \zeta} \int_{D^r(z)}
G(\zeta,\eta)\Delta\phi(\eta)\bigr|\leq
\int_{D^r(z)} \frac{2r\rho(z)}{| \zeta-\eta|} 
\Delta\phi(\eta) .
\]
Take $R$ (depending on $r$) such that 
$D^r(z)\subset D^R(\zeta)$.
From $\Delta\phi\simeq 1/\rho^2$ we deduce
\[
\int_{D^r(z)} \frac{2r\rho(z)}{| \zeta-\eta|} 
\Delta\phi(\eta) \lesssim
\frac 1{\rho(\zeta)}\int_{D^R(\zeta)} 
\frac{d m(\eta)}{|\zeta-\eta|}\simeq 1 .
\]
Since  $|\nabla(|f|e^{-\phi})|=|f'-2f \partial\phi/\partial z|e^{-\phi}$, 
we have 
\begin{equation}\label{loc1}
|\nabla(fe^{-H_z})(z)|=|f'(z)-2f(z) h_z'(z)|  \simeq 
|\nabla(|f|e^{-\phi})(z)|e^{\phi(z)}   .
\end{equation}
On the other hand,
\begin{eqnarray*}
|\nabla(fe^{-H_z})(z)| \lesssim \bigl| \int_{|z-\zeta|=\rho(z)}
\frac{f(\zeta)e^{-H_z(\zeta)}}{(z-\zeta)^2}\, d\zeta \bigl| \simeq 
\frac{1}{\rho^2(z)} \int_{|z-\zeta|=\rho(z)} |f(\zeta)|e^{-h_z(\zeta)}|d\zeta|.
\end{eqnarray*}
From (a), for $|z-\zeta|=\rho(z)$
\begin{equation*}
|f(\zeta)|e^{-\phi(\zeta)}\lesssim \frac1{\omega^{p}(z)} 
\bigl(\int_{D^r(z)}  |f|^p e^{-p\phi}\omega^{p}d\sigma  \bigr)^{1/p} .
\end{equation*}
By Lemma \ref{cota-disc} we have then
 \begin{equation*}
|\nabla(fe^{-H_z})(z)| \lesssim  \frac1{\omega(z)\rho(z)} 
\bigl(\int_{D^r(z)}  |f|^p e^{-p\phi} \omega^{p}
d\sigma\bigr)^{1/p} e^{\phi(z)} ,
\end{equation*}
which together with (\ref{loc1}) concludes the proof.

(c) As (a), using the subharmonicity of $|fe^{-H_z}|^p$.
\end{proof}

\begin{lemma}\label{derivative}
Let $1\leq p<\infty$ and $\omega\in\Omega_\phi$.
For any entire function 
$g$ with $g(\lambda)=0$ we have
\begin{equation*}
|g'(\lambda)| e^{-\phi(\lambda)}
\lesssim \frac1{\omega(\lambda)\rho(\lambda)}
\left(\int_{D(\lambda) }|g|^p e^{-p\phi} \omega^p d\sigma\right)^{1/p}.
\end{equation*}
\end{lemma}

\begin{proof}
Lemma \ref{estimacio-puntual}(c) with $r=1/2$ and $R=1$ applied to the function 
$g(z)/(z-\lambda)$ yields
\begin{eqnarray*}
|g'(\lambda)|^p e^{-p\phi(\lambda)}& \lesssim &\frac 1{\omega^p(\lambda)}
\int_{D(\lambda) \setminus D^{1/2}(\lambda) }
\frac{|g(z)|^p}{|z-\lambda|^p} e^{-p\phi(z)}\omega^p d\sigma\\
&\lesssim&
\frac 1{\omega^p(\lambda)\rho^{p}(\lambda)} 
\int_{D(\lambda)} |g(z)|^p e^{-p\phi(z)}
\omega^p d\sigma .
\end{eqnarray*}
\end{proof}

\subsection{H\"ormander type estimates}
This section is devoted to the proof of the $\db$-estimates of Theorem C in
the introduction.

\begin{TheoremC}
Let $\phi$ be a subharmonic function
such that $\Delta\phi$ is a doubling measure. For any $\omega\in \Omega_\phi$, 
there is a solution $u$ to the equation $\db u = f$ such that $\| u
e^{-\phi}\omega\|_{L^p(\C)}  \lesssim \| f e^{-\phi}\omega\rho
\|_{L^p(\C)}$ for any $1\le p\le \infty$.  
\end{TheoremC}

\begin{proof}
By Lemma~\ref{estimacio-puntual}(b), there exists $r>0$ such that 
$|P_\eta(z)| \gtrsim e^{\varepsilon({\phi(z)-\phi(\eta)})}$ on  $D^r(\eta)$,
for all $\eta \in \C$. Take a sequence $\Lambda$ such that
$\{D^r(\lambda)\}_{\lambda\in\Lambda}$ covers $\C$ and the disks 
$\{D^{r/5}(\lambda)\}_{\lambda\in\Lambda}$ are pairwise disjoint, which exist
by a standard covering Lemma, see \cite[Theorem~2.1]{Mattila}. 
Let $\{\chi_\lambda\}\subset
\mathcal C_0^\infty$ be a partition of unity associated to
$\{D^r(\lambda)\}_\lambda$. 

Decompose the datum $f=\sum
f_\lambda$, with $f_\lambda(z)= f(z) \chi_{\lambda}(z)$. 
By Theorem~\ref{weight}, for any $\lambda$ there exists 
an entire function $m_\lambda(z)=P_\lambda(z) e^{-\phi(\lambda)}$ such that 
\[
|m_\lambda(z)| \lesssim e^{\phi(z)} \frac 1{d_\phi^{M}(z,\lambda)+1}
\frac{\omega(\lambda)}{\omega(z)}.
\]
The radius $r$ has been chosen so that 
$|m_\lambda(\zeta)|\gtrsim e^{\phi(\zeta)}$ if $\zeta\in D^r(\lambda)$.
Define 
\[
u_\lambda(z)= m_\lambda(z) \frac 1\pi \int_{D^r(\lambda)}
\frac{f_\lambda(\zeta)/m_\lambda(\zeta)}{\zeta -z}\, dm(\zeta).
\]
Clearly $\db u_\lambda = f_\lambda$, thus  $u=\sum_{\lambda\in\Lambda}
u_\lambda$ is as a solution to $\db u = f$.  We must prove the size estimates.
As we have used a linear operator to construct $u$ from the datum  $f$, we
only need to check that $\| u e^{-\phi}\omega\|_{L^\infty} \lesssim \| f
e^{-\phi}\omega\rho  \|_{L^\infty}$ and $\| u e^{-\phi}\omega\|_{L^1} \lesssim
\| f e^{-\phi} \omega \rho \|_{L^1}$. The estimates for $1<p<\infty$ follow
then by Marcinkiewicz interpolation theorem.

Assume that $z\in D^r(\lambda)$ and take $K,K'>0$ such that $D^K(z)\subset
D^{K'}(\lambda)$. Then 
\begin{eqnarray*}
|u_\lambda(z) e^{-\phi(z)}\omega(z)| &\lesssim& \int_{D^K(z)}
\frac{|f(\zeta)|e^{-\phi(\zeta)}\omega(\zeta)\rho(\zeta)}
{\rho(z)|\zeta-z|}\, dm(\zeta)\\
&\lesssim&
\int_{D^{K'}(\lambda)}
\frac{|f(\zeta)|e^{-\phi(\zeta)}\omega(\zeta)\rho(\zeta)}
{\rho(\lambda)|\zeta-z|}\, dm(\zeta).
\end{eqnarray*}
On the other hand, if $z\notin D^{K}(\lambda)$
\begin{eqnarray*}
|u_\lambda(z) e^{-\phi(z)}\omega(z)| &\lesssim &
d_\phi^{-M}(z,\lambda) \int_{D^{r}(\lambda)}
\frac{|f(\zeta)| e^{-\phi(\zeta)}}{|\zeta-z|}\omega(\zeta) \, dm(\zeta)\\
&\lesssim & \frac{d_\phi^{-M}(z,\lambda)}{\rho^2(\lambda)}
 \int_{D^{r}(\lambda)}
|f(\zeta)| e^{-\phi(\zeta)}\omega(\zeta)\rho(\zeta) \, d m(\zeta) .
\end{eqnarray*}
Therefore, applying Lemma~\ref{suma} 
\[
\|u e^{-\phi}\omega\|_{L^\infty} \lesssim  \|f e^{-\phi}\omega \rho 
\|_{L^\infty} \sup_{z\in \C}
\Bigl( \int_{D^K(z)} 
\frac{d m(\zeta)}{\rho(z)|z-\zeta|} + \sum_{\lambda:z\notin D^{r}(\lambda)}
d_\phi^{-M}(z,\lambda) 
\Bigr)\lesssim 
 \|f e^{-\phi}\omega \rho \|_{L^\infty}.
\]
In the $L^1$ norm we get
\[
\begin{split}
\|u e^{-\phi}\omega\|_{L^1}\lesssim \sum_{\lambda\in \Lambda}
&\Bigl( 
\int_{z\in D^{r}(\lambda)} \int_{\zeta\in D^{K'}(\lambda)}
\frac{|f(\zeta)| e^{-\phi(\zeta)}\omega(\zeta) \rho (\zeta)}{\rho(\lambda)
|\zeta-z|}\, 
dm(\zeta)\, d m(z)+\\
&\int_{z\notin D^{r}(\lambda)}\frac{d_\phi^{-M}(z,\lambda)}{\rho(\lambda)^{2}}
\int_{D^{r}(\lambda)}
|f(\zeta)| e^{-\phi(\zeta)}\omega(\zeta)\rho(\zeta) \, dm(\zeta)\, d m(z)
\Bigr).
\end{split}
\]
Reversing the order of integration we immediately get 
$\|u e^{-\phi}\omega\|_{L^1}\lesssim 
\|f e^{-\phi}\omega \rho\|_{L^1}$.
\end{proof}

\subsection{Bergman kernel estimates}\label{bergman} 
Let $K_{\phi,\omega}(z,\zeta)$ denote the
Bergman kernel for $\Fdosa$, i.e, for any $f\in\Fdosa$
\[
f(z)=\langle   K_{\phi,\omega}(z,\cdot),f \rangle = 
\iC K_{\phi,\omega}(z,\zeta) f(\zeta) e^{-2\phi(\zeta)} 
\omega^2(\zeta)  d\sigma(\zeta) .
\]
By definition
\[
K_{\phi,\omega}(z,z)=\int_{\C} |K_{\phi,\omega}(z,\zeta)|^2 e^{-2\phi(\zeta)} 
\omega^2(\zeta) d\sigma(\zeta) .
\]

\begin{lemma}\label{bergman-diagonal}
There exists $C>0$ such that 
\[
C^{-1} (e^{\phi(z)}/\omega(z) )^2\leq 
K_{\phi,\omega}(z,z)\leq C
 ( e^{\phi(z)}/\omega(z) )^2\qquad z\in\C .
\]
\end{lemma}

\begin{proof}
We use the identity
\[
\sqrt{K_{\phi,\omega}(z,z)}=\sup\{|f(z)| : f\in \Fdosa 
\; ,\; \|f\|_{\Fdosa}
\leq 1\} .
\]
The estimate $\sqrt{K_{\phi,\omega}(z,z)}\lesssim e^{\phi(z)}/\omega(z)$ 
is immediate
from Lemma \ref{estimacio-puntual}(a).
In order to prove the reverse estimate we construct $f\in \Fdosa$
with $\|f\|_{\Fdosa}\leq 1$ and
$|f(z)|\geq C e^{\phi(z)}/\omega(z) $, 
for some constant $C$ independent of $z$.
 
By Theorem~\ref{weight}, for every $m\in\N$ there exists $P_z$ entire such that
\[
|P_z(\zeta)|\leq C
e^{\phi(\zeta)-\phi(z)} \frac{\omega(z)}{\omega(\zeta)} 
\frac{1}{1+d_\phi^m(z,\zeta)} ,
\]
with $C$ independent of $z$. Define
$f_z(\zeta)=c_0\; e^{\phi(z)}/\omega(z)\; P_z(\zeta)$, where $c_0$ is a 
positive
constant to be chosen later. Now
$f_z(z)=c_0 e^{\phi(z)}/\omega(z)$ and
\[
|f_z(\zeta)|^2 e^{-2\phi(\zeta)} \omega^2(\zeta)\rho^{-2}(\zeta)
\leq  
\frac{c_0\; C}{1+d_\phi^m(z,\zeta)}\Delta\phi(\zeta) ,
\]
hence by Lemma \ref{integrals}(b) there exist $c_0$ and $C$ independent of $z$
so that $\|f_z\|_{\Fdosa}\leq 1$.
\end{proof} 
\begin{remark}\label{general-bergman}
This argument and Lemma~\ref{estimacio-puntual}(a)
show that for any $p\in [1,\infty]$
\[
\sup\{|f(z)| : f\in \Fpa \; ,\; \|f\|_{\Fpa}
\leq 1\}\simeq e^{\phi(z)}/\omega(z)  .
\]
\end{remark}

\subsection{Scaled translations and invariance}\label{invariance}
In this section we introduce the scaled translation and its main properties.

Given $\phi$ consider the class $W_\phi$ of subharmonic functions 
$\psi$ such that
\begin{itemize}
\item[(i)] $\Delta\psi$ doubling with 
$C_{\Delta\psi}\leq C_{\Delta\phi}$.
\item[(ii)] $\int_{D_\phi(0)}\Delta\psi\simeq 1$.
\item[(iii)] $\psi(0)=0$.
\end{itemize}

An important property of $W_\phi$ is that there exists $\eta$ such that
$\Delta\psi(z) \lesssim |z|^{2\eta}$ for all regular
$\psi\in W_\phi$. This is a consequence
of (\ref{cota-rho}) and the fact that $\Delta\psi \simeq 1/\rho_\psi^2$.

Fix $q>2\eta+1$ and consider the kernel
\begin{equation*}
\kappa (z,\zeta):=\frac 1{2\pi}\left[\log|1-\frac{z}{\zeta}| - 
\Re(P_q(\frac{z}{\zeta}))\chi_{\C\setminus D(0,1)}(\zeta)\right],
\end{equation*}
where $P_q$ is the Taylor polynomial of degree $q$ of $\log(1+x)$ around $x=0$,
and its associated integral operator
\[
K[f](z)=\int_{\C} \kappa (z,\zeta) f(\zeta)\; dm(\zeta) .
\]
This operator solves the Poisson equation, that is $\Delta K[f]=f$.

For every $x \in \C$, consider the scaled translation 
\[
\tau_x(z) = x+z\rho_{\phi}(x),
\] 
the associated subharmonic function
\begin{equation*}
\phi_x(z) = 
K[\Delta(\phi\circ\tau_x)](z)- K[\Delta(\phi\circ\tau_x)](0),
\end{equation*}
and the associated weight
\[
\omega_x(z)=\omega(\tau_x(z))/\omega(x).
\]
Define also $\mathsf h_x:= \phi\circ\tau_x - \phi_x$. It is clear that 
$\mathsf h_x$ is harmonic. Take then $H_x$ holomorphic having 
$\mathsf h_x$ as real part and consider the scaled translation operator 
\begin{equation*}
T^{\phi,\omega}_x f(z) =f(\tau_x(z))e^{-H_x(z)}\omega (x).
\end{equation*}

\begin{lemma}\label{translations}
For every $x\in\C$,
\begin{itemize}
\item[(a)] The subharmonic function $\phi_x $ belongs to $ W_\phi$, 
and the weight $\omega_x$
satisfies $\omega_x(0)=1$ and $\omega_x\in \Omega_{\phi_x}$.
\item[(b)] $T_x^{\phi,\omega}$ is an isometry from $\Fpa$ to 
$\mathcal{F}_{\phi_x,\omega_x}^{p}$, for $1\leq p\leq \infty$.
\end{itemize}
\end{lemma}

\begin{proof}
Note first that from the identity
\begin{equation*}
1=\int_{D_{\phi_x}(z)}\Delta\phi_x=\int_{D_{\phi_x}(z)}\rho^2_\phi(x)
\Delta\phi(\tau_x(\zeta))=
\int_{D(\tau_x(z),\rho_{\phi_x}(z)\rho_{\phi}(x))}\Delta\phi
\end{equation*}
it follows that
\begin{equation}\label{identite_rho_x}
\rho_\phi(\tau_x(z)) = \rho_{\phi_x}(z) \rho_\phi(x).
\end{equation}
This implies that the mapping $\tau_x$ is actually an isometry between $\C$
endowed with the distance $d_{\phi_x}$ and $\C$ with $d_{\phi}$, 
that is
\begin{equation}\label{isometria}
d_{\phi_x}(z,\zeta)=d_{\phi}(\tau_x(z),\tau_x(\zeta))\quad \forall z,\zeta \in
\C.
\end{equation}
(a) By definition  $\phi_x(0)=0$, and by \eqref{identite_rho_x}
$\rho_{\phi_x}(0)=1$. This gives properties (ii) and (iii) of $W_\phi$.

It is also clear that $\Delta \phi_x$ is doubling and
$C_{\Delta\phi_x}=C_{\Delta\phi}$, since for any $a\in\C$ and $r>0$
\begin{multline*}
\int_{D(a,2r)}\Delta\phi_x= \int_{D(\tau_x(a),2r\rho_{\phi}(x) )}
\Delta\phi\leq C_{\Delta\phi} \int_{D(\tau_x(a),r\rho_{\phi}(x) )} 
\Delta\phi\leq C_{\Delta\phi}  \int_{D(a,r)}\Delta\phi_x.
\end{multline*}
That $\omega_x(0)=1$ and
$\omega_x(\zeta)/\omega_x(z)=\omega(\tau_x(z))/\omega(\tau_x(\zeta))$
follows from the definition. 
This and \eqref{isometria} imply that $\omega_x\in \Omega_{\phi_x}$ . 

(b) For $p<\infty$ we use the change of variable $\zeta=\tau_x(z)$ and
\eqref{identite_rho_x}:
\begin{multline*}
\int_{\C}|T_x^{\phi,\omega}(f)(z)|^p e^{-p\phi_x(z)} 
\omega_x^p(z) \frac{ dm(z)}{\rho_{\phi_x}^{2}(z)}\\
=\int_{\C} |f(\tau_x(z))|^p e^{-p\phi( \tau_x(z))}
\omega^{p}(x) 
\bigl(\frac{\omega(\tau_x(z))}{\omega(x)}\bigr)^{p}
\bigl(\frac{\rho_\phi(\tau_x(z))}{\rho_\phi(x)}\bigr)^{-2}
d m(z)\\
=\int_{\C} |f(\zeta)|^p e^{-p\phi(\zeta)}\omega^p(\zeta) 
\frac{d m(\zeta)}{\rho_\phi^{2}(\zeta)} .
\end{multline*}
The case $p=\infty$ is straightforward from \eqref{identite_rho_x}.
\end{proof}

Given a sequence $\Lambda$ and $x\in \C$  let 
\[
\Lambda_x:=(\tau_x)^{-1}(\Lambda) .
\]
Given a sequence $\Lambda$ and $z\in\C$, denote
$n_\Lambda(z,r)=\#(\Lambda\cap \overline{D(z,r)})$, for any $r>0$.
\begin{lemma}\label{properties} Let $\Lambda$ be a sequence in $\C$.
\begin{itemize}
\item[(a)] $\Lambda$ is  $\rho$-separated if and only if 
$\Lambda_x$ is  $\rho_{\phi_x}$-separated.
\item[(b)] $\Lambda\in\Int \Fpa$  
if and only if  $\Lambda_x\in\Int \mathcal{F}_{\phi_x,\omega_x}^{p}$. 
Similarly,
$\Lambda\in\Samp \Fpa$  
if and only if  $\Lambda_x\in\Samp \mathcal{F}_{\phi_x,\omega_x}^{p}$.
Furthermore, the interpolation and sampling
constants remain the same.
\item[(c)] The densities are stable: 
$\mathcal D_{\Delta\phi}^+(\Lambda)=\mathcal D_{\Delta\phi_x}^+(\Lambda_x)$,
and  
$\mathcal D_{\Delta\phi}^-(\Lambda)=\mathcal D_{\Delta\phi_x}^-(\Lambda_x)$.
\end{itemize}
\end{lemma}

\begin{proof}
(a) is an immediate consequence of (\ref{identite_rho_x}).

(b) is a consequence of Lemma \ref{translations} and the identity 
$\|f|\Lambda\|_{\ellpa(\Lambda)}=
\|T_x^{\phi,\omega} f|\Lambda_x\|_{\ell_{\phi_x,\omega_x}^{p}(\Lambda_x)}$.

(c) 
Define
\begin{equation}\label{densite}
\mathcal D_{\Delta\phi}(z,r,\Lambda) = \displaystyle 
\frac{n_\Lambda(z,r\rho(z))}{\int_{D_{\phi}^r(z)}\Delta\phi} .
\end{equation}
By a change of variables, it is clear that
\begin{equation*}
\mathcal D_{\Delta\phi}(z,r,\Lambda)=
\mathcal D_{\Delta\phi_x}((\tau_x)^{-1}(z),r,\Lambda_x).
\end{equation*}
Taking the supremum over $z \in \C$ and passing to the
limsup  we get the result for the upper density. The lower density is 
dealt with similarly.
\end{proof}

\subsection{Weak limits.}\label{weak-limits}
In this section we study weak limits of sequences 
$\Lambda$ and their properties. 
 
\begin{definition}
A sequence of closed sets $Q_j$ \emph{converges strongly} to $Q$, denoted
$Q_j \rightarrow Q$ if $[Q,Q_j]  \rightarrow 0$; here $[Q,R]$ denotes the
Fr\'echet distance between $Q$ and $R$. We say that $Q_j$ 
\emph{converges compactwise}  to $Q$, denoted $Q_j \rightharpoonup Q$, if for
every compact set $K$ we have $(Q_j\cap K) \cup \partial K \rightarrow (Q\cap
K )\cup \partial K$.
\end{definition}

\begin{definition}
A set 
$\Lambda^*$ is a \emph{weak limit} of $\Lambda$ if there exists
a sequence $\{x_n\}_{n\in \N}$ in $\C$ such that $\Lambda_{x_n} \rightharpoonup
\Lambda^*$. 
\end{definition}
Given a $\rho$-separated sequence $\Lambda$, and a 
sequence $\{x_n\}_{n\in \N}$ it is always possible to extract a
subsequence of $\Lambda_{x_{n_j}}$ such that $\Lambda_{x_{n_j}}\rightharpoonup
\Lambda^*$ for some $\Lambda^*$. We need also a normal family argument for
the translated weights that define the space.

\begin{lemma}\label{weak_convergence}
Let $\{x_n\}_{n\in\N}$ be a sequence in $\C$. There exist a subharmonic 
function $\phi^*$, a weight $\omega^*\in\Omega_{\phi^*}$  and a
subsequence  $\{x_{n_k}\}_k$ such that  $\{\phi_{x_{n_k}}\}_k$,
$\{\omega_{x_{n_k}}\}_k$ and  $\{\Delta\phi_{{x_n}_k}\}_k$
converge uniformly on
compact sets to $\phi^*$, $\omega^*$ and $\Delta\phi^*$ respectively. 
Furthermore, $\Delta\phi^*$ is a doubling measure and $C_{\Delta\phi^*}
\le C_{\Delta\phi}$.
\end{lemma}

\begin{proof} 
Take $\eta$ and $q>2\eta+1$ as in the definition of the kernel $\kappa$ (see
previous section). Denote $\mu_n=\Delta\phi_{x_n}$. 

Since $|\nabla \mu_n| \lesssim \rho_{\phi_{x_n}}^{-3}$
(Theorem~\ref{regularitzacio}) and  $\rho_{\phi_{x_n}}(0)=1$, for any compact
set $K$ there exits $C_K>0$ such that  $|\nabla \mu_n(z) |\leq C_K$. By the
Arzel\`a-Ascoli theorem, we can extract a subsequence $\{\mu_{n_k}\}_{k}$
converging  uniformly on compact sets of $\C$ to a function $\mu^*$. It
follows immediately that the measure with density $\mu^*$ is doubling and
$C_{\mu^*}\le C_{\mu_n}=C_{\Delta\phi}$. Furthermore, this implies that
$\rho_{\phi_{x_n}}\to \rho^*$ uniformly on compacts.

Let now $\phi^*= K[\mu^*]- K[\mu^*](0) $, and denote 
$\phi_k:=\phi_{x_{n_k}}$, $\mu_k:=\mu_{n_k}$. We will show that 
$\{\phi_k\}_k$ converges uniformly on compact sets to $\phi^*$. 

By definition $\phi_p(z)= K[\mu_p](z)- K[\mu_p](0)$, thus we only have to
prove that $ K[\mu_p]$ converges uniformly on compacts set to $ K[\mu^*]$.
Take $z \in D(0,R)$ and $t>R$. Then
\begin{multline*}\label{ineq_weak}
| K[\mu_p](z) -  K[\mu^*](z)| \leq  \Bigl|\int_{\C\setminus D(0,t)}  
\kappa(z,\zeta) (\mu_p(\zeta) - \mu^*(\zeta))dm(\zeta)  \Bigr| \\
 +   \Bigl|\int_{D(0,t)}  \kappa (z,\zeta) (\mu_p(\zeta) - \mu^*(\zeta))
 dm(\zeta) \Bigr| .
\end{multline*}
Let $I_1$ be the first integral. By construction of $\kappa$ we have
\begin{equation*}
|\kappa (z,\zeta)| \lesssim \Bigl(\frac{R}{|\zeta|}\Bigr)^q.
\end{equation*}
Also, $(|\mu_p(\zeta)| + |\mu^*(\zeta)|)dm(\zeta)$ is a doubling 
measure with doubling  constant less than $C_{\Delta\phi}$. By (\ref{cota-rho}) 
$|\mu_p(\zeta)| + |\mu^*(\zeta)|\lesssim |\zeta|^{2\eta}$, and therefore
\begin{equation*}
I_1 \lesssim \int_{|\zeta|>t} \Bigl(\frac{R}{|\zeta|}\Bigr)^q |\zeta|^{2\eta}
dm(\zeta).
\end{equation*}
This is smaller than $\varepsilon$ for $t$ big enough.

Let $I_2$ be the second integral in the estimate above. We have
\begin{multline*}\label{ineq_weak_2}
I_2 \lesssim \int_{D(0,1)} \Bigl|\log\bigl|\frac{z-\zeta}{\zeta}\bigr| \Bigr|
|\mu_p(\zeta) - \mu^*(\zeta)|dm(\zeta) + 
\int_{D(0,t)\setminus D(0,1) }|P_q(\frac{z}{\zeta})\|\mu_p(\zeta) - 
\mu^*(\zeta)|dm(\zeta)
\end{multline*}
For all $z \in D(0,R)$ and $\zeta \in D(0,t)\setminus D(0,1)$ we have
$|P_q( z/\zeta)| \leq C(R,t)$, hence the uniform convergence of 
$\mu_p$ implies that for $p$ big enough the second integral here 
is smaller 
than $\varepsilon$.  It remains to prove the convergence
of the first term. Take $C(t)$ such that$ \int_{D(0,t)}
|\log| z-\zeta / \zeta\|dm(\zeta) \leq C(t)$ and choose $p$ big
enough so that $ |\mu_p(\zeta) - \mu^*(\zeta)| \leq \varepsilon /C(t) $
uniformly on $D(0,t)$. Then the estimate follows. 

We know that the sequence of distance functions $d_{\phi_{x_n}}$ has a
subsequence converging to $d_{\phi^*}$ uniformly on compact sets of
$\C\times\C$, because the $\rho_{x_{n_k}}$ converge uniformly. By construction
$\omega_{x_n}(0)=1$. On the other hand, the definition of flat weight implies
that they are equibounded on any compact. Moreover, the regularity given by
\eqref{fw2} makes them  equicontinous on compact sets. We can thus extract
again a convergent subsequence. 
\end{proof}

\begin{corollary}\label{sequence-weak-limit} 
Given a subharmonic function $\phi$ with doubling Laplacian, $\Lambda$ a
$\rho$-separated  subsequence, $\omega\in \Omega_\phi$ and $\{z_n\}_{n\in\N }$ 
a sequence of complex numbers,  there exist a subharmonic function $\phi^*$, a
$\rho_{\phi^*}$-separated sequence $\Lambda^*$, a weight $\omega^*\in 
\Omega_{\phi^*}$ and a subsequence  $\{x_n\}_{n\in\N }$ of $\{z_n\}_{n\in\N }$ 
such that $\Lambda_{x_n} \rightharpoonup \Lambda^*$,  and $\phi_{{x_n}}
\rightarrow \phi^*$, $\omega_{x_n}\rightarrow \omega^*$  and
$\Delta\phi_{{x_n}}  \rightarrow \Delta\phi^*$ uniformly on compact sets.  
\end{corollary}

We will write $(\Lambda_{x_n},\phi_{{x_n}},\omega_{x_n})\rightarrow 
(\Lambda^*,\phi^*,\omega^*)$. The set of all such  weak limits will be denoted
by $W(\Lambda,\phi,\omega)$. 

Let us prove now the stability of the upper and lower densities with respect
to weak limits.

\begin{lemma} Let $\Lambda$ be a $\rho$-separated sequence, 
$\{x_n\}_n\subset \C$, and assume that $(\Lambda_{x_n},\phi_{x_n},\omega_{x_n})
\rightarrow (\Lambda^*,\phi^*,\omega^*)$. Then
\begin{itemize}
\item[(a)] $\mathcal D^{+}_{\Delta\phi}(\Lambda)<1/2\pi$  implies 
$\mathcal D^{+}_{\Delta\phi^*}(\Lambda^*)<1/2\pi$.
\item[(b)] $\mathcal D^{-}_{\Delta\phi}(\Lambda)>1/2\pi$ implies 
$\mathcal D^{-}_{\Delta\phi^*}(\Lambda^*)>1/2\pi$.
\end{itemize}
\end{lemma}

\begin{proof}
Denote $\Lambda_n=\Lambda_{x_n},\phi_n=\phi_{x_n}$
and $\rho_n=\rho_{x_n}$. By hypothesis 
$\{\Delta\phi_n\}_n\rightarrow \Delta\phi^*$ uniformly on compact sets, 
and therefore $\{\rho_n\}_n\rightarrow \rho^*$ also uniformly on 
compact sets. Thus, for any $\epsilon(r)>0$,
\[
\begin{split}
\frac{n_{\Lambda^*}(z,(r-\epsilon(r))\rho_{\phi^*}(z))}
{\int_{D_{\phi^*}^r(z)}\Delta\phi^*}\le
\liminf_{n\to\infty}
\frac{n_{\Lambda_n}(z,r\rho_n(z))}{\int_{D_{\phi_n}^r(z)} \Delta\phi_n}\le \\
\leq\limsup_{n\to\infty}
\frac{n_{\Lambda_n}(z,r\rho_n(z))}{\int_{D_{\phi_n}^r(z)} \Delta\phi_n}\le
\frac{n_{\Lambda^*}(z,(r+\epsilon(r))\rho_{\phi^*}(z))}
{\int_{D_{\phi^*}^r(z)} 
\Delta\phi^*} .
\end{split}
\]
(a) Since $\mathcal D_{\Delta\phi}^+(\Lambda)<1/2\pi$, there
exist $\varepsilon,R_0>0$ such that, if
$w= \tau_{x_n}^{-1}(z)$
\[
\frac{n_{\Lambda_n}(w,r\rho_n(w))}{\int_{D_{\phi_n}^r(w)} \Delta\phi_n}=
\frac{n_{\Lambda}(z,r\rho(z))}{\int_{D_{\phi}^r(z)} \Delta\phi}\leq 
1/2\pi-\varepsilon
\qquad \forall r>R_0\; ,\; \forall n\in\N\; ,
\; \forall w\in\C .
\]
Taking limits as $n\to\infty$ and picking $\epsilon(r)$ so 
that $\epsilon(r)/r\to 0$  we see, using Lemma~\ref{corona-zero}, that
$\mathcal D_{\Delta\phi^*}^+(\Lambda^*)<1/2\pi$.

(b) is proved similarly.
\end{proof}

\section{Preliminary properties of sampling and interpolating
sequences}\label{basic-int-samp} 

This section is devoted to prove auxiliary results on interpolating and
sampling sequences. A main result is that there do not exist
sequences which are
simultaneously sampling and interpolating.  We also prove some results on
inclusions between spaces of sampling and interpolating sequences for various
weights.

An easy consequence of Lemma \ref{estimacio-puntual} is that we only need to
deal with $\rho$-separated sequences.

\begin{lemma}\label{separacio} Let 
$\Lambda\subset\C$.
\begin{itemize}
\item[(a)] If $\Lambda\in \Int  \Fpa$, 
then $\Lambda$ is $\rho$-separated.
\item[(b)] If $\Lambda\in \Samp  \Fpa$, there exists a 
$\rho$-separated subsequence
$\Lambda'\subset\Lambda$ such that $\Lambda'\in \Samp  \Fpa$.
\item[(c)] If $p<\infty$ and $\Lambda \in \Samp \Fpa$, then $\Lambda$ is a 
finite union of $\rho$-separated sequences.
\item[(d)] Let $\Lambda\in\Samp\Fpa$ be $\rho$-separated. There exists 
$r>0$ such 
that $\C=\cup_{\lambda\in\Lambda} D^r(\lambda)$.
\end{itemize}
\end{lemma}

\begin{proof}
(a) Assume that $\lambda,\mu\in\Lambda$ with $|\lambda-\mu|\leq
\rho(\lambda)$ and take $f\in \Fpa$ such that
$f(\lambda)= 
e^{\phi(\lambda)}/\omega(\lambda)$, $f(\mu)=0$ and $\|f\|_{\Fpa}\lesssim 1$.
Then
\[
\omega^{-1}(\lambda)=\left| |f(\lambda)|
e^{\phi(\lambda)}-|f(\mu)|
e^{-\phi(\mu)}\right|\lesssim |\nabla(|f| e^{-\phi})(\zeta)|
|\mu-\lambda|  .
\]
The result follows then from Lemma \ref{estimacio-puntual}(b).

(b) As in the proof of \cite[Theorem 2, p. 344]{Br}, using here
Lemma~\ref{estimacio-puntual}(b) instead of Bernstein's theorem, we get
\begin{equation*}
\left|\frac1{L^{p}_{\phi,\omega}(\Lambda)} - 
\frac1{L^{p}_{\phi,\omega}(\Lambda')} \right| \leq C[\Lambda,\Lambda'].
\end{equation*}

(c) It is enough to show that there exists $r>0$ and $M$ such that 
$\#(D^\eta(z)\cap \Lambda)\leq M$ for all $z \in \C$. 
To this end, consider the function 
$f_z(\zeta)=e^\phi(z)/\omega(z) P_z(\zeta)$, where 
$P_z$ is given by Theorem~\ref{weight} (with $\varepsilon=1$ and $\omega=0$). 
We have $\|f_z\|_{\Fpa} \leq C$, and for $r$ small enough 
$|f_z(\zeta)| \gtrsim  e^\phi(\zeta)/\omega(\zeta)$ in $D^r(z)$.
So the left sampling inequality (see (\ref{sampling1})) yields 
\begin{equation*}
\#(D^r(z)\cap \Lambda) \leq \|f_z|\Lambda\|_{\ellpa(\Lambda)} 
\leq CL_{\phi,\omega}^{p}(\Lambda).
\end{equation*}

(d) It is enough to see that for $R$ big enough $\Lambda\cap D^R(z)
\neq\emptyset$ for all $z\in\C$.

Take
$f_z$ as in (c). Let $\varepsilon>0$ be the $\rho$-separation of $\Lambda$.
Since
\[
|f_z(\zeta)|^p e^{-p\phi(\zeta)}\omega^p(\zeta) \rho^{-2}(\zeta)
\lesssim \frac{\Delta\phi(\zeta)}{1+d_\phi^m(z,\zeta)},
\]
Lemma~\ref{estimacio-puntual}(a) and Lemma~\ref{flors} lead to
\begin{eqnarray*}
\sum_{\lambda\notin D^R(z)} \omega^{p}(\lambda)|f_z(\lambda)|^p 
e^{-p\phi(\lambda)}
&\lesssim & \sum_{\lambda\notin D^R(z)}\int_{D^\varepsilon(\lambda)}
\frac{\Delta\phi(\zeta)}{1+d_\phi^m(z,\zeta)}\\
&\lesssim & \int_{\zeta\notin D^{R-\epsilon(R)}(z)}
\frac{\Delta\phi(\zeta)}{1+d_\phi^m(z,\zeta)}\\
\end{eqnarray*}

According to Remark~\ref{int-to-zero} this tends to 0 uniformly in $z$ 
as $R$ goes to
$\infty$. Thus, for $R$ big enough the sampling inequality gives 
\[
1\leq C \sum_{\lambda\in\Lambda\cap D^R(z)}
\omega^p(\lambda)|f_z(\lambda)|^p e^{-p\phi(\lambda)} .
\]
In particular $\Lambda\cap D^R(z)\neq \emptyset$, as desired.
\end{proof}

\subsection{Weak limits and interpolating and sampling sequences}
In this section $\tau_x^\phi$ will denote the scaled translation associated to
the weight $\phi$, as described in Section~\ref{invariance}. The main result is
as follows.

\begin{proposition} \label{weak_stability} Let $\phi$ a 
subharmonic function with doubling Laplacian, 
$\omega\in \Omega_{\phi}$ and 
$\Lambda$ be a $\rho$-separated sequence. Assume  
$(\Lambda^*,\phi^*,\omega^*)\in W(\Lambda,\phi,\omega)$.
\item[ (a)] If $\Lambda \in \Samp \Fpa $ then 
$\Lambda^* \in\Samp\calf_{\phi^*,\omega^*}^{p}$.
\item[ (b)] If $\Lambda \in \Int \Fpa $ then 
$\Lambda^* \in\Int\calf_{\phi^*,\omega^*}^{p}$.
\end{proposition}
\begin{proof}
(a) We argue by contradiction. Otherwise there exist 
$\varepsilon_n>0$ decreasing to zero and functions 
$f_n \in \calf_{\phi^*,\omega^*}^{p}$
such that $\|f_n\|_{\calf_{\phi^*,\omega^*}^{p}}=1$ and 
$\| f_n|\Lambda^* \|_{\ell_{\phi^*,\omega^*}^{p}(\Lambda^*)}  
\leq \varepsilon_n$.

By Corollary \ref{sequence-weak-limit} there exists a
sequence $\{x_j\}_{j\in \N}$ in $\C$ such that
$(\phi_j,\omega_j,\Lambda_j)\rightarrow(\phi^*,\omega^*,\Lambda^*)$, where
we denote $\Lambda_j:=\Lambda_{x_j}$, $\omega_j:=\omega_{x_j}$ 
and $\phi_j:=\phi_{x_j}$. 

For every $n$ consider $R_n$ big enough so that if 
$D_n:=D_{\phi^*}^{R_n}(0)$ then $\|f_n|D_n\|_{\calf_{\phi^*,\omega^*}^{p}} 
\geq 1-\varepsilon_n$. Set $\widetilde{D}_n:=D_{\phi^*}^{R_n^2}(0)$.

We claim that there exists a smooth cut-off function $\mathcal X_n$ such that
$\mathcal X_n(\zeta)=1$  in $ D_n$,  
$\mathcal X_n(\zeta)=0$ in $\C\setminus \widetilde{D}_n$ 
and $|\displaystyle \bar\partial \mathcal X_n | \leq 
\varepsilon_n / \rho_{\phi^*} $. 
To see this start with a smooth 
$\mathcal X_n$  depending linearly on $|\zeta|$ on 
$R_n\leq |\zeta|\leq R_n^2$. Then 
\begin{equation*}
  |\bar\partial\mathcal X_n(\zeta)|\leq \frac{1}{\rho_{\phi^*}(0)(R_n^2-R_n)}.
\end{equation*}
By Lemma \ref{christ} $\rho_{\phi^*}(\zeta)/\rho_{\phi^*}(0) 
\leq R_n^{2(1-\delta)}$ for some $\delta \in(0,1)$. Thus, if $R_n$ is big 
enough
\begin{equation*}
  |\bar\partial\mathcal X_n(\zeta)|\leq 
  \frac{ R_n^{2(1-\delta)}}{\rho_{\phi^*}(\zeta)(R_n^2-R_n)} \leq  
  \frac{\varepsilon_n}{\rho_{\phi^*}(\zeta)}.
\end{equation*}
Take now $j_n$ big enough so that $\rho_{\phi_{j_n}}/\rho_{\phi^*} \leq 2$ on 
$\widetilde{D}_n$ and
\begin{eqnarray*}
\bigl| \|f_n | \widetilde{D}_n\|_{\calf_{\phi_{j_n},\omega_{j_n}}^{p}} -  
\|f_n | \widetilde{D}_n \|_{\calf_{\phi^*,\omega^*}^{p}} 
 \bigr| &\leq &\varepsilon_n,\\
\bigl| \| f_n | \Lambda_{j_n}\cap\widetilde{D}_n
\|_{\ell^{p}_{\phi_{j_n},\omega_{j_n}}} -  
\| f_n| \Lambda^* \cap \widetilde{D}_n \|_{\ell^{p}_{\phi^*,\omega^*}} 
\bigr| &\leq &\varepsilon_n.
\end{eqnarray*}
Define $g_n=f_n\mathcal X_n$. Then $\bar\partial g_n$ is supported on 
$C_n:=\{R_n\leq|\zeta|\leq R_n^2\}$ and 
$  |\bar\partial g_n(\zeta)| \leq 
\varepsilon_n|f_n(\zeta)| / \rho_{\phi^*}(\zeta)$, so by Theorem 
\ref{estimacions-dbar} there exists $u_n$ solution to 
$\bar\partial u_n=\bar\partial g_n $ with
\begin{equation*}
\|u_n\|_{\calf_{\phi_{j_n},\omega_{j_n}}^{p}} \lesssim \|\bar\partial g_n 
\rho_{\phi_{j_n}} 
\|_{\calf_{\phi_{j_n},\omega_{j_n}}^{p}} \lesssim \varepsilon_n  \|f_n | 
\widetilde{D}_n\|_{\calf_{\phi_{j_n},\omega_{j_n}}^{p}} 
\lesssim \varepsilon_n.
\end{equation*}
The function $G_n=g_n-u_n$ is holomorphic and satisfies
\begin{equation*}\label{G_n}
\|G_n\|_{\calf_{\phi_{j_n},\omega_{j_n}}^{p}}\geq \|f_n 
|D_n\|_{\calf_{\phi_{j_n},\omega_{j_n}}^{p}} - 
\|u_n\|_{\calf_{\phi_{j_n},\omega_{j_n}}^{p}} \geq 1-C\varepsilon_n \simeq 1 .
\end{equation*}
We will check now that $G_n|\Lambda_{j_n}$ is small. Split
$\Lambda_{j_n}$ into $\widetilde{\Lambda}_{j_n}=\Lambda_{j_n} \cap \{D_n \cup 
(\C\setminus \widetilde{D}_n))\}$ and 
$\widehat{\Lambda}_{j_n}=\Lambda_{j_n} \setminus \widetilde{\Lambda}_{j_n}$. On
the one hand 
\begin{equation*} \label{G_n_lambda}
\|G_n |\widetilde{\Lambda}_{j_n}
\|_{\ell_{\phi_{j_n},\omega_{j_n}}^{p}(\widetilde{\Lambda}_{j_n})} \leq 
\| f_n |\widetilde{D}_n\cap \widetilde{\Lambda}_{j_n} \|_{\ell_{\phi_{j_n},
\omega_{j_n}}^{p}(\widetilde{\Lambda}_{j_n})} + 
\|u_n |\widetilde{\Lambda}_{j_n} \|_{\ell_{\phi_{j_n},\omega_{j_n}}^{p}
(\widetilde{\Lambda}_{j_n})} .
\end{equation*}
From $\|u_n |\widetilde{\Lambda}_{j_n} \|_{\ell_{\phi_{j_n},\omega_{j_n}}^{p}
(\widetilde{\Lambda}_{j_n})} 
\leq 
\|u_n\|_{\calf_{\phi_{j_n},\omega_{j_n}}^{p}} \leq \varepsilon_n$  (by Lemma 
\ref{estimacio-puntual} 
for the case $p<\infty$, since $u$ is holomorphic in $D_n \cup (\C\setminus
\widetilde{D}_n$ )) we deduce that  
$\|G_n | \widetilde{\Lambda}_{j_n}\|_{\ell_{\phi_{j_n},\omega_{j_n}}^{p}
(\widetilde{\Lambda}_{j_n})} 
\lesssim 
\varepsilon_n$. 
On the other hand 
\[
\|G_n |\widehat{\Lambda}_{j_n}
\|_{\ell_{\phi_{j_n},\omega_{j_n}}^{p}(\widehat{\Lambda}_{j_n})} \lesssim
\|G_n|(\widetilde{D}_n\setminus D_n)\|_{\calf_{\phi_{j_n},\omega_{j_n}}^{p}}
\lesssim
\||f_n|+|u_n||(\widetilde{D}_n\setminus D_n)
\|_{\calf_{\phi_{j_n},\omega_{j_n}}^{p}}
\lesssim \varepsilon_n.
\]
This together with the above and the fact that the sampling 
constants of $\Lambda$ and  $\Lambda_{j_n}$ coincide (Lemma 
\ref{properties}(b)) leads to 
contradiction.

(b) Assume that $\Lambda^*=\{\lambda_k^*\}_k$, and let $v \in 
\ellpa(\Lambda^*)$ with
$\|v\|_{\ellpa(\Lambda^*)} \leq 1$. Let also 
$\Lambda_j=\{\lambda_{k}^{j}\}_k $ be such that 
$\Lambda_j \rightarrow \Lambda^*$ uniformly on compact sets. 
For $\varepsilon_n$ decreasing to zero and $R_n$ big enough 
(to be chosen later) there exists $j_n$ such that
$\|v\|_{\ell_{\phi_{j_n}}^{p}(\Lambda_{j_n}\cap D^{R_n}_{\phi^*}(0))} 
\leq 2$ 
and 
\begin{equation}\label{choix2}
\frac{e^{-\phi^*}\omega^*\rho_{\phi^*}^{-2/p}}{e^{-\phi_{j_n}}
\omega_{j_n}^*\rho_{\phi_{j_n}}^{-2/p}} \leq 2 \qquad \text{on} \qquad 
D_{\phi^*}^{R_n^2}(0).
\end{equation}

Since the interpolation constant $M(\Lambda_j)$ does not depend on $j$
 there exist 
$f_n \in \calf^{p}_{\phi_{j_n},\omega_{j_n}}$ with 
$\|f_n\|_{\calf^{p}_{\phi_{j_n},\omega_{j_n}}} \leq 2M(\Lambda)$ and 
\begin{equation*}
f_n(\lambda_{k}^{j_n})=
\begin{cases} 
v_k &\text{ if $ \lambda_{k}^{j_n} \in D_{\phi^*}^{R_n} (0)$} \\
0 &\text{otherwise.} 
\end{cases}
\end{equation*}
We will now use the same technique as in (a) to modify $f_n$ so that
it falls in $\calf_{\phi^*,\omega^*}^{p}$. Take the cut-off function 
$\mathcal X_n$
constructed above, define $g_n=f_n\mathcal X_n$ and consider 
a solution $u_n$ to $\bar\partial u_n=f_n\bar\partial(\mathcal X_n) $
such that:
\begin{eqnarray*}
\|u_n\|_{\calf_{\phi^*,\omega^*}^{p}} &\lesssim&
\|f_n\bar\partial(\mathcal X_n)\rho_{\phi^*}\|_{\calf_{\phi^*,\omega^*}^{p}} 
\lesssim \varepsilon_n\|f_n|D^{R_n^2}_{\phi^*}(0)
\|_{\calf_{\phi^*,\omega^*}^{p}} 
\lesssim \varepsilon_n\|f_n\|_{\calf_{\phi_{j_n},\omega_{j_n}}^{p}} 
\lesssim \varepsilon_n, \\
\|u_n\|_{\calf_{\phi^*,\omega^*}^{\infty}} &\lesssim& 
\|f_n\bar\partial(\mathcal X_n)\rho_{\phi^*}\|_{\calf_{\phi^*,
\omega^*}^{\infty}} 
\lesssim \varepsilon_n
\|f_n\|_{\calf_{\phi^*,
\omega^*}^{p}} \lesssim \varepsilon_n.
\end{eqnarray*} 
According to Theorem C and \eqref{choix2} such a solution always exists.

The entire function $G_n=f_n\bar\partial(\mathcal X_n)-u_n$
is $\calf_{\phi^*,\omega^*}^{p}$  and $\|G_n\|_{\calf_{\phi^*,\omega^*}^{p}} 
\leq CM$. By Montel's theorem we may assume that $G_n$ converges to a 
function $G \in \calf_{\phi^*,\omega^*}^{p}$. Notice that
$G_n(\lambda_{k}^{j_n})=v_k-u_n(\lambda_{k}^{j_n})$ for  
$\lambda_{k}^{j_n} \in D^{R_n}_{\phi^*}(0)$,
and by the $L^\infty$ estimates,
$|u_n(\lambda_{k}^{j_n}) | $
 tends to zero as $n$ goes to infinity.
Therefore $G$ interpolates $v$.
\end{proof}

\begin{lemma}
Suppose that for every weak limit $(\Lambda^*,\phi^*,\omega^*)\in
W(\Lambda,\phi,\omega)$ the sequence $\Lambda^*$
is a uniqueness set for
$\calf^{\infty}_{\phi^*,\omega^*}$. 
Then there exists $\varepsilon>0$ such that 
$\Lambda$ is sampling for 
$\calf^{\infty}_{(1+\varepsilon)\phi,\omega}$.
\end{lemma}

\begin{proof}
If this is not the case there exist $\varepsilon_n > 0$
decreasing to 0,  
$f_n \in\calf_{(1+\varepsilon_n)\phi,\omega}^{\infty}$ and $z_n\in\C$ such that
$|f_n(z_n)|e^{-(1+\varepsilon_n)\phi(z_n)}\omega(z_n)=1$, 
$\|f_n\|_{\calf_{(1+\varepsilon_n)\phi,\omega}^{\infty}}\leq 2$ and 
$\|f_n |\Lambda \|_{\ell^{\infty}_{(1+\varepsilon_n)\phi,\omega}(\Lambda)} \leq 
\varepsilon_n$.

Denote $\psi_n=(1+\varepsilon_n)\phi$. Let 
$\Lambda_n=(\tau_{z_n}^{\psi_n})^{-1}(\Lambda)$, $\omega_n=\omega_{z_n}$ and 
$g_n=T_{z_n}^{\psi_n,\omega_n} f_n$.
Then, denoting $\psi_{n,z_n} =(1+\varepsilon_n)\phi_{z_n}$, we have
$|g_n(0)|=1$ and
$\|g_n|\Lambda_n\|_{\ell^{\infty}_{\psi_{n,z_n},\omega_n}(\Lambda_n)  } = 
\|{f_n |\Lambda}\|_{\ell^{\infty}_{\psi_n,\omega_n}
(\Lambda)} \leq \varepsilon_n$.
Taking a subsequence if necessary, we can assume that $\Lambda_n$ converges 
weakly to $\Lambda^*$, $\psi_{n,z_n} \rightarrow \phi^*$, $\omega_n\to\omega^*$
 uniformly on compact
sets and $g_n \rightarrow g^* \in \calf^{\infty}_{\phi^*,\omega^*}$
(by Montel's Theorem). So $g^*$ vanishes on $\Lambda^*$ and $|g^*(0)|=1$, 
contradicting the fact that $\Lambda^*$ is a uniqueness sequence.
\end{proof}

\begin{corollary}\label{weak-limit}
Let $\phi$ a subharmonic function with doubling Laplacian, let
$\omega\in \Omega_{\phi}$ and let 
$\Lambda$ be a $\rho$-separated sequence. The sequence
$\Lambda$ is in $\Samp  \Finfa$ if and only if for all weak limit 
$(\Lambda^*,\phi^*,\omega^*)\in W(\Lambda,\phi,\omega)$, the sequence
$\Lambda^*$ is a  
uniqueness set for $\mathcal F_{\phi^*,\omega^*}^{\infty}$. 
\end{corollary}

\subsection{Non-existence of simultaneously sampling and interpolating
sequences} An important result in the proof of Theorems A and B is the
following.

\begin{theorem}\label{no-sampling-interpolation}
There is no sequence $\Lambda$ both sampling and 
interpolating  for $\Fpa$, 
$p\in [1,\infty]$.
\end{theorem}
\begin{proof}
Assume that such sequence $\Lambda$ exists. We claim that
\begin{equation}\label{key}
\sup_{\lambda^*\in\Lambda}\sum_{\lambda\in\Lambda\setminus
\lambda^*}\frac{\rho(\lambda)\rho(\lambda^*)}{|\lambda-\lambda^*|^2} < \infty.
\end{equation}
Let $p\in[1,\infty)$. Given any $\lambda^*\in \Lambda$ take a
function $g$ such that $g(\lambda^*)=1$, $g(\lambda)=0$ for
$\lambda\neq\lambda^*$ and $\|g\|_{\Fpa}\lesssim
e^{-p\phi(\lambda^*)}\omega^{p}(\lambda^*)$. Such $g$ exists 
because $\Lambda$ is interpolating. Consider the function
\[
F(z)=\sum_{\lambda\in\Lambda\setminus\lambda^*}\rho(\lambda)
\frac{g(z)(z-\lambda^*)}{(z-\lambda)\overline{(\lambda^*-\lambda)}}.
\]
The sampling inequality shows that $F\in \Fpa$. 
Moreover, since $|F(\lambda)| = |g'(\lambda)|\rho(\lambda)$ for 
all $\lambda\in \Lambda\setminus \lambda^*$ and $F(\lambda^*)=0$, we have
\[
\|F\|^p_{\Fpa}\lesssim \sum_{\lambda\in\Lambda\setminus \lambda^*}
 |g'(\lambda)|^p
\rho^p(\lambda)e^{-p\phi(\lambda)}\omega^{p}(\lambda).
\]
We use now Lemma~\ref{derivative} and the fact that $\Lambda$ is
$\rho$-separated (since it is interpolating):
\[
\|F\|^p_{\Fpa} \lesssim \sum_{\Lambda\setminus\lambda^*}
\int_{D(\lambda)} |g|^p e^{-p\phi}\omega^p d\sigma 
\lesssim \|g\|^p_{\Fpa}\lesssim e^{-p\phi(\lambda^*)}
\omega^p(\lambda^*).
\]
We want to estimate $|F'(\lambda^*)|$. Using again Lemma~\ref{derivative}
\[
|F'(\lambda^*)|^p e^{-p\phi(\lambda^*)}\omega^p(\lambda^*)
\rho^{p}(\lambda^*)\lesssim
\int_{D(\lambda)} |F|^p e^{-p\phi}\omega^p d\sigma \lesssim
e^{-p\phi(\lambda^*)}\omega^{p}(\lambda^*).
\]
Therefore $|F'(\lambda^*)|\rho(\lambda^*)\lesssim 1$. On the other hand
\[
F'(\lambda^*)=\sum_{\lambda\in\Lambda\setminus \lambda^*}\frac{\rho(\lambda)}
{|\lambda-\lambda^*|^2}.
\]
This yields \eqref{key}. The obvious modifications give \eqref{key} in 
the case $p=\infty$.

According to Lemma~\ref{separacio}(d) there exists $r>0$ with 
$\C=\cup_{\lambda\in\Lambda} 
D^r(\lambda)$.
Also, there exists $r_0>0$ depending on $r$ such that,
\[ 
\int_{ D^r(\lambda) \setminus D^{r_0}(\lambda^*)}
\frac{dm(z)}{1+|z-\lambda^*|^2}
\le C(r) \frac{\rho^2(\lambda)}{|\lambda-\lambda^*|^2}
\qquad \forall \lambda \notin D^{r_0}(\lambda^*) .
\]
We may now finish by taking a big disk $D(0,M)$ and $\lambda_M^*\in
D(0,M)$ in
such a way that $\rho(\lambda_M^*)\ge \rho(\lambda)$ for all $\lambda\in
\Lambda\cap D(0,M)$. In this case
\[
\int_{D(0,M)\setminus D^{r_0}(\lambda_M^*)}
\frac{dm(z)}{1+|z-\lambda_M^*|^2}\lesssim 
\sum_{
\substack
{\lambda\in \Lambda \\ 
\lambda\notin D^{r_0}(\lambda^*) } } 
\int_{D^r(\lambda)} \frac{dm(z)}{1+|z-\lambda_M^*|^2}
\lesssim
\sum_{\lambda\in\Lambda\setminus\lambda_M^*}
\frac{\rho(\lambda)\rho(\lambda_M^*)}{|\lambda-\lambda_M^*|^2} < C .
\]
This is a contradiction, since $\lim_{M\to\infty}\rho(\lambda_M^*)/M =0$ 
and the left hand side of the previous inequality tends to
$\infty$ as $M$ goes to $\infty$.
\end{proof}

\begin{corollary}\label{add-delete-point}
Any sequence obtained by deleting a finite number of points of  
$\Lambda\in\Samp\Fpa$ is still in $\Samp\Fpa$.
\end{corollary}

We want to prove next an analogue  for interpolating sequences: adding
a finite number of points to an interpolating sequence gives again an
interpolating sequence.
 
Given $\Lambda$ and a point $z$ define, following 
\cite[p.352--354]{Br}
\begin{equation*}
\sigma_{\phi,\omega}^{p}(z,\Lambda) := 
\sup\bigl\{|f(z)|e^{-\phi(z)}\omega(z), \ \|f\|_{\Fpa} 
\leq 1,\ f|\Lambda\equiv 0 \bigr\}.
\end{equation*}
Notice first that if $\Lambda$ is interpolating and $z \notin \Lambda$ 
this is strictly positive. Indeed, $\Lambda$ is not a uniqueness
sequence, otherwise $\Lambda$ would be also sampling, contradicting 
Theorem \ref{no-sampling-interpolation}. Thus there exists $f \in \Fpa$,
$f\neq 0$ with 
$f|\Lambda\equiv 0$ and, eventually dividing $f$ by a power of $(\zeta-z)$,
$f(z) \neq 0$. Hence $\sigma_{\phi,\omega}^{p}(z,\Lambda) >0$.

\begin{lemma}\label{add-point}
Let $\Lambda\in\Int \Fpa$. Then 
$\Lambda \cup\{z\}\in\Int\Fpa$ for all $z \notin \Lambda$. Furthermore, 
for all $\varepsilon >0$ there exists $C>0$ such that  
$d_\phi(\Lambda,z)\geq \varepsilon$ implies 
$M^{p}_{\phi,\omega}(\Lambda \cup\{z\}) \leq C M^{p}_{\phi,\omega}(\Lambda)$.
\end{lemma}
\begin{proof}
As in the proof of \cite[Lemma 4, p.233]{Br}, we have
\begin{equation*}
M^{p}_{\phi,\omega}(\Lambda \cup\{z\}) 
\leq \frac{1+2 M^{p}_{\phi,\omega}(\Lambda)}
{\sigma_{\phi,\omega}^{p}(z,\Lambda)}.
\end{equation*}
Thus we will be done if we prove that
there exists $A>0$ such that $d_\phi(z,\Lambda) \geq \varepsilon$ implies 
$\sigma_{\phi,\omega}^{p}(z,\Lambda) \geq A$.

If this is not true, there exists a sequence $\{z_n\} \in \C$ with 
$d_\phi(z_n,\Lambda) \geq \varepsilon $ and  
$\sigma_{\phi,\omega}^{p}(z_n,\Lambda) \leq 1/n$.  Transferring $z_n$ to the
origin by $\tau_{z_n}^{-1}$ (see Section~\ref{translations}), we get a
sequence $\Lambda_n:=\Lambda_{z_n}$ such that  $|\lambda|\geq \varepsilon$
for all $\lambda\in \Lambda_n$
and  $\sigma_{\phi_n,\omega_n}^{p}(0,\Lambda_n) \leq 1/n$, where
$\phi_n=\phi_{z_n}$ and $\omega_n=\omega_{z_n}$. 

Taking a subsequence if necessary, assume that $(\Lambda_n,\phi_n,\omega_n)$ 
converges to $(\Lambda^*,\phi^*,\omega^*)$. 
By Proposition \ref{weak_stability},  $\Lambda^*\cup
\{0\}\in\Int\mathcal{F}_{\phi^*,\omega^*}^{p}$, so there exists $f \in
\mathcal{F}_{\phi^*,\omega^*}^{p}$ with $f|\Lambda^*=0$ and $|f(0)|=1$. Arguing 
as in
the proof of Proposition \ref{weak_stability} we see that there exist $f_n \in
\mathcal{F}_{\phi_n,\omega}^{p}$ 
and $\varepsilon_n$ decreasing to zero such that 
\begin{equation*}
\|f_n|\Lambda_n\|_{\ell^{p}_{\phi_n,\omega_{n}}(\Lambda_n)} 
\leq \varepsilon_n,\qquad |f_n(0)|\geq c \qquad \text{and}\qquad 
\|f_n\|_{ \mathcal{F}_{\phi_n,\omega_n}^{p}} \leq C. 
\end{equation*}
Since $\Lambda_n$ is interpolating, there exist also 
$g_n \in \mathcal{F}_{\phi_n,\omega_n}^{p}$ with 
\begin{equation*}
g_n|\Lambda_n=f_n|\Lambda_n \qquad \text{and} 
\qquad \|g_n\|_{\calf_{\phi_n,\omega_n}^{p}}\leq 
M_{\phi_n,\omega_n}^{p}(\Lambda_n)
\|f_n|\Lambda_n\|_{\ell^{p}_{\phi_n,\omega_n}(\Lambda_n)} \leq 
\varepsilon_n M(\Lambda) .
\end{equation*}
Then $h_n:=f_n-g_n \in \calf_{\phi_n,\omega_n}^{p}$ 
vanishes on $\Lambda_n$ and
$\| h_n \|_{\calf_{\phi_n,\omega_n}^p}
\leq 2C$, therefore $|h_n(0)|\lesssim 1/n$.  On the other hand
$|g_n(0)|\lesssim \varepsilon_n$ and therefore $|h_n(0)|\geq c/2$, thus 
contradicting the previous estimate.
\end{proof}

\subsection{Inclusions between various spaces of interpolating sequences}
We want to study next the relationship between
the spaces of interpolating sequences for various weights. We will use the
techniques already exploited in \cite{MaPa00}.

We start with the construction of a sort of peak-functions
associated to an interpolating sequence.
Let $\delta_\lambda^{\lambda'}$ denote the Kroenecker indicator, i.e.
$\delta_\lambda^{\lambda'}=1 $ if $\lambda=\lambda'$ and 
$\delta_\lambda^{\lambda'}=0$ otherwise.
 
\begin{lemma}\label{peak-functions} Let $\Lambda\in\Int \Fpa$, 
$1\leq p\leq \infty$.
Given $\varepsilon>0$ and $\widetilde{\omega}\in \Omega_{\phi}$
there exist $m\in\N$, $C>0$ and functions 
$g_\lambda\in \mathcal F_{(1+\varepsilon)\phi,\tilde{\omega}}^{p}$  such that
\begin{itemize}
\item[(a)] $g_\lambda(\lambda')=\delta_\lambda^{\lambda'}\; $ for all 
$\lambda,\lambda'\in \Lambda$.
\item[(b)] $\|g_\lambda\|_{\calf_{(1+\varepsilon)\phi,\tilde{\omega}}^{p}}
\simeq
\widetilde{\omega}(\lambda) e^{-(1+\varepsilon)\phi(\lambda)}$.
\item[(c)] $\displaystyle |g_\lambda(z)|\lesssim
\frac{\widetilde{\omega}(\lambda)}{\widetilde{\omega}(z)}
e^{(1+\varepsilon)(\phi(z)-\phi(\lambda))}\frac{1}{1+d_\phi^m(z,\lambda)}$.
\item[(d)] For all $v\in
\ell_{(1+\varepsilon)\phi,\tilde{\omega}}^{p}(\Lambda)$, 
$\|v\|_{\ell_{(1+\varepsilon)\phi,\tilde{\omega}}^{p}(\Lambda)}
\lesssim \bigl\|\sum_{\lambda\in\Lambda} 
v_\lambda g_\lambda\bigr\|_{\calf_{(1+\varepsilon)\phi,\tilde{\omega}}^{p}}
\lesssim 
\|v\|_{\ell_{(1+\varepsilon)\phi,\tilde{\omega}}^{p}(\Lambda)}\ $.
\item[(e)] $\displaystyle \lim\limits_{r\to\infty}
\sup_{\lambda\in\Lambda}
\frac{e^{p(1+\varepsilon)\phi(\lambda)}}{\widetilde{\omega}^{p}(\lambda)}
\int\limits_{\C\setminus D^r(\lambda)}
|g_\lambda(z)|^p e^{-p(1+\varepsilon)\phi(z)} 
\widetilde{\omega}^{p}(z) d\sigma(z)=0$.
\end{itemize}
\end{lemma}

\begin{proof}
By hypothesis, there exist functions $f_\lambda\in\Fpa$ such that
$f_\lambda(\mu)=\delta_\lambda^\mu$ for all $\lambda,\mu\in\Lambda$ and
$\|f_\lambda\|_{\Fpa}\leq M(\Lambda) \omega(\lambda) e^{-\phi(\lambda)}$.
Consider the weights $P_\lambda$ given by Theorem~\ref{weight} for 
the weight $\widetilde{\omega}/\omega$, and define
$g_\lambda=f_\lambda P_\lambda$.
By construction we have (a) and (c).

(b) When $p=\infty$ we have
$\widetilde{\omega}(\lambda) e^{-(1+\varepsilon)\phi(\lambda)}=
\widetilde{\omega}(\lambda) e^{-(1+\varepsilon)\phi(\lambda)} 
|g_\lambda(\lambda)|\leq
\|g_\lambda\|_{\calf_{(1+\varepsilon)\phi,\tilde{\omega}}^{\infty}}$. 
The remaining inequality is immediate from (c).

Let $p<\infty$. On the one hand, Lemma \ref{estimacio-puntual}(a) gives
\[
\widetilde{\omega}(\lambda) e^{-(1+\varepsilon)\phi(\lambda)}=
\widetilde{\omega}(\lambda) e^{-(1+\varepsilon)\phi(\lambda)} 
|g_\lambda(\lambda)
|\lesssim
\bigl(\int_{D(\lambda)}|g_\lambda|^p
e^{-p(1+\varepsilon)\phi} \widetilde{\omega}^p d\sigma\bigr)^{1/p}\lesssim 
\|g_\lambda\|_{\calf_{(1+\varepsilon)\phi,\tilde{\omega}}^{p}} .
\]
On the other hand,  (c) and Lemma \ref{integrals}(b) show that for $m$ big 
enough
\begin{eqnarray*}
\int_{\C}|g_\lambda|^p
e^{-p(1+\varepsilon)\phi} \widetilde{\omega}^pd\sigma 
&\lesssim& \widetilde{\omega}^p(\lambda)
e^{-p(1+\varepsilon)\phi(\lambda)} 
\left[\int_{D(\lambda)} d\sigma(z) 
+ \int_{\C\setminus D(\lambda)}
\frac{\Delta\phi(z)}{d_\phi^{pm}(z,\lambda)}\right]\\
&\lesssim& \widetilde{\omega}^p(\lambda)
e^{-p(1+\varepsilon)\phi(\lambda)} 
\end{eqnarray*}

(d) Denote $f=\sum_{\lambda} v_\lambda g_\lambda$. The left inequalities are
proved similarly to (b), for 
\[
\widetilde{\omega}(\lambda) e^{-(1+\varepsilon)\phi(\lambda)}|v_\lambda|=
\widetilde{\omega}(\lambda) e^{-(1+\varepsilon)\phi(\lambda)}|f(\lambda)| .
\]
For $p=\infty$ and $v\in\ell_{(1+\varepsilon)\phi,\tilde{\omega}}^{\infty,}
(\Lambda)$ 
Lemma \ref{suma} and (c) yield
\[
\widetilde{\omega}(z) e^{-(1+\varepsilon)\phi(z)}\bigl(\sum_{\lambda\in\Lambda}
|v_\lambda|
|g_\lambda(z)|\bigr)\lesssim 
\|v\|_{\ell_{(1+\varepsilon)\phi,\tilde{\omega}}^{\infty}(\Lambda)}
\sum_{\lambda\in\Lambda} 
\frac 1{1+d^m_\phi(\lambda,z)}
\lesssim \|v\|_{\ell_{(1+\varepsilon)\phi,\tilde{\omega}}^{\infty}(\Lambda)}.
\]

Let now $p<\infty$. Using the estimate (c) and Jensen's inequality for
convex functions (which is legitimate thanks to Lemma~\ref{suma}) we have 
\begin{eqnarray*}
|f(z)|^p e^{-p(1+\varepsilon)\phi(z)} \widetilde{\omega}^p(z) 
\rho^{-2}(z)&\lesssim&  \frac 1{\rho^2(z)}
\Bigl[
\sum_{\lambda\in\Lambda} \widetilde{\omega}(\lambda) |v_\lambda|
e^{-(1+\varepsilon)\phi(\lambda)}\frac 1{1+d_\phi^m(z,\lambda)}
\Bigr ]^p \\
&\lesssim& 
\frac 1{\rho^2(z)}\sum_{\lambda\in\Lambda} 
\widetilde{\omega}^p(\lambda) |v_\lambda|^p
e^{-p(1+\varepsilon)\phi(\lambda)}
\frac 1{1+d_\phi^{m}(z,\lambda)}.
\end{eqnarray*}
Now we apply Lemma \ref{integrals}(b) and obtain
\[
\int_{\C}|f|^p e^{-p(1+\varepsilon)\phi } 
\widetilde{\omega}^p\ d\sigma
\lesssim
\sum_{\lambda\in\Lambda} \widetilde{\omega}^p(\lambda) |v_\lambda|^p
e^{-p(1+\varepsilon)\phi(\lambda)}
\int_{\C} 
\frac {\Delta\phi(z)}{1+d_\phi^{m}(z,\lambda)}
\lesssim \|v\|^p_{\ell_{(1+\varepsilon)\phi,\tilde{\omega}}^{p}(\Lambda)}.
\]

(e) This follows from (c) and Remark \ref{int-to-zero}, since
\[
\frac{e^{p(1+\varepsilon)\phi(\lambda)}}{\widetilde{\omega}^p(\lambda)}
\int_{\C\setminus D^r(\lambda)}
|g_\lambda(z)|^p e^{-p(1+\varepsilon)\phi(z)} \widetilde{\omega}^p(z)d\sigma(z) 
\lesssim \int_{\C\setminus D^r(\lambda)}
\frac {\Delta\phi(z)}{d_\phi^m(z,\lambda)}.
\]
\end{proof}

\begin{theorem}\label{inclusions} For all $\varepsilon>0$, $1\leq p,p'\leq
\infty$ and $\omega,\widetilde{\omega}\in \Omega_{\phi}$, 
the following inclusions hold
\[
\Int \Fpa \subset  
\Int \mathcal F_{(1+\varepsilon)\phi,\tilde{\omega}}^{p'} .
\]
\end{theorem}

\begin{proof}
It will be enough to prove that for all $\varepsilon>0$, $1\leq p\leq\infty$ 
and  $\omega,\widetilde{\omega}\in \Omega_{\phi}$,
\begin{eqnarray*}
\textrm{(a)}\quad \Int \Fpa\subset\Int 
\mathcal F_{(1+\varepsilon)\phi,\tilde{\omega}}^{\infty}\qquad\qquad
\textrm{(b)}\quad \Int \Finfa\subset
\mathcal F_{(1+\varepsilon)\phi,\tilde{\omega}}^{p} .
\end{eqnarray*}

(a) Take the functions $g_\lambda$ given by Lemma \ref{peak-functions}. For
$v\in\ell_{(1+\varepsilon)\phi,\tilde{\omega}}^{\infty}(\Lambda)$ we consider
the interpolating function 
\[
f(z)=\sum_{\lambda\in\Lambda} v_\lambda g_\lambda(z) \ 
\]
A direct estimate using Lemma~\ref{peak-functions}(c) yields
\[
\widetilde{\omega}(z)|f(z)| e^{-(1+\varepsilon)\phi(z)}\lesssim
\sum_{\lambda\in\Lambda} 
\frac 1{1+d_\phi^m(z,\lambda)},
\]
which is bounded, by Lemma \ref{suma}.

(b) Given
$v\in\ell_{(1+\varepsilon)\phi,\tilde{\omega}}^{p}(\Lambda)$, take
$f=\sum_{\lambda} v_\lambda g_\lambda$ as before and estimate as in the proof of
Lemma~\ref{peak-functions}(d).
\end{proof}

\subsection{Inclusions between various spaces of sampling sequences} 
In this section we want 
to prove some inclusions between various spaces of
sampling sequences. Unlike in the corresponding result 
for interpolating sequences, for the spaces of sampling sequences 
there is a gain, in the sense that
any sampling sequence 
is actually sampling for a slightly bigger space. This will allow us to pass
from the non-strict to the strict inequality of Theorem A.

\begin{theorem}\label{inclusions-sampling}
Let $\Lambda\in\Samp\Fpa$ be  $\rho$-separated. There exists
$\varepsilon>0$ such that for all $p'\in[1,\infty]$ and $\widetilde{\omega}\in
\Omega_\phi$, the sequence
$\Lambda\in\Samp \mathcal F^{p'}_{(1+\varepsilon)\phi,\tilde{\omega}}$.
\end{theorem}

\begin{proof} The proof is divided in three steps.

(a) \emph{If $\Lambda\in\Samp\Fpa$, then $\Lambda\in\Samp \Finfa$}.
We know from Proposition~\ref{weak_stability} that for all weak limit
$(\Lambda^*,\phi^*,\omega^*)$ the sequence  
$\Lambda^* $ is in $\Samp \mathcal F^{p}_{\phi^*,\omega^*} $, and
by Lemma \ref{weak-limit} it will be enough to see that all
weak limit $\Lambda^*$
is a  uniqueness set for $\mathcal F^{\infty}_{\phi^*,\omega^*}$.
 
If this is not the case, there exists
$f \in \mathcal F^{\infty}_{\phi^*,\omega^*}$ with 
$f|\Lambda^* \equiv 0$, $f\neq 0$. 

We claim that for $m$ large enough 
\begin{equation*}
g(z):=\frac{f(z)}{(z-\lambda_1^* )\ldots(z-\lambda_m^*) } \in 
\mathcal F^{p}_{\phi^*,\omega^*} .
\end{equation*}

It is clear that Lemma~\ref{derivative} gives the $p$-integrability on 
$\cup_{j=1}^m D(\lambda_j^*)$. On the other hand, by 
Lemma~\ref{estimacio-puntual}
\begin{equation*}
\int_{z\notin\cup_j D(\lambda_j^*)} 
\frac{|f|^p e^{-p\phi^*} {\omega^*}^{p}
 \rho_{\phi^*}^{-2}}{|z-\lambda_1^* |^p\ldots|z-
\lambda_m^*|^p } \leq C\int_{z\notin\cup_j D(\lambda_j^*)} 
\frac{\|f\|_{\calf_{\phi^*,\omega^*}^{\infty}}^p \Delta\phi^*}{|z-\lambda_1^* 
|^p\ldots|z-\lambda_m^*|^p } .
\end{equation*}
Since $\Delta \phi^*$ is doubling there exists 
$m$ such that this integral converges (Lemma \ref{integrals}(b)).

By Corollary \ref{add-delete-point},  
$\Lambda^*\setminus \{\lambda_1^*\ldots\lambda_m^*\}\in\Samp 
\mathcal F^{p}_{\phi^*,\omega^*} $. As $f$ vanishes on this sequence 
we deduce that $f \equiv 0$, which is a contradiction. 

(b) \emph{If $\Lambda\in\Samp \Finfa$ there exists $\varepsilon>0$ such that 
$\Lambda\in\Samp \mathcal F^{\infty}_{(1+\varepsilon)\phi,\omega}$}. If this
is not the case for any sequence $\{\varepsilon_n\}\searrow 0$ there exist 
functions $f_n \in \mathcal  F^{\infty}_{(1+\varepsilon_n)\phi,\omega}$ and
$\delta_n>0$ decreasing to  $0$ with $\|
{f_n|\Lambda}\|_{\ell^{\infty}_{(1+\varepsilon_n)\phi,\omega}(\Lambda)} \leq 
\delta_n$  and $ |f_n(z_n)|=1$.

Let $\Lambda_n=\tau_{z_n}^{-1}(\Lambda)$, $\phi_n=
(1+\varepsilon_n)\phi_{z_n}$,
$\omega_n=\omega_{z_n}$ 
and  $\tilde f_n=T_{z_n}^{\phi,\omega} f_n$. Then $|\tilde f_n(0)|=1$, $\|
{\tilde f_n}{|\Lambda_n}\|_{\ell^{\infty}_{\phi_n,\omega_n}}\leq 
\delta_n$, and there exist a sequence $\Lambda^*$ and functions
$\phi^*, f^*,\omega^*$ such that 
\[
(\Lambda_n,\phi_n,\omega_n)\rightarrow 
(\Lambda^*,\phi^*,\omega^*)\in W(\Lambda,\phi,\omega)
\]
and $\{f_n\}_n\rightarrow f^*\in \mathcal{F}^\infty_{\phi^*,\omega^*}$ 
uniformly on compact sets.  So we have
$|f^*(0)|=1$ and $f^* |\Lambda^*=0$, i.e.  $\Lambda^*$ is not a uniqueness
sequence for $\mathcal{F}^\infty_{\phi^*,\omega^*}$, a contradiction
with Lemma \ref{weak-limit}.

(c) \emph{If $\Lambda\in\Samp \mathcal 
F^\infty_{(1+\varepsilon)\phi,\omega}$ for some $\varepsilon>0$, then 
$\Lambda\in\Samp \mathcal F^{p'}_{\phi,\tilde{\omega}}$, 
for all $\widetilde{\omega}\in\Omega_\phi$, $1\le p'\le\infty$}.
Consider the spaces
\begin{eqnarray*}
\mathcal F_{(1+\varepsilon)\phi,\omega}^{\infty,0} &= &
\{ f\in \mathcal 
F^\infty_{(1+\varepsilon)\phi,\omega} :
\lim_{|z|\to\infty} \omega(z) |f(z)| e^{-(1+\varepsilon)\phi(z)}=0\}, \\
\ell_{(1+\varepsilon)\phi,\omega}^{\infty,0}(\Lambda) &= &
\{ v\in \ell^{\infty}_{(1+\varepsilon)\phi,\omega} :
\lim_{|\lambda|\to\infty} \omega(z) 
|v_\lambda| e^{-(1+\varepsilon)\phi(\lambda)}=0\} .
\end{eqnarray*}

There is a
sequence of functions $\{g(z,\lambda)\}_{\lambda\in\Lambda}$ such that for all
$f\in\mathcal F^{\infty,0}_{(1+\varepsilon)\phi,\omega}$
\[
\omega(z) e^{-(1+\varepsilon)\phi(z)} f(z)=\sum_{\lambda\in\Lambda} 
\omega(\lambda) 
e^{-(1+\varepsilon)\phi(\lambda)} f(\lambda)\; g(z,\lambda) ,
\]
and $\sum_\lambda |g(z,\lambda)|\leq K$ uniformly in $z$. 
This is so by a duality argument, because 
\[
\{f(\lambda)\}_{\lambda\in\Lambda}\mapsto \omega(z)
e^{-(1+\varepsilon)\phi(z)} f(z) \qquad \text{with } 
f\in \mathcal F^{\infty,0}_{(1+\varepsilon)\phi,\omega}
\]
is a bounded linear functional from a closed subspace of
$\ell_{(1+\varepsilon)\phi,\omega}^{\infty,0}(\Lambda)$  whose norm is bounded
independently of $z$. This is an argument from \cite[p.348--358]{Br} (see also
\cite[p.36]{Se94}).

Consider now $f\in \mathcal F^{p}_{\phi,\tilde{\omega}}\subset \mathcal
F^{\infty,0}_{\phi,\tilde{\omega}}$. 
Given $z\in\C$ take the function 
$P_z$ of Theorem~\ref{weight}, with weight
$\omega/\widetilde{\omega}\in\Omega_\phi$.  
Then $f P_z\in \mathcal
F^{\infty,0}_{(1+\varepsilon)\phi,\omega}$, and by the representation above
\[
\omega(z) e^{-(1+\varepsilon)\phi(z)} f(z)=\sum_{\lambda\in\Lambda} 
\omega(\lambda) 
e^{-(1+\varepsilon)\phi(\lambda)} f(\lambda)\; P_z(\lambda) g(z,\lambda) .
\]
Hence
\begin{eqnarray*}
\omega(z)|f(z)|e^{-\phi(z)}&\lesssim&
\sum_{\lambda\in\Lambda}
\omega(\lambda) |f(\lambda)| e^{-\phi(\lambda)} 
|P_z(\lambda)|
e^{\varepsilon(\phi(z)-\phi(\lambda))}|g(z,\lambda)|\\
&\lesssim& \frac{\omega(z)}{\widetilde{\omega}(z)}\sum_{\lambda\in \Lambda}
\widetilde{\omega}(\lambda) |f(\lambda)| e^{-\phi(\lambda)}
\frac{|g(z,\lambda)|}{1+d_\phi^m(z,\lambda)}.
\end{eqnarray*}

The case $p=\infty$ is clear, so assume that $p<\infty$.
Since $\sum_\lambda |g(z,\lambda)|\leq K$, we  may apply Jensen's inequality
and obtain
\[
\widetilde{\omega}^p(z)|f(z)|^p e^{-p\phi(z)} \rho^{-2}(z)
\lesssim\rho^{-2}(z)
\sum_{\lambda\in \Lambda}
\widetilde{\omega}^p(\lambda) |f(\lambda)|^p e^{-p\phi(\lambda)}
\frac{|g(z,\lambda)|}{1+d_\phi^{mp}(z,\lambda)}.
\]

Now integrate, use that $|g(z,\lambda)|\leq K$ and apply Lemma
\ref{integrals}(b) to finally obtain the sampling inequality
\[
\int_{\C} |f(z)|^p e^{-p\phi(z)}\widetilde{\omega}^p(z)
d\sigma(z)\lesssim
\sum_{\lambda\in\Lambda} \widetilde{\omega}^p(\lambda) |f(\lambda)|^p
e^{-p\phi(\lambda)}.
\]
\end{proof}

\subsection{Nets}
We finish this section by giving useful examples of interpolating  
and sampling sequences.

\begin{lemma}\label{xarxa}
Let $f$ be the multiplier associated to $\phi$, as constructed in the
proof of Theorem~\ref{multiplier}, and let $\Lambda=\mathcal Z(f)$. Then
$\mathcal D_{\Delta\phi}^+(\Lambda)=\mathcal D_{\Delta\phi}^-(\Lambda)=1/2\pi$.
We  say that $\Lambda$ is a \emph{net} associated to $\phi$.
\end{lemma}

\begin{proof} 
The construction of $f$ is made with quasi-squares $R_p$ of $\mu(R_p)=2\pi 
mN$ and $mN$
associated points in a dilated $CR_p$ that made up $\Lambda$. Thus, for
$z\in\C$ and $t>0$:
\begin{eqnarray*}
n(z,t\rho(z))&\geq& mN\#\{p : 
CR_p\subset D^t(z)\}=
\frac 1{2\pi}\mu(\bigcup_{p : CR_p\subset D^t(z)} R_p),\\
n(z,t\rho(z))&\leq& mN\#\{p : CR_p\cap D^t(z)\neq\emptyset\}=
\frac 1{2\pi}\mu(\bigcup_{p : CR_p\cap D^t(z)\neq\emptyset} R_p).
\end{eqnarray*}

By Corollary \ref{flors-quadrats}, 
\[
D^{t-\epsilon(t)}(z)\subset \bigcup_{p : CR_p\subset D^t(z)} 
R_p\subset D^t(z)
\subset
\bigcup_{p : CR_p\cap D^t(z)\neq\emptyset} 
R_p\subset D^{t+\epsilon(t)}(z)  ,
\]
whence
\[
\frac 1{2\pi} \mu(D^{t-\epsilon(t)}(z))\leq n(z,t\rho(z))\leq 
\frac 1{2\pi} \mu(D^{t+\epsilon(t)}(z)) .
\]
The result is then an application of Lemma
\ref{corona-zero}.
\end{proof}

\begin{lemma}\label{nets}
Let $\Lambda$ be a net associated to $\phi$. Then
$\Lambda\in\Int \mathcal F_{(1+\varepsilon)\phi,\omega}^{p}$ 
and $\Lambda\in\Samp \mathcal F_{(1-\varepsilon)\phi,\omega}^{p}$ 
for all $\varepsilon>0$,  $1\leq p\leq \infty$ and $\omega\in \Omega_\phi$.
\end{lemma}

\begin{proof}
By Theorems~\ref{inclusions} and \ref{inclusions-sampling}, 
it is enough to consider the case $\omega=\rho$.

Let $f$ be a multiplier associated to $\phi$
such that $\Lambda=\mathcal Z(f)$.  

Let us start by proving that $\Lambda$ is interpolating. By 
Theorem~\ref{inclusions} it is enough to prove that
$\Lambda\in\Int \mathcal F_{(1+\varepsilon)\phi,\rho}^{p}$ for
all $\varepsilon>0$.
For each $\lambda\in\Lambda$ define
\[
g_\lambda(z)=\frac{f(z)}{z-\lambda}\frac 1{f'(\lambda)} .
\]
We want to see that these functions play a similar role to the 
peak-functions of Lemma~\ref{peak-functions}.
Clearly $g_\lambda(\lambda')=\delta_{\lambda}^{\lambda'}$. 
The growth condition of the multiplier gives
$|f'(\lambda)|\simeq e^{\phi(\lambda)}/\rho(\lambda)$, and then
\[
\rho(z)|g_\lambda(z)| e^{-\phi(z)}\lesssim
\frac{|z-\Lambda|}{|z-\lambda |}
\frac 1{|f'(\lambda)|}\lesssim \rho(\lambda) e^{-\phi(\lambda)} .
\]
Hence
$\|g_\lambda\|_{\calf_{\phi,\rho}^{\infty}}\lesssim 
\rho(\lambda) e^{-\phi(\lambda)}$.

As seen in the proof of Theorem~\ref{inclusions} this is enough to
construct, for any $\varepsilon>0$, an interpolation operator for 
$\mathcal F_{(1+\varepsilon)\phi,\rho}^{p}$. 

Let us see next that
$\Lambda\in\Samp \mathcal F_{(1-\varepsilon)\phi,\omega}^{p}$.
By Theorem \ref{inclusions-sampling} it is enough to consider the
case $p=\infty$ and $\omega=1$, and
by Corollary \ref{weak-limit} it will be enough
to see that every weak limit $(\Lambda^*,(1-\varepsilon)\phi^*,\omega^*)
\in W(\Lambda,(1-\varepsilon)\phi,\omega)$ 
is a uniqueness sequence for $\calf_{(1-\varepsilon)\phi^*,1}^{\infty}$.

Let $(\Lambda_{z_n},\phi_{z_n},\omega_{z_n})
\rightarrow(\Lambda^*,\phi^*,\omega^*)$ and let
$f_{z_n}$ be the corresponding multipliers, with
$\mathcal Z(f_{z_n})=\Lambda_{z_n}$ and
$|f_{z_n}(z)|\simeq e^{\phi_{z_n}(z)}
d_{\phi_{z_n}}(z,\Lambda_{z_n})$.
By Montel's theorem let $\{f_{z_n}\}_n\rightarrow f^*$ with $\mathcal
Z(f^*)=\Lambda^*$ and
$|f^*(z)|\simeq e^{\phi^*(z)} d_{\phi^*}(z,\Lambda^*)$,
i.e, $f^*$ is a multiplier for $\phi^*$.

Consider also a multiplier
$g$ associated to $\varepsilon\phi^*$. In particular
$|g(z)|\simeq e^{\varepsilon\phi^*(z)} d_{\phi^*}(z,\mathcal Z(g))$.
In order to see that $\Lambda^*$ is a uniqueness sequence
assume that $h\in \mathcal F_{(1-\varepsilon)\phi^*,1}^{\infty}$ and 
$h|\Lambda^*=0$. 
Then $hg\in \mathcal F_{\phi^*,1}^{\infty}$, by construction. 
On the other hand, the function
$F:=hg/f^*$
is entire, because $h$ vanishes on $\Lambda^*$.
It is also bounded when $z$ is far from $\Lambda^*$, since $|hg|\lesssim
e^{\phi^*}$ and $|f^*|\gtrsim e^{\phi^*}$. By the maximum principle
$F$ is bounded globally, and 
by Liouville's theorem there exists
$c\in\C$ such that $hg=c f^*$. Since $g$ vanishes in some points outside
$\Lambda^*$ we have $c=0$, hence $h\equiv 0$. 
\end{proof}

\section{Sufficient conditions for sampling}\label{sufisampling}
We prove here the sufficiency part of Theorem A. Assume that 
$\mathcal D_{\Delta\phi}^-(\Lambda)>1/2\pi$. By
Lemma~\ref{separacio} we can assume that $\Lambda$ is $\rho$-separated,
and according to Theorem~\ref{inclusions-sampling} it will be enough to
prove that $\Lambda\in\Finfa$. By Corollary~\ref{weak-limit} this will be done 
as soon as we show that every weak limit $\Lambda^*$ is a uniqueness
sequence for $\mathcal F^{\infty}_{\phi^*,\omega^*}$.

Recall the notation  $n_{\Lambda}(z,r)=\# [\Lambda\cap \overline{D(z,r)}]$.

Assume thus that we have $f\in\calf_{\phi^*,\omega^*}^{\infty}$ with 
$f|\Lambda^*\equiv 0$ and
$\|f\|_{\calf_{\phi^*,\omega^*}^{\infty}}=1$. 
There is no loss of generality in assuming that
$f(0)\neq 0$. 
Applying Jensen's formula to $f$ on $D_{\phi^*}(0)$
\begin{multline*}
\int_0^{r\rho_{\phi^*}(0)} \frac{n_{\Lambda^*}(0,t)}t dt=
\frac 1{2\pi}\int_0^{2\pi} \log|f({r\rho_{\phi^*}(0)} 
e^{i\theta})|d\theta -\log|f(0)| \\
 \leq\frac 1{2\pi}
\int_0^{2\pi} (\phi^*({r\rho_{\phi^*}(0)} e^{i\theta})- 
\log\omega^*(r\rho_{\phi^*}(0) e^{i\theta})) d\theta-\log|f(0)|\\
= \bigl[\frac 1{2\pi}\int_0^{2\pi} 
\phi^*({r\rho_{\phi^*}(0)} e^{i\theta})d\theta -\phi^*(0)\bigr] 
+ \bigl[\log\omega^*(0)-
\frac 1{2\pi}\int_0^{2\pi} \log\omega^*({r\rho_{\phi^*}(0)}e^{i\theta})
d\theta\bigr]\\
+\phi^*(0)-\log\omega^*(0)-\log|f(0)| .\qquad\qquad
\end{multline*}
By definition of flat weight and by Lemma~\ref{distance},
$\omega^*( r\rho_{\phi^*}(0)e^{i\theta}) /\omega^*(0)
\lesssim r^\gamma $ for some $\gamma>0$. Then, Green's
identity yields
\begin{eqnarray*}
\int_0^{r\rho_{\phi^*}(0)}
\frac{n_{\Lambda^*}(0,t)}{t} dt&\leq & \frac 1{2\pi}
\int_{D(0,r\rho_{\phi^*}(0))}\log\frac {r\rho_{\phi^*}(0)}{|\zeta|} 
\Delta\phi^*(\zeta) + \textrm{O}(\log r)\\ 
&=& \frac 1{2\pi}\int_0^{r\rho_{\phi^*}(0)}
\Delta\phi^*(D(0,t))\frac {dt}t +\textrm{O}(\log r) ,
\end{eqnarray*}
for all $r$ big enough.
This contradicts the hypothesis, which
implies in particular that for some $\varepsilon>0$ and all $t$ big enough 
$n_{\Lambda^*}(0,t)\geq (1/2\pi+\varepsilon) \Delta\phi^*(D(0,t))$.

\section{Necessary conditions for sampling}\label{necsampling}
This section contains the proof of the necessity part of Theorem A. 
By Lemma~\ref{separacio}(b) and
Theorem~\ref{inclusions-sampling} it will be enough to prove the following 
result.

\begin{theorem}\label{11}
Let $\Lambda$ be $\rho$-separated. If $\Lambda\in
\Samp \Fdosa$ then
$\mathcal D_{\Delta\phi}^-(\Lambda)\geq 1/2\pi$.
\end{theorem}

We use a result comparing the densities between interpolating and sampling
sequences, as in \cite{RaSt}. We do that by adapting Lemma 4 in \cite{OrSe}
to our setting. 

\begin{lemma}\label{ramanathan}
Let $\varepsilon>0$. Assume $I\in\Int 
\mathcal F_{(1-\varepsilon)\phi,\omega}^{2}$
and $S\in\Samp\Fdosa$ is  $\rho$-separated.  There exists a positive function
$\epsilon(R)$ such that $\lim\limits_{R\to\infty}\epsilon(R)/R=0$
and for every $\epsilon>0$ there is $R_0>0$ with
\[
(1-\epsilon)\; n_I(z,R\rho(z))\leq n_S (z,(R+\epsilon(R))\rho(z))
\qquad z\in\C  .
\]
\end{lemma}

\begin{proof} The proof is as in \cite[Lemma 4]{OrSe} with minor
modifications, so we keep it short.

According to our definition, if $S$ is sampling then 
$\{k(z,s)=K_{\phi,\alpha}(z,s) e^{-\phi(s)}\omega(s)\}_{s\in S}$
is a frame in $\Fdosa$ ($K_{\phi,\omega}$ denotes 
the Bergman kernel, as in Section~\ref{bergman}). That is, for $f\in\Fdosa$
\[
\|f\|_{\Fdosa}^2\simeq\sum_{s\in S}|\left\langle k(z,s),f(z)\right\rangle|^2 .
\]
A consequence is that 
\[
f(z)=\sum_{s\in S}\langle k(\xi,s), f(\xi)\rangle \tilde k(z,s)=
\sum_{s\in S} f(s) e^{-\phi(s)}\omega(s) \tilde k(z,s) ,
\]
where $\tilde k(z,s)$ is the dual frame of $k(z,s)$.

Consider also the  functions $g_i$ given by Lemma \ref{peak-functions} for the
weight $(1-\varepsilon)\phi$. Lemma~\ref{peak-functions}(d) implies that
the normalised functions $\kappa(i,z):=g_i(z)
e^{\phi(i)}/\omega(i)$ form a Riesz basis in the closed linear span $H$
of $\{\kappa(i,z)\}_{i\in I}$ in $\Fdosa$.

Given $z\in\C$ and $R,r>0$  ($R$ much bigger that $r$) consider 
the following two finite
dimensional subspaces of $\Fdosa$:
\begin{eqnarray*}
W_S&=&<\tilde k(\xi,s) : s\in S\cap D^{R+r}(z)>\\
W_I&=&< \kappa(\xi,i) : i\in I\cap D^R(z)  > .
\end{eqnarray*}

Let $P_S$ and $P_I$ denote the orthogonal projections of $\Fdosa$ on $W_S$ and
$W_I$ respectively.  We estimate the trace of the operator
$T=P_IP_S$ in two different ways.
To begin with
\[
\textrm{tr}(T)\leq \textrm{rank}\; W_S\leq\#\{ S\cap D^{R+r}(z) \} .
\]
On the other hand
\[
\textrm{tr}(T)=\sum_{i\in I\cap D^R(z) }
\langle T(\kappa(\xi,i)), P_I \kappa^*(\xi,i)\rangle ,
\]
where $\{\kappa^*(\xi,i)\}$ is the dual basis of $\kappa(\xi,i)$ in $H$. 
Using that $P_I$ and $P_S$ are projections one deduces that
\[
\textrm{tr}(T)\geq \#\{i\in I\cap D^R(z) \}
\bigl(1-\sup_i |\langle P_S(\kappa(\xi,i))-\kappa(\xi,i), 
\kappa^*(\xi,i)\rangle|\bigr) .
\]
Since $\|\kappa(\xi,i)\|_{\Fdosa}\simeq 1$, also
$\|\kappa^*(\xi,i)\|_{\Fdosa}\simeq 1$. Thus we will be done as soon as we
show that $\|P_S(\kappa(\xi,i))-\kappa(\xi,i)\|_{\Fdosa} \leq\varepsilon$ for a
suitable $r\simeq \epsilon(R)$.

We have
\begin{eqnarray*}
\left\|P_S(\kappa(\xi,i))-\kappa(\xi,i)\right\|_{\Fdosa}^2 \lesssim  
\sum\limits_{s\notin D^{R+r}(z) }
|\langle\tilde k(\xi,s),\kappa(\xi,i)\rangle|^2 
=
\sum\limits_{s\notin D^{R+r}(z)}
\bigl|\kappa (s,i)e^{- \phi(s)}\omega(s)\bigr|^2 .
\end{eqnarray*}

Since $S$ is $\rho$-separated, there exists $\eta>0$ such that the disks
$D^\eta(s)$ are pairwise disjoint. Using 
Lemma \ref{estimacio-puntual}(a) we get, for some $c>0$ depending on 
$\phi$ and $\eta$
\begin{eqnarray*}
\left\|P_S(\kappa(\xi,i))-\kappa(\xi,i)\right\|_{\Fdosa}^2  \lesssim 
\int\limits_{\bigcup\limits_{s\notin D^{R+r}(z) } D^\eta(s)} 
|\kappa(\zeta,i)|^2 e^{-2\phi(\zeta)}\omega^2(\zeta)
d\sigma(\zeta) .
\end{eqnarray*}

Applying Lemma \ref{flors} with $r^k=R^\tau$ and $\tau$ so that
$0<(\varepsilon-\tau)\gamma<1$, we see that
there exist $\delta\in(0,1)$, $c>0$ and a function 
$\epsilon(R)=cR^{1-\delta}$ such that
\[
\bigcup\limits_{s\notin D^{R+c \epsilon(R)}(z)} D^\eta(s)\subset  
\C\setminus D^{\varepsilon \epsilon(R)}(i) .
\]
Finally, for $R$ is big enough Lemma \ref{peak-functions}(e) yields
\[
\left\|P_S(k(\xi,i))-\kappa(\xi,i)\right\|_{\Fdosa}^2\lesssim
\int_{\C\setminus D^{\varepsilon \epsilon(R)}(i)}
|\kappa(\xi,i)|^2 e^{-2\phi(\xi)} \omega^2 (\xi) d\sigma(\xi)\lesssim \epsilon .
\]
\end{proof}

\begin{proof}[Proof of Theorem~\ref{11}]
Given $\varepsilon>0$ consider a net $I$ associated to $(1-2\varepsilon)\phi$.
By Lemma \ref{xarxa}  $I\in\Int\mathcal F_{(1-\varepsilon)\phi}^{2,\alpha}$,
and by Lemma~\ref{nets} $\mathcal D_{\Delta\phi}^+(I)=\mathcal
D_{\Delta\phi}^+(I)=(1 -2\varepsilon)/2\pi$. Apply now Lemma \ref{ramanathan}:
there exist $R_0$ and  $\epsilon(R)$ such that for $R>R_0$ 
\begin{eqnarray*}
n_\Lambda(z,R\rho(z))\geq (1-\varepsilon)\; n_I(z,(R-\epsilon(R) )\rho(z)) 
\geq
\frac{(1-\varepsilon)^3}{2\pi}\;  \mu(D^{R-\delta(R)}(z)) ,  
\end{eqnarray*} 
where
$\delta(R)=R-\epsilon(R)-\epsilon(R-\epsilon(R))$. This estimate together with
Lemma \ref{corona-zero} finish the proof. 
\end{proof}

\section{Sufficient conditions for interpolation}\label{sufinterpolation}

Taking into account Theorem~\ref{inclusions} and Lemma~\ref{nets}, in order
to prove the sufficiency part of Theorem B it is enough to prove the
following.

\begin{theorem}\label{equilibrat}
If $\Lambda$ is $\rho$-separated and $\mathcal
D_{\Delta\phi}^+(\Lambda)<1/2\pi$ there exist $\varepsilon>0$ and a sequence 
$\Sigma$ such that  $\Lambda \cup \Sigma$ is a $\rho$-separated net associated
to $(1-\varepsilon)\phi$.
\end{theorem} 

In the proof of this result we need to express the density condition in terms
of the quasi-squares appearing in Theorem~\ref{partition}. this will be done in
Theorem~\ref{density-rectangles}; before we need some preliminaries.

Denote $\phi_r=e^{-r}\phi$.
\begin{lemma}\label{landau}
Let
\[
I_r(\zeta) =\int_{|z-\zeta| < \rho_{\phi_r}(z)/r} \frac {r^2\, dm(z)}
{\pi \rho^2_{\phi_r}(z)}.
\]
Then $\sup\limits_{\zeta\in\mathbb C}|I_r(\zeta) -1|<1/r$. 
\end{lemma}
\begin{proof}
We estimate  $I_r$ using the change of variable 
$w= (z-\zeta)/\rho_{\phi_r}(z)$, whose Jacobian 
is 
\[
\rho^{-2}_{\phi_r}(z) \Bigl| 1- \frac {\langle \nabla \rho_{\phi_r}(z), 
z-\zeta \rangle}{\rho_{\phi_r}(z)} \Bigr| .
\]
From \eqref{lipschitz} it follows
that $|\nabla \rho_{\phi_r}|\leq 1$, hence the Jacobian is bounded
above by $\rho^{-2}_{\phi_r}(z)(1+1/r)$ and below by 
$\rho^{-2}_{\phi_r}(z)(1-1/r)$. Then
\[
1-\frac 1r\int_{|w|\le 1/r} \frac{r^2}{\pi}(1-\frac 1r)\, dm(w) 
\le I_r(\zeta) \le 
\int_{|w|\le 1/r} \frac{r^2}{\pi}(1+\frac 1r)\, dm(w)
=1+\frac 1r .
\]
\end{proof}

It follows immediately from \eqref{cota-inferior} that
there exist $0<\varepsilon<m$ such that 
\[
t^\varepsilon\rho_\phi\lesssim \rho_{\phi/t} \lesssim t^m\rho_\phi .
\]
This implies, with $t=e^r$, 
\begin{equation}\label{exponencial}
\lim_{r\to\infty}\frac {\rho_{\phi_r}(z)}{r\rho_\phi(z)}= \infty
\end{equation}
uniformly in $z\in \mathbb C$.

Let $R_\alpha^s(z)$ denote the rectangle with vertices
$z+s\rho(z)(1+i\alpha)$, $z+s\rho(z)(1-i\alpha)$, $z-s\rho(z)(1+i\alpha)$
 and $z-s\rho(z)(1-i\alpha)$, where
$\alpha\in [e^{-1},e]$ and $e$ is the constant of Theorem~\ref{partition}(b).

\begin{theorem}\label{density-rectangles}
Let $\mu=\Delta\phi$ and let $\nu$ be a positive measure such
that 
\begin{equation}\label{densitat-disc}
\nu(D_\phi^r(z))\leq (1-\varepsilon)\mu(D_\phi^r(z))\quad \forall r\ge r_0,\ 
\forall
z\in \mathbb C.
\end{equation}
There exists $s_0>0$ such that for any $\alpha\in [e^{-1},e]$
\[
\nu( R_\alpha^s(z))\leq (1-\frac\varepsilon{2})
\mu( R_\alpha^s(z))\quad \forall s\ge s_0,\ \forall
z\in \mathbb C.
\]
\end{theorem}
\begin{proof}
Fix  $r$ big enough so that $\rho_{\phi_r}/r >r_0\rho_\phi$
and $(1+1/r)(1-\varepsilon)<(1-1/r)(1-3\varepsilon/4)$. This can
be done because of \eqref{exponencial}. By hypothesis
\[
\nu(D_{\phi_r}^{1/r}(z)) \le (1-\varepsilon)\; 
\mu(D_{\phi_r}^{1/r}(z))\quad  \forall
z\in \mathbb C, 
\]
and if $s$ is much bigger than $r$ we get
\begin{equation*}
\int_{z\in R_s^\alpha(w)} \frac{r^2}{\pi\rho^{2}_{\phi_r}(z)}
\nu(D_{\phi_r}^{1/r}(z))
\, dm(z)\le  (1-\varepsilon)
\int_{z\in R_s^\alpha(w)} \frac{r^2}{\pi\rho^{2}_{\phi_r}(z)}
\mu(D_{\phi_r}^{1/r}(z))
\, dm(z).
\end{equation*}
Denote  
\begin{eqnarray*}
\Omega_r(\zeta)&=&\{z\in \mathbb C,\ |z-\zeta|< \rho_{\phi_r}(z)/r\}\\
F_{r}(w,s)&=&\{\zeta \in \mathbb C,\ \Omega_r(\zeta)\subset R_\alpha^s(w)\} \\
G_{r}(w,s)&=&\bigcup_{\zeta\in R_\alpha^s(w)} \Omega_r(\zeta).
\end{eqnarray*}
Reversing the order of integration and using the previous Lemma we deduce that
\[
\nu( F_{r}(w,s)) \leq (1-\frac 34 \varepsilon)\; \mu( G_{r}(w,s)) . 
\]  
It is clear that $F_r(w,s)\subset R_\alpha^s(w)\subset G_r(w,s)$.  Similarly
to the proof of Lemma \ref{flors}, there exists $\epsilon(s)$  such that
$R^{s-\epsilon(s)}_\alpha(w)\subset F_r(w,s)$ and $ G_r(w,s)\subset
R^{s+\epsilon(s)}_\alpha(w)$.

By Remark~\ref{corona-rectangles}
\[
\lim_{s\to\infty}
\frac{\mu (R^{s+\epsilon(s)}_\alpha(w))}
{\mu(R^{s-\epsilon(s)}_\alpha(w))}=1
\]
uniformly in $z$, and therefore there exists $s_0$ such that for $s>s_0$
\begin{eqnarray*}
\nu(R^{s-\epsilon(s)}_\alpha(w))&\leq& (1-\frac 34 \varepsilon)\;
\mu(G_{r}(w,s) )\leq (1-\frac 34 \varepsilon)\;
\mu(R^{s+\epsilon(s)}_\alpha(w))\\
&\leq& (1-\frac \varepsilon 2 )\; 
\mu(R^{s-\epsilon(s)}_\alpha(w)) .
\end{eqnarray*} 
\end{proof}

\begin{proof}[Proof of Theorem \ref{equilibrat}]
Take an entire function $g$ vanishing exactly on $\Lambda$. 
We will construct a sequence $\Sigma$ and an entire function $h$ such that
for some $\varepsilon>0$,
\begin{itemize}
\item[(i)] $\Lambda\cup \Sigma$ is
$\rho$-separated.
\item[(ii)] $h$ vanishes exactly on $ \Sigma$.
\item[(iii)] For any $\tau>0$,  $|\log|h(z)|- (1-\varepsilon)\phi(z) + 
\log |g(z)|| \le C_\tau$ 
if $D^{\tau}(z)\cap ( \Lambda\cup\Sigma)=\emptyset$. 
\end{itemize}

Accepting this we reach the result by taking $f=g h$. This is so because  the
separateness of $\Lambda\cup\Sigma$ and (iii) imply that $f$ is a multiplier
for $(1-\varepsilon)\phi$.  
\end{proof}  

\begin{proof}[Construction of $\Sigma$ and $h$] 
To avoid the repetition of the factors $2\pi$ and $1-\varepsilon$, 
denote here $\mu=(1-\varepsilon)\Delta\phi/2\pi$.
Let  
\[
\tilde\mu=\mu-\sum_{\lambda\in\Lambda}\delta_\lambda=
\frac 1{2\pi} \Delta\bigl( (1-\varepsilon)\phi -\log |g|\bigr) .
\]
Following Theorem
\ref{partition} and the Remark thereafter, given $n,M\in\N$ we can take a
system of  quasi-squares $\{R_k\}_k$ such that, denoting  $\mu_k=\mu_{|R_k}$,
we have $\mu=\sum_k\mu_k$ and  $\mu_k(\C)=Mn$.  Then
$\tilde\mu=\sum_k\tilde\mu_k$, being
\[
\tilde\mu_k=\mu_k-\sum\limits_{\lambda\in\Lambda\cap R_k}\delta_\lambda .
\]
By hypothesis there exists $\varepsilon>0$ such that
$\mathcal D_{\Delta\phi}(\Lambda)<1/2\pi -4\varepsilon$. Therefore, 
there exists $r_0>0$
such that
\[
\tilde\mu (D^r(z))\geq 3\varepsilon \mu (D^r(z))\qquad 
\textrm{for all}\quad z\in\C, r\geq r_0 .
\]
Also, Theorem~\ref{density-rectangles} implies that for 
$M\geq m/(2\varepsilon) $ and $n$ big enough:
\[
Mn\geq \tilde\mu(R_k)\geq 2\varepsilon\mu(R_k) =
2{\varepsilon} Mn\geq  mn .
\]
Let $\tilde\mu(R_k)=m_k n$, with $m\leq m_k\leq M$. Notice that $m_k\in\N$,
since $\mu(R_k)\in\N$. Applying Lemma~\ref{sigma} we obtain a sequence 
$\Sigma$ made of points $\sigma_1^k,\dots,\sigma_{m_k n}^k\in CR_k$ so that the
first $m$ moments of the
measures $\nu_k=\tilde\mu_k-\sum\limits_{j=1}^{m_kn} \delta_{\sigma_j^k}$ 
vanish. Furthermore, it is clear that we can choose the 
$\tau_j^k$ in the proof of Lemma~\ref{sigma} so that $\Lambda\cup\Sigma$ is
$\rho$-separated.

Let 
\[
\nu=\sum_k\nu_k=\frac 1{2\pi} \Delta ( (1-\varepsilon)\phi -\log |g| )
-\sum_{\sigma\in\Sigma} \delta_\sigma .
\] 
In order to  prove (iii) consider $v= (1-\varepsilon)\phi -\log |g|-w$, where
\[
w(z) = \int_{\C} \log|z-\zeta| \, d\nu(\zeta) .
\]
Since
\[
\Delta v=2\pi \sum_{\sigma\in\Sigma} \delta_\sigma,
\]
there exists $h$ entire (vanishing exactly on $\Sigma$) such that $\log|h|=v$.

We need to estimate $|w(z)|$ when
$|z-\Lambda\cup\Sigma|\geq \tau \rho(z)$. Given $z\in\C$, let $k_0\in\N$ be
such that $z\in R_{k_0}$. By Theorem~\ref{partition}(c), there exists $r_0>0$
independent of $z$ such that  $R_{k_0}\subset D^{r_0}(z)\subset C R_{k_0}$. We
have
\[
 w(z)=\iC \log|z-\zeta| d\nu(\zeta)=\iC \log|z-\zeta| d\nu_{k_0}(\zeta)+
\sum_{k:k\neq k_0} \iC\log|z-\zeta| d\nu_k(\zeta) ,
\]
and we estimate the two terms separately.

Let $C>0$ be the constant of Lemma \ref{sigma}.
Since the first $m$ moments of $\nu_{k_0}$ vanish, 
\begin{eqnarray*}
\Bigl |\iC \log|z-\zeta| d\nu_{k_0}(\zeta)\Bigr |
&=&\Bigl |\iC \log\frac{|z-\zeta|}{r_0 \rho (z)} d\nu_{k_0}(\zeta)\Bigr |
\lesssim\Bigl|\int_{C R_{k_0}} \log\frac{r_0\rho(z)}{|z-\zeta|}
d\mu(\zeta)\Bigr| + K |\log\tau| \\
&\lesssim&
\int_{D^{c r_0}(z)} \log\frac{c r_0\rho(z)}{|z-\zeta|} 
d\mu(\zeta) + K|\log \tau|\le C_\tau.
\end{eqnarray*}
 
The other integral is estimated using the moment
condition for each $\nu_k$, as in the estimate of $I_1$ in
Theorem~\ref{regularitzacio}. 
\end{proof}

\section{Necessary conditions for interpolation}\label{neceinterpolation}

Let us start by proving the non-strict density inequality. By
Theorem~\ref{inclusions}, it is enough to consider the case $p=2$.

\begin{theorem}\label{necessary-dos}
If $\Lambda\in\Int  \Fdosa$ then 
$\mathcal D^+_{\Delta\phi}(\Lambda)\leq 1/2\pi$.
\end{theorem}
\begin{proof}
Given $\varepsilon>0$, take a net $S$ associated to $(1+2\varepsilon)\phi$,
as described in Lemma \ref{xarxa}. Lemma \ref{nets} implies that 
$S\in\Samp \calf_{(1+\varepsilon)\phi,\omega}^{2}$, and by Lemma
\ref{ramanathan}, there exists $R_0>0$  such that
\[
n_\Lambda(z,R\rho(z))\leq (1+\varepsilon)\;
n_S(z,(R+\epsilon(R))\rho(z))\qquad z\in\C\; ,\; R\geq R_0 .
\]
Since $S$ is a net of density $(1+2\varepsilon)/2\pi$, the radius
$R_0$ can be taken so that for $R\geq R_0$
\[ 
n_S(z,(R+\epsilon(R))\rho(z))\leq \frac{1+3\varepsilon}{2\pi}\; 
\mu(D^{R+\epsilon(R)}(z))  .
\]
This and Corollary \ref{flors-quadrats} give the result.
\end{proof}

Let us see now that the inequality is strict.

\begin{proof}[Proof of the necessity part in Theorem B]
Assume that $\Lambda\in\Int \Fpa$. We know that  $\mathcal D_{\Delta\phi}^+
(\Lambda)\leq 1/2\pi$. 
In order to see that
$\mathcal D_{\Delta\phi}^+(\Lambda)< 1/2\pi$ consider, given $\varepsilon>0$,
a net $\Sigma$ associated to  $2\varepsilon\phi$
such that $Z:=\Lambda\cup\Sigma$ is
$\rho$-separated. 

\begin{lemma}\label{fk}
Denote $Z=\{z_k\}_k$. 
For every $m\in\N$ and $\varepsilon>0$ there exist $C>0$ and functions 
$f_k\in \Finfa$ such that 
\begin{itemize} 
\item[(a)] $f_k(z_k)=1$.
\item[(b)] $f_k(z_j)=0$ for all $z_j\in D^{1/\varepsilon}(z_k )$.
\item[(c)] $|f_k(z)|\leq C M(\Lambda) e^{\phi(z)-\phi(z_k)}
\displaystyle\frac{\omega(z_k)}{\omega(z)}\frac{1}
{1+d_\phi^m(z,z_k)}$.
\item[(d)] 
$\|f_k\|_{\phi,\infty}\leq C M(\Lambda) e^{-\phi(z_k)}\omega(z_k)$.
\end{itemize}
\end{lemma}

\begin{proof}
Assume first that $z_k=\lambda_k\in\Lambda$. By hypothesis there exists
$g_k\in\Fpa\subset\Finfa$ with  $g_k(\lambda_k)=1$, $g(\lambda_j)=0$, and 
$\|g_k\|_{\Finfa}\leq M(\Lambda) e^{-\phi(\lambda_k)}\omega(\lambda_k)$. 
Since $\Lambda$ plus a finite number of points is still in
$\Int \Fpa$ (Lemma \ref{add-point}), we can take $g_k$ so that 
moreover $g_k(\sigma_j)=0$ if $|\lambda_k-\sigma_j|\leq  
1/\varepsilon \rho(\lambda_k)$ and $g_k(c_j)=0$, $j=1,\dots,M$, where
$c_j=\lambda_k+2\delta\rho(\lambda_k) e^{j\frac{2\pi i}M}$ and $\delta>0$
is taken so that the balls $\{B(\lambda,10\delta)\}_\lambda$ are pairwise 
disjoint. 

By construction of the nets there exists $C$ independent of $z$ and
$\varepsilon$ such that $\# \Sigma\cap D^{1/\varepsilon}(z)\leq C$ for any 
$\Sigma$ net of density  $\varepsilon/\pi$.

Define then
\[
f_k(z)=(2\delta)^{-M}\frac{g_k(z)}{(z-c_1)\cdots (z-c_M)} \rho^M(z_k) .
\]
It is clear that $f_k\in\Finfa$ satisfies (a) and (b). 

For $z\notin \cup_{j=1}^M D^\delta(c_j)$, 
\[
|f_k(z)|\leq C |g_k(z)|\bigl(\frac{\rho(z_k)}{|z-z_k|}\bigr)^M\leq C M(\Lambda)
e^{\phi(z)-\phi(z_k)} \frac{\omega(z_k)}{\omega(z)}
\bigl(\frac{\rho(z_k)}{|z-z_k|}\bigr)^M ,
\]
and the estimate follows from Lemma~\ref{distance}.

For $z\in D^\delta(c_j)$ we have
\[
|f_k(z)|\leq C  \bigl| \frac{g_k(z)}{z-c_j} \bigr| \rho(z_k) .
\]
Estimating like in (iii) in the proof of  Theorem~\ref{weight} we get 
$|f_k(z)|\leq C M(\Lambda) e^{\phi(z)-\phi(z_n)}$,
as desired.

In case $z_k=\sigma_k\in\Sigma$, use again that $\Lambda$ plus one point is
$\Fdosa$-interpolating and start with  $g_k\in\Fdosa\subset\Finfa$ such that 
$g_k(\sigma_k)=1$,
$g_k(\lambda_j)=0$ for all $j$.  Then proceed as before.
\end{proof}

\begin{lemma}
$Z\in\Int  \Finfa$.
\end{lemma}

\begin{proof}
Given $v=\{v_k\}_k\in\ell_{\phi}^{\infty,\alpha}(Z)$ consider the 
pseudo-extension
\[
E(v)(z)=\sum_{k=1}^\infty v_k f_k(z)  .
\]
Let us see first that $E(v)\in\Finfa$. 
By (c) above and Lemma \ref{suma} we see that for any $z\in\C$
\begin{eqnarray*}
\omega (z) e^{-\phi(z)}|E(v)(z)| &\lesssim& 
\sum_{k=1}^\infty \omega(z_k) |v_k| e^{-\phi(z_k)}\frac{1}
{1+d_\phi^M(z,z_k)}\lesssim
\|v\|_{\ell_{\phi}^{\infty,\alpha}(Z)} .
\end{eqnarray*}

Let $R$ denote the restriction operator from $\Finfa$ 
to $\ell_{\phi,\omega}^{\infty}(Z)$. In order
to see that $Z$ is in $\Int  \mathcal F_{(1+\varepsilon)\phi,\omega}^{\infty}$ 
it will
be enough to prove that $\|RE-I\|_{op}<1$, since then
$(RE)^{-1}=I+\sum_{k=1}^\infty
(RE-I)^k$ converges and $E(RE)^{-1}$ defines an inverse to $R$.

By Lemma \ref{fk}(b) and (c)
\begin{multline*}
\|RE(v)-v\|_{\ell_{\phi,\omega}^{\infty}(Z)}=
\Bigl\|\bigl\{\sum_{k:k\neq j} v_k f_k(z_j) 
\bigr\}_{j\in\N}\Bigr\|_{\ell_{\phi,\omega}^{\infty}(Z)}\\
\leq \sup_{j\in\N} \omega(z_j) e^{-\phi(z_j)}
\sum_{k: z_j\notin D^{1/\varepsilon}(z_k)}|v_k| |f_k(z_j)|
\leq C M(\Lambda) \|v\|_{\ell_{\phi,\omega}^{\infty}(Z)}
\sum_{k: z_j\notin D^{1/\varepsilon}(z_k)} 
\frac 1{d_\phi^m(z_j,z_k)} 
\end{multline*}
By Lemma~\ref{distance} and Corollary \ref{suma-zero}, if $m$ is big 
and $\varepsilon$
is small enough we have 
\[
\|RE(v)-v\|_{\ell_{\phi,\omega}^{\infty}(Z)}\leq 1/2\;  
\|v\|_{\ell_{\phi,\omega}^{\infty }(Z)} ,
\]
thus $\|RE-I\|_{op}<1/2$,
as desired.
\end{proof}

By this Lemma and the results above we have 
$\mathcal D_{\Delta\phi}^+(Z)\leq 1/2\pi$, i.e
for all $\delta>0$ there exists $R_0$ such that for all $z\in\C$ and $R>R_0$
\[
n_\Lambda(z,R\rho(z))+n_\Sigma(z,R\rho(z))\leq (1/2\pi+\delta) 
\mu(D^R(z)) .
\]
By Lemma \ref{xarxa}, $\mathcal D_{\Delta\phi}^-(\Sigma)=\varepsilon/\pi$,
thus for all $\delta>0$
there exists $R_0$ such that for all $ R>R_0$
\[
n_\Sigma(z,R\rho(z))\geq (\varepsilon/\pi-\delta) 
\mu(D^R(z))\qquad  z\in\C\   .
\]
This shows that for $\delta>0$ and $R$ big enough
\[
n_\Lambda(z,R\rho(z))\leq (\frac{1-2\varepsilon}{2\pi}+2\delta)\;  
\mu(D^R(z))\qquad z\in\C  ,
\]
hence $\mathcal D_{\Delta\phi}^+(\Lambda)<1/2\pi$.
\end{proof}

\section*{Appendix. Alternative construction of
peak functions.}
As seen at the end of the proof of Theorem~\ref{peak-functions}, it is enough
to consider the case $\omega=\rho$. Also, 
it will be enough to prove 
that for any $\phi$ there exist $C,\delta>0$ such that for all $\eta\in\C$
there is $P_\eta$ holomorphic with 
$P_\eta(\eta)=1$ and
\[
|P_\eta(z)|\leq C
e^{\phi(z)-\phi(\eta)}
\min\bigr\{1,\bigl(\frac{\rho(\eta)}{|z-\eta|} \bigr)^\delta\bigr\} ,
\]
since then we can apply this to $\varepsilon\delta /m\; \phi(z)$, take the 
$m$-th power and use Lemma~\ref{distance} to conclude.

We claim that there exists $h_\eta$ holomorphic with 
$h_\eta(\eta)=0$, $h_\eta^\prime(\eta)=1$ and
$|h_\eta(z)|\lesssim e^{\phi(z)-\phi(\eta)}
\rho^2(\eta)/ \rho(z)$.

Once this is proved we take 
$w_\eta(z)=h_\eta(z)/(z-\eta)$ and use Lemma \ref{christ} 
to deduce that
\[
|P_\eta(z)|\lesssim e^{\phi(z)-\phi(\eta)}
\frac{\rho(\eta)}{|z-\eta|}
\bigl(\frac{|z-\eta|}{\rho(\eta)}\bigr)^{1-\delta}=
e^{\phi(z)-\phi(\eta)}
\bigl(\frac{\rho(\eta)}{|z-\eta|}\bigr)^{\delta}\qquad z\notin D(\eta) .
\]

In order to construct the function $h_\eta$ define first
\[
F(z)=(z-\eta)\mathcal X\bigl(\frac{|z-\eta|^2}{\rho^2(\eta)} \bigr)
e^{H_\eta(z)} ,
\]
where $H_\eta$ is a holomorphic function such that $\Re 
H_\eta =h_\eta$ (see Lemma \ref{cota-disc}) and $\mathcal X$
is a smooth cut-off function 
with 
$\mathcal X\equiv 1$ for $|\zeta|<1$, 
$\mathcal X\equiv 0$ for $|\zeta|\geq 2$ and $|\mathcal X'|$ bounded.

Notice that by construction and by Lemma \ref{cota-disc}, we have
\[
\rho(z)|F(z)|e^{-\phi(z)}\lesssim \rho^2(\eta)
e^{-\phi(\eta)} .
\]

\begin{lemma} There exists $u$ solution to $\db u=\db F$ such that
$u(\eta)=\partial u(\eta)=0$ and 
$\|u \|_{\mathcal F_{\phi,\rho}^\infty} \leq C\rho^2(\eta) e^{-\phi(\eta)} $
\end{lemma}

Once this is proved we take $h_\eta=F-u$ and we are done.

\begin{proof}
First we show that there exists a solution $u$ as in the statement but
satisfying an analogous $L^2$ estimate instead of the $L^\infty$ one.
We use H\"ormander's theorem \cite{Ho}: for every $\psi$ subharmonic in $\C$
there exists a solution $u$ to $\db u=\db F$ such that
\[
\int_{\C} |u|^2 e^{-2\psi}\leq C \int_{\C} |\db F|^2 
\frac{e^{-2\psi}}{\Delta\psi} .
\]

Define $\psi= \phi+2 v$, where
\[
v(z)=\log|z-\eta|-\frac 1{\mu(D^s(\eta))}\int_{D^s(\eta)}\log |z-\zeta|
\Delta\phi(\zeta) dm(\zeta) .
\]
Take $s$ so that $\mu(D^s(\eta))=8\pi$. By the  doubling condition 
there exists $c$ depending only on the doubling
constant $C_{\Delta\phi}$ such that $s\leq c$. Then
\[
\Delta\psi \geq \Delta\phi 
-\frac{4\pi}{\mu(D^s(\eta))} \Delta\phi = \frac 12 \Delta\phi\simeq
\rho^{-2} .
\]

By construction $v$ is bounded above. Notice also that there exists $C>0$
(independent of $\eta$) such that $-v(z)\leq C$ for all $z\in\textrm{supp}(\db
F)$. Since $|\db F|\lesssim e^{h_\eta}$, we deduce from H\"ormander's estimate
and Lemma~\ref{cota-disc} that
\[
\|u\|_{\mathcal F_{\phi,\rho}^2}\leq \iC |u|^2 e^{-2\psi} 
\leq C \int_{D^2(\eta)\setminus D(\eta)} e^{2h_\eta} e^{-2\psi} \rho^2
\lesssim \rho^4(\eta) e^{-2\phi(\eta)}.
\]
On the other hand 
\[
e^{-2\psi(z)}\simeq |z-\eta|^{-4}\qquad\textrm{for}\quad |z-\eta|\leq
\epsilon\rho(\eta) ,
\]
thus necessarily $u(\eta)=\partial u(\eta)=0$.

Let us see now that $u$ satisfies also the $L^\infty$ estimate.
For any $z\in supp(\db F)$ define
\[
U(\zeta)=\frac{K\rho(z)}{\rho^2(\eta) e^{-\phi(\eta)} } u(\zeta) ,
\]
where $K>0$ will be chosen later on. Then
\[
\int_{D(z)} |U(\zeta)|^2 e^{-2\phi(\zeta)}\leq
\frac{\rho^2(z)}{\rho^4(\eta) e^{-2\phi(\eta)}}
\| u \|_{L^2(e^{-\phi})}^2\lesssim\rho^2(z) .
\]

Also, since $\rho(\zeta)\simeq\rho(\eta)$ on $supp(\db F)$, we have
\begin{equation*}
\rho(z)\sup_{\zeta\in D(z)}  |\db U(\zeta)| e^{-\phi(\zeta)}=
\sup_{\zeta\in D(z)}\frac{K\rho^2(z)}{\rho^2(\eta) e^{-\phi(\eta)} }
|\db F(\zeta)| e^{-\phi(\zeta)}\lesssim 1 .
\end{equation*}

We choose $K$ (independent of $z$) so that
\begin{itemize}
\item[(a)] $\displaystyle \frac 1{\rho^2(z)}
\displaystyle\int_{D(z)} |U(\zeta)|^2 e^{-2\phi(\zeta)}\leq 1$ ,
\item[(b)] $\rho(z)\sup\limits_{\zeta\in D(z)} |\db U (\zeta)| 
e^{-\phi(\zeta)}\leq 1$.
\end{itemize}
We will be done as soon as we prove that 
\[
|U(z)|e^{-\phi(z)}\leq C .
\]
This is consequence of \cite[Lemma 3.1]{Be} applied to
the function
$V(\zeta)= u(\rho(z)\zeta+z) $. Defining 
$\phi_z(\zeta)=\phi(\rho(z)\zeta+z)$ and changing to the
variable
$w=\rho(z)\zeta+z$ we see that
\[
\int_{\D} |V(\zeta)|^2 e^{-2\phi_z(\zeta)} dm(\zeta)=
\int_{D(z)} |U(w)|^2 e^{-2\phi(w)}
\frac{dm(w)}{\rho^2(z)}\leq 1
\]
and
\[
\sup_{\zeta\in\D} |\db V(\zeta)|^2 e^{-2\phi_z(\zeta)}=
\sup_{w\in D(\eta)}|\db U(w)|^2 e^{-2\phi(w)} \rho(z)\leq 1 .  
\]
Thus, by \cite[Lemma 3.1]{Be}
$|V(0)|^2 e^{-\phi_z(0)}\leq C e^{-a_{\phi_z}}$,
where 
\[
a_{\phi_z}=\sup \{ \psi(0): \psi\leq 0\; ,  \Delta \psi=\Delta\phi_\eta\}.
\]
Defining $v$ so that $\psi(z)=v(\rho(z)\zeta+z)$ we see that
\[
a_{\phi_z}=\sup\{ v(z) : v\leq 0 : \Delta v=\Delta\phi\}.
\]
The function $v(w)=\phi(w)-h_z(w)-\phi(z)-A$ is negative
in $D(z)$ if $A$ is big enough (Lemma \ref{cota-disc}) and 
$v(z)=-A$. Hence
$a_{\phi_z}\geq -A$ and
$|U(z)|^2 e^{-2\phi(z)}=|V(0)|^2 
e^{-2\phi_z(0)}\leq Ce^A$, as desired.
\end{proof}
%
%
\providecommand{\bysame}{\leavevmode\hbox to3em{\hrulefill}\thinspace}
\renewcommand{\MR}[1]{MR #1}

\end{document}